\documentclass[10pt] {amsart}

\usepackage{amssymb}
\usepackage{pb-diagram}
\usepackage{Bdiagx}
\usepackage{graphicx}

\numberwithin{equation}{section}
\newcommand{\TryPackage}[3]{\IfFileExists{#1.sty}{\usepackage{#1}#2}{#3}
}
\TryPackage{mathrsfs}{\renewcommand{\mathcal}{\mathscr}}{%
      \TryPackage{eucal}{}{}}

\newcommand{\lto}{\longrightarrow}
\newcommand{\al}{\alpha}
\newcommand{\be}{\beta}
\newcommand{\de}{\delta}
\newcommand{\ga}{\gamma}
\newcommand{\ep}{\epsilon}
\newcommand{\varep}{\varepsilon}
\renewcommand{\th}{\theta}
\newcommand{\la}{\lambda}
\newcommand{\om}{\omega}
\newcommand{\si}{\sigma}
\newcommand{\De}{\Delta}
\newcommand{\Ga}{\Gamma}
\newcommand{\Om}{\Omega}
\newcommand{\La}{\Lambda}
\newcommand{\Si}{\Sigma}
\newcommand{\Th}{\Theta}
\newcommand{\ba}{{\boldsymbol{a}}}
\newcommand{\bb}{{\boldsymbol{b}}}
\newcommand{\bc}{{\boldsymbol{c}}}
\newcommand{\ZZ}{{\mathbb Z}}
\newcommand{\RR}{{\mathbb R}}
\newcommand{\CC}{{\mathbb C}}

\newcommand{\cA}{{\mathcal A}}
\newcommand{\cB}{{\mathcal B}}\newcommand{\cC}{{\mathcal C}}
\newcommand{\cG}{{\mathcal G}}
\newcommand{\cM}{{\mathcal M}}

\newcommand{\cP}{{\mathcal P}}
\newcommand{\cQ}{{\mathcal Q}}
\newcommand{\cU}{{\mathcal U}}
\newcommand{\cZ}{{\mathcal Z}}
\newcommand{\red}{{\hbox{\scriptsize \sl red}}}
\newcommand{\ab}{{\hbox{\scriptsize \sl ab}}}

\newcommand{\im}{\operatorname{im}}
\newcommand{\Image}{\operatorname{Image}}
\newcommand{\hol}{\operatorname{{\it hol}}}
\newcommand{\tr}{\operatorname{\it tr}}
\newcommand{\Hess}{\operatorname{Hess}}

\newcommand{\Spec}{\operatorname{Spec}}
\newcommand{\End}{\operatorname{End}}
\newcommand{\Hom}{\operatorname{Hom}}

\newcommand{\Mas}{\operatorname{Mas}}
\newcommand{\Xmin}{X_\ell}
\newcommand{\Xmax}{X_u}
\newcommand{\Ymin}{Y_\ell}
\newcommand{\Ymax}{Y_u}
\newcommand{\Zmin}{Z_\ell}
\newcommand{\Zmax}{Z_u}

\def\mapright#1{\smash{
\mathop{\longrightarrow}\limits^{#1}}}
\def\del{\partial}

\def\del{\partial}

\newtheorem{thm}{Theorem}[section]
\newtheorem{lem}[thm]{Lemma}
\newtheorem{cor}[thm]{Corollary}
\newtheorem{prop}[thm]{Proposition}

\theoremstyle{definition}
\newtheorem{defn}[thm]{Definition}

\theoremstyle{remark}
\newtheorem{remark}[thm]{Remark}

\newtheorem{claim}[thm]{Claim}

\begin{document}
\title{The Integer Valued SU(3) Casson Invariant for Brieskorn
spheres}
\author{Hans U. Boden}
\address{Department of Math \& Stats, McMaster University,
Hamilton, Ontario L8S 4K1 Canada}
\email{boden@mcmaster.ca}

\author{ Christopher M. Herald}
\address{Department of Mathematics, University of Nevada, Reno Nevada
89557}
\email{herald@unr.edu}

\author{Paul A. Kirk}
\address{Department of Mathematics, Indiana University, Bloomington,
Indiana
47405}
\email{pkirk@indiana.edu}
\date{Sept. 1, 2003}
\thanks{The first author was partially supported by a grant from the Natural Sciences and Engineering Research Council of Canada. 
The second author acknowledges support from NSF under
grant no.~DMS-0072564.
The third author would like to thank the NSF for its support under
grant no.~DMS-0202148}
\begin{abstract}
We develop techniques for computing
the integer valued $SU(3)$ Casson invariant
defined in \cite{BHK}.
Our method involves
  resolving the singularities in the flat moduli space
using a {\it twisting} perturbation and
analyzing its effect on the topology of  the perturbed flat moduli
space.
These techniques, together with Bott-Morse theory and the splitting principle for
  spectral flow, are applied to calculate $\tau_{SU(3)}(\Sigma)$
for all Brieskorn homology spheres.
\end{abstract}
\maketitle
\section{Introduction}

In this article we compute the integer valued $SU(3)$ Casson
invariant $\tau_{SU(3)}$
for Brieskorn spheres $\Si(p,q,r)$.  Computations of
$\tau_{SU(3)} (\Si(2,q,r))$  appear in \cite{BHK},
and  we extend those computations to all
Brieskorn spheres.

If $\Si$ is a 3-dimensional homology sphere whose
flat $SU(3)$ moduli space
is nondegenerate and  0-dimensional, then
the integer valued $SU(3)$ Casson invariant $\tau_{SU(3)}(\Si)$
is simply a signed count of the points in the irreducible stratum
of the flat moduli space.
On the other hand, if the
moduli space has positive dimension and
is nondegenerate in the sense of Bott and Morse
(or more generally if its lift to the based
moduli space is nondegenerate),
then one can apply standard (equivariant) Morse theoretic techniques to
compute the invariant $\tau_{SU(3)}(\Si)$ (see \cite{BH2}).

The family of computations given here represents the
first successful attempt to
compute the invariant $\tau_{SU(3)}(\Si)$ for manifolds $\Si$ with truly singular
moduli spaces.  Even in the connected sum theorem of \cite{BH2}
where one encounters components of  mixed type in the moduli space
(i.e., components containing both irreducible
and reducible gauge orbits), when lifted to the based moduli space,
these components become nondegenerate and one can apply
equivariant Bott-Morse theory to determine the invariant
$\tau_{SU(3)}$.
In contrast, the flat $SU(3)$ moduli space of the  Brieskorn spheres
considered in this paper are singular even when lifted to the based
moduli space.
Thus  the  perturbation
techniques presented here go beyond the standard theory and in fact
provide a new approach to transversality issues that may well
apply more generally.

The new approach involves a combination of manifold
decomposition and Mayer-Vietoris techniques and traditional holonomy
perturbations. Simply put, our idea is to construct
a special type of perturbation (called the
{\it twisting perturbation}) and analyze its effect
on the moduli space.
We prove that under such perturbations, the moduli space
becomes nondegenerate and we express the invariant
$\tau_{SU(3)}$ in terms of the topology of
the perturbed moduli space
and the spectral flow of the odd signature operator.
 
The remainder of this paper is divided into five sections.
Section \ref{s2} presents a
detailed description of the $SU(3)$ representation varieties
of Brieskorn spheres. Corresponding results for
knot complements are given in Section \ref{s3}.
Section \ref{s4} introduces the twisting perturbations and
describes their effect on the moduli spaces.
Section \ref{s5} presents spectral flow computations based
on a splitting argument, and
Section \ref{s6} presents a lattice point count which provides
numerical calculations of
$\tau_{SU(3)}$ for families of  Brieskorn spheres $\Si(p,q,r)$,
including all homology 3-spheres
obtained by Dehn surgery on a $(p,q)$ torus knot.
The rest of the introduction is devoted to
outlining  the main argument.

Recall first that if $\pi$ is a (finitely presented) group,  a
representation $\al\colon \pi\to SU(3)$ is {\em irreducible} if no
nontrivial linear subspace of $ \CC^3$ is invariant under $\al(g)$
for all $g\in \pi$.  This is equivalent to the condition
that the stabilizer    of $\al$  with respect to the conjugation action
equals the center of $SU(3)$.   Otherwise, $\al$ is reducible and
its image can be conjugated to lie in the subgroup
$S(U(2)\times U(1))$.

Suppose that $\Si$ is a homology 3-sphere and let
$R(\Si, SU(3))$ be the set of  conjugacy classes of
  representations $\al\colon \pi_1(\Si) \to SU(3)$. Then
  $R(\Si, SU(3))$ is a real algebraic variety
homeomorphic to the
moduli space $\cM(\Si)$  of flat $SU(3)$ connections on
$\Si$.  We denote by $R^*(\Si, SU(3))$ the subspace of
conjugacy classes of irreducible
representations and by $\cM^*(\Si)$ the subspace of
irreducible flat connections.

The integer valued $SU(3)$ Casson invariant $\tau_{SU(3)}(\Si)$
is defined in \cite{BHK} and gives an algebraic count of the
conjugacy classes of  irreducible
representations of $\pi_1(\Si)$, with a correction term involving the
reducible representations. More precisely, the flatness equations are
perturbed so that the flat moduli space becomes nondegenerate, and
gauge orbits of irreducible perturbed flat connections are counted
with sign given by
the spectral flow of the $su(3)$
odd signature operator. The resulting integer depends
on the perturbation used, and to compensate for this we add a
correction term defined in terms of the reducible stratum.

For $\Si=\Si(p,q,r)$ the Brieskorn sphere, the analysis of
\cite{B2} shows that $R(\Si, SU(3))$ is a union of path
components, each of which is homeomorphic to either an isolated point
or a 2-sphere. More precisely, we will show that
each  path component is
one of the following four types:

\begin{description}
\item[Type Ia] Isolated conjugacy classes of irreducible
representations.
\item[Type IIa]  Smooth 2-spheres, each parameterizing a family of
conjugacy classes of irreducible representations.
\item[Type Ib]  Isolated conjugacy classes of nontrivial reducible
representations.
\item[Type IIb]   {\em Pointed 2-spheres}, each parameterizing a family
of conjugacy classes
of representations, exactly one of which is reducible.
\end{description}

The main result in this paper is the following theorem (Theorem 6.2),
which describes
how each of the component types contributes to the $SU(3)$ Casson
invariant.  This, together with enumerations of the components of each
type, enable us to calculate the invariant for a variety of Brieskorn
spheres $\Si (p,q,r)$.  The results of these computations can be found in Tables
1 and 2.

\medskip \noindent
{\bf Theorem.}  {\it   Type Ia, IIa, Ib, and IIb components each contribute +1,
+2, 0, and +2, respectively, to the integer valued $SU(3)$ Casson
invariant $\tau_{SU(3)}(\Si(p,q,r)) $.}
\medskip

We conclude the introduction by outlining the  proof of this theorem.
Components of Type Ia are regular and remain so after small
perturbations.
The sign attached to each such point is positive by the results of
\cite{B2}, and so computing the contribution of the Type Ia points to
$\tau_{SU(3)}(\Si)$ reduces to an enumeration problem. This
is carried out in Section \ref{s6}.

Components of  Type IIa are nondegenerate
critical submanifolds of the Chern-Simons function.
Bott-Morse theory,
together with a spectral flow computation, implies
that each such component contributes $\chi(S^2)=2$ to
$\tau_{SU(3)}(\Si)$.  Thus the computation of
the contribution of the Type IIa
components to
$\tau_{SU(3)}(\Si)$  is also reduced to an enumeration problem
which is solved
in Section \ref{s6}.

Components of  Type Ib do not contribute to $\tau_{SU(3)}(\Si)$
(although they do enter into the calculations of
the invariant $\la_{SU(3)}$ given in \cite{BHKK}).

The only remaining issue is to calculate the
contribution of components of Type IIb.
This requires some sophisticated techniques
that go beyond those of \cite{BHK}, where one can find
computations
of $\tau_{SU(3)}$ for
Brieskorn spheres of the form $\Si(2,q,r)$ (whose representation
varieties do not contain any Type IIb components).
The problem is that Type IIb components are singular in a strong sense:
even their lifts to the
based moduli space are singular.
We introduce a perturbation  which resolves these
singularities and then carefully analyze its effect on
the topology of the moduli space.  We prove that after applying 
the perturbation,  each pointed 2-sphere resolves into two
pieces, one isolated gauge orbit of reducible connections
and the other a smooth, nondegenerate 2-sphere of
gauge orbits of irreducible connections (similar to a Type IIa
component).

In defining the perturbation, we regard one of the singular
fibers of the Seifert fibration $\Si \to S^2$
as a knot in $\Si$ and perturb
the flatness equations in a small
neighborhood of this knot.
Consequently, perturbed flat connections are seen to be
flat on the knot complement, and we
study the perturbed flat moduli space in terms
of the
$SU(3)$ representation space  of this knot complement.
Basically, the  perturbed flat moduli space on $\Si$  is
obtained from the flat moduli space of the knot complement by replacing
the condition ``meridian is sent to the identity'' by a condition
of the form ``the meridian and longitude are related by a certain
equation.''

Having resolved the singularities in the Type IIb components,
we then determine the contribution of the reducible, perturbed flat
connection to the correction
term. This is given by the spectral flow (with $\CC^2$ coefficients) of
the
odd signature operator. To calculate this we
prove a splitting theorem for spectral flow determined by the decomposition
of $\Si$  into a knot complement and a solid torus.

\bigskip
\noindent{\bf Notation.}
If $\pi$ is a discrete group and
$\al\colon \pi \to G$ is a representation, we denote
the stabilizer subgroup of $\al$ by
$$\Ga_\al = \{ g \in G \mid g \al g^{-1} = \al\}.$$
If
$G$ is a Lie group, the orbit  of $\al$ under conjugation
is smooth and
diffeomorphic to the homogeneous manifold $G/\Ga_\al$.
We denote the
representation variety
$$R(\pi,G) = \Hom(\pi,G)/\text{conjugation}.$$
Given a representation $\al\colon \pi \to G,$ we denote its
conjugacy class by $[\al].$
Given a manifold $X$, we denote by $R(X,G)$  the
representation variety of the fundamental group
$\pi_1(X).$

\section{SU(3) representation spaces of Brieskorn spheres}
\label{s2}
In this section, we identify the components of  the $SU(3)$
representation
varieties of Brieskorn spheres $\Si$, both as topological spaces and as
varieties with their Zariski tangent spaces.
Crucial to our discussion are computations of the
   twisted cohomology groups which reflect the local structure
   of the representation varieties. The global structure of
   the representation variety is presented in
Theorem \ref{sthm}, which gives
a complete classification of
the different path components of $R(\Si,SU(3))$.

\begin{figure}[t]
\begin{center}
\leavevmode\hbox{}
\includegraphics[width=2.5in]{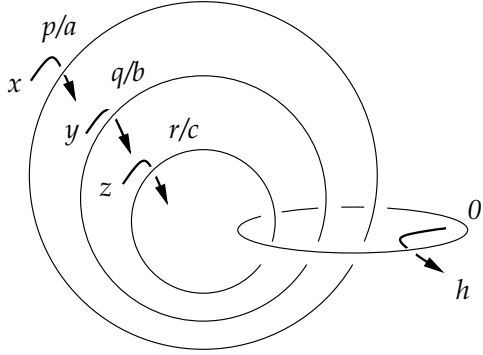}
\caption{A surgery description of the Brieskorn manifold $\Si(p,q,r)$
indicating the Wirtinger generators $x,y,z,$ and $h$ for $\pi_1(\Si)$.}
\label{earring}
\end{center}

\end{figure}

\subsection{Brieskorn spheres}
\label{BS}

Given integers $p,q,r$, set
$$\Si (p,q,r) = \{(x,y,z) \in \CC^3 \mid x^p + y^q+z^r = 0\} \cap S^5.$$
If $p,q,r$ are pairwise relatively prime then
$\Si (p,q,r)$ is a homology 3-sphere and has surgery
description
in Figure \ref{earring} (see \cite{S} for details).  Here,   $a,b,c$
satisfy
\begin{equation}\label{hsc}
a qr   +  b pr + c pq = 1.
\end{equation}
The resulting manifold $\Si (p,q,r)$ is independent
of $a,b,c$, up to orientation preserving homeomorphism.
  Without loss of generality we assume that $p$ and $q$ are odd.

\begin{prop}  \label{forceit}The numbers $a$ and $b$ can be chosen to
be equal.
\end{prop}
\begin{proof}  Since $p,q,$ and $r$ are pairwise relatively prime,
$r(p+q)$ and $pq$ are relatively prime.  Thus
there are integers $a$ and $c$ such that
$$ar(p+q) + cpq = 1,$$
which is equivalent to the condition (\ref{hsc}) with $b=a$.
\end{proof}

Fix integers $a$ and $c$ as above.
Note that since $p$ and $q$ are both odd, $c$ must
also be odd.
A presentation for the fundamental group of $\Si(p,q,r)$ is
\begin{equation} \label{pres}
\pi_1\big( \Si(p,q,r) \big)= \langle x,y,z,h\mid ~ x^p=y^q=h^a,~ z^r =
h^c,~ xyz=1,
h\mbox{ is
central} \rangle,
\end{equation}
   where $x,y,z$ and $h$ are the Wirtinger generators indicated
in Figure \ref{earring}.

Whenever $p,q,$ and $r$ are clear from  the context, we drop
them from the notation and   denote the Brieskorn
sphere by $\Si$.
A  regular neighborhood of the singular $r$-fiber in $\Si$
is a solid torus
whose boundary torus $T$ splits the Brieskorn
sphere $\Si=Y\cup_{T} Z$,
where $Y=D^2 \times S^1 $ is the solid torus and $Z = \Si - Y$ is its
complement.
Alternatively, $Z$ is the  complement of an open tubular neighborhood of 
the core of the $\left(\frac{r}{c}\right)$ curve in $\Si$ and depicted in Figure
\ref{earring}.
With regard to the natural peripheral structure thus obtained on
$Z$, its fundamental group has presentation
\begin{equation} \label{pres-Z}
\pi_1(Z) = \langle x,y, h\mid ~ x^p = y^q= h^a, h \mbox{ is central}
\rangle.
\end{equation}
In terms of these generators, the meridian and longitude are
represented by
\begin{equation}\label{mu}
\mu  = (xy)^r h^c \qquad \hbox{and} \qquad \la = (xy)^{pq} h^{-(p+q)a}.
\end{equation}
Then $\mu$ generates the  abelianization of $\pi_1(Z)$, and one can
check that in  $H_1(Z)$,
\begin{equation}\label{homolog}
  [x]=aq [\mu], \ [y]= ap  [\mu], \ [h]=pq [\mu], \text{ and }
[\la]=0.\end{equation}

   \subsection{Decompositions of $\boldsymbol{su(3)}$}
In this subsection, we  examine the restriction of
the adjoint representation of $SU(3)$ on its Lie algebra $su(3)$ to
various subgroups.

   Consider first the subgroup
$$SU(2)\times \{1\}=\left\{
\begin{pmatrix} A &0\\ 0&1\end{pmatrix}\,
\Bigm|  A \in SU(2)       \right\} \subset SU(3).$$
Then $su(3)$ decomposes invariantly with respect to the
adjoint action of $SU(2)\times
\{1\}$ as
\begin{equation}\label{dec1} su(3) = su(2) \oplus \CC^2
\oplus \RR,\end{equation}
   where $SU(2)\times \{1\}$ acts by the adjoint action on $su(2)$,
   by the defining representation on $\CC^2$, and trivially on $\RR$.

More generally, consider the subgroup
$$ S(U(2)\times U(1))=\left\{
\begin{pmatrix} A&0\\  0&\det{A}^{-1}\end{pmatrix}
   \,\Bigm|\,
   A \in U(2) \right\} \subset SU(3). $$
The  decomposition of the Lie algebra $su(3)$ takes the form
\begin{equation}\label{dec2} su(3)  = s(u(2)\times u(1))
\oplus
\CC^2,\end{equation} where $S(U(2)\times U(1))$ acts on the first
factor  via the adjoint representation and on the second factor by
$$\left[\begin{array}{ccc} t a & t b &0\\
-t \bar{b} & t \bar{a} & 0 \\
0&0& t^{-2}\end{array} \right] \cdot
\left[\begin{array}{c} z_1 \\ z_2 \end{array}
\right] = t^3 \left[\begin{array}{c} a z_1 +b z_2 \\
-\bar{b} z_1+\bar{a} z_2  \end{array}\right] \ (|t|=1,\  |a|^2 +
|b|^2=1).$$

There is a
canonical isomorphism
$ S(U(2) \times U(1)) \cong U(2).$ However,
the action of $ S(U(2) \times U(1))$
on $\CC^2$ is not the
standard $U(2)$ action, even though its restriction to the subgroup
$SU(2)\times \{1\}$ is standard.

Every $SU(3)$ matrix is diagonalizable. We
parameterize the diagonal matrices using the map
$\Phi\colon \RR^2 \to  SU(3)$ given by
\begin{equation} \label{defnPhi}
\Phi(u,v)=\left[ \begin{array}{ccc} e^{i(u+v)}&0&0\\
0& e^{i(-u+v)}&0\\ 0&0& e^{-2iv}\end{array} \right].
\end{equation}
With respect to the decomposition (\ref{dec2}), the matrix $\Phi(u,v)$
acts on $\CC^2$ by
$$ \Phi(u,v) \left[\begin{array}{c} z_1 \\z_2  \end{array}
\right]
= e^{3iv} \left[\begin{array}{c} e^{iu} z_1 \\   e^{-iu} z_2
\end{array}
\right] .$$
Note that the centralizer
of $S(U(2)\times U(1))$ is $\{ \Phi (0,v)
\}$, and this circle acts trivially on $s(u(2)\times u(1))$ and with
weight three on $\CC^2$.

\subsection{Cohomology calculations}
\label{CC}

In this subsection, we give computations of  $H^i(\Si;su(3)_\al)$,
where $\al\colon  \pi_1(\Si) \to SU(3)$ is a representation and
$SU(3)$ acts on its
Lie algebra $su(3)$ via the adjoint representation.
First, we establish some notation and recall some basic facts.
Let $X$ be a cell complex,   $G$  a Lie group, $V$ a vector space
on which $G$ acts, and
$\al\colon \pi_1(X) \to G$   a representation.
Denote by $V_\al$
the  local coefficient system determined by $\al$
and by $H^i(X;V_\al)$ the $i$-th cohomology group
of $X$ with twisted coefficients in $V_\al.$
Although some of the cohomology groups we consider have natural
complex structures, we  use the notation  $\dim(H)$ to
refer to the dimension of $H$ as a real vector space.

Given a finite complex $X$ and representation $\al \colon \pi=\pi_1(X)\lto GL(V)$,
we can identify $H^i(X;V_\al)\cong H^i(\pi;V_\al)$ for $i=0,1$.  Group cohomology
$H^*(\pi;V_\al)$ can be computed from the reduced bar resolution.
In this model, the space of  (twisted) $i$-cochains is
given by a set of functions:
  $$C^0(\pi;V)=V, \ C^i(\pi;V)=\{f\colon \stackrel{\text{\tiny $i$ terms}}{\pi \times \cdots\times \pi} \lto V\}, i>0.$$
We will only need the formulas for the $0$-th and $1$-st coboundary
  operators, $$d^0(v)(\gamma)=(\gamma-1)\cdot v, \ d^1(f)(\ga_1,\ga_2)=
f(\ga_1)+\ga_1\cdot f(\ga_2) -f(\ga_1\ga_2).$$  The  {\em Fox calculus}
provides a means to calculate the 1-cocycles,
i.e.   the solutions $f \in C^1(\pi;V_\al)$ to  the equation $d^1(f)=0$.
 Given a presentation
$\pi=\langle x_1,\ldots, x_m \mid r_1,\ldots, r_n\rangle$,
    the Fox derivative of a relation $r_j$ with respect to a cocycle
$f$ is  the element of $V$ obtained by using the equation
$f(\ga_1\ga_2)=f(\ga_1)+\ga_1\cdot f(\ga_2)$ inductively to express
 $0=f(1)=f(r_j)$ in terms of $X_i=f(x_i)$. Note that
the map $f\mapsto (f(x_1),\ldots , f(x_n))$ is injective on 1-cocyles.
Most of our computations
 are given in terms of group cohomology, but  occasionally we
 will make use of topological
tools such as Poincar\'e duality and the Euler characteristic.

We now explain the  relationship  between representation varieties and
these cohomology groups.
Suppose that $G$ is a compact Lie group, acting on its
Lie algebra ${\mathfrak g}$  via the adjoint action,   and $\pi$  is a
finitely presented group. Then the Zariski tangent space to (the
algebraic variety) $R(\pi,G)$ at the conjugacy class of a
representation $\al\colon \pi\to G$ is
isomorphic to   $H^1(\pi; {\mathfrak g}_\al)$.  Moreover, $\dim
H^0(\pi;{\mathfrak g}_\al) = \dim \Ga_\al,$
where $\Ga_\al\subset G$ denotes the stabilizer subgroup of $\al$ under
the conjugation action of $G$
   on $\Hom(\pi,G)$. Equivalently,  $\Ga_\al$ equals
the centralizer of $\im(\al)$.

The Kuranishi map embeds a neighborhood of $[\al]$ in $R(\pi,G)$ into its
Zariski tangent space modulo $\Ga_\al$. In particular if $H^1(\pi;{\mathfrak g}_\al)=0$, then
$[\al]$ is an isolated point in $R(\pi,G)$ (although the converse is
sometimes false). We say  that $[\al] \in R(\pi,G)$ is a {\em smooth
point} if  a neighborhood of $[\al]$ in $R(\pi,G)$ is homeomorphic to
$H^1(\pi, {\mathfrak g}_\al)$; otherwise $[\al]$ is called a {\em singular
point}.

 We are mostly interested in the case when $G=SU(3)$ and $V=su(3)$. 
 For reducible representations, we are interested in the case $G=SU(2)$
 and $V=su(2)$ or $\CC^2$.  
 To see why, note that up to conjugation,
 any reducible  representation $\al\colon \pi_1(\Si) \to SU(3)$ has
 image in the subgroup $S(U(2)\times U(1)).$
 Since $\Si$ is a   homology sphere, it follows that
    $\al$ has image in $SU(2) \times \{1\}.$
    Using the decomposition
(\ref{dec1}) we conclude that
$$ H^i(\Si;su(3)_\al)=  H^i(\Si;su(2)_\al)\oplus
   H^i(\Si;\CC^2_\al)\oplus  H^i(\Si;\RR),$$
where the first cohomology group has  coefficients $su(2)$
twisted via
the adjoint action (viewing $\al$ as an $SU(2)$ representation),
the second has coefficients $\CC^2$ twisted by
the standard representation, and the last has untwisted
real coefficients.

\begin{prop} \label{coho}
Suppose $\al\colon \pi_1(\Si) \to SU(3)$  is a nontrivial
representation. Then
$\al$ has nonabelian image. Moreover:
\begin{enumerate}
\item[(i)] If $\al$ is irreducible, then
   $\al(h) =e^{2\pi i k/3} I$ for an integer $k$  and
   $$H^1(\Si;su(3)_\al) =
\begin{cases} \RR^2 & \text{if $\al(x), \al(y),$ and $\al(z)$
each have three distinct eigenvalues,} \\
0 & \text{otherwise.}
\end{cases} $$
\item[(ii)] If $\al $ is reducible and has been
conjugated to take values in $SU(2) \times \{1\},$
then
$$\al(h) = \left[ \begin{array}{ccc} \pm 1 &0& 0 \\
   0& \pm 1 & 0 \\
   0 &0&1
   \end{array} \right].$$
With respect to the splitting $su(3) = su(2) \oplus \CC^2 \oplus \RR$
(see equation (\ref{dec2})), we have
$H^0(\Si;su(2)_\al)=0, \ H^0(\Si;\CC^2_\al)=0,  \ H^1(\Si;su(2)_\al)=0$
and
$$H^1(\Si;\CC^2_\al) =
\begin{cases} \CC^2 & \text{if $\al(h)=I$,} \\
0 & \text{otherwise.}
\end{cases} $$
\end{enumerate}
\end{prop}

   \begin{proof}

   First note that since $h$ is central in $\pi_1(\Si)$,
   $\al(h)$  lies in the centralizer of $\im(\al)$.

We now prove (i). Suppose $\al$ is
   irreducible. Then $\Ga_\al$ is the center of $SU(3)$, and hence is
discrete.
Thus  $\al(h)$ is central  and   $\dim H^0(\Si;su(3)_\al)=0$.

   Set $\pi = \pi_1(\Si)$
   and let
$ B^i (\pi;su(3)_\al)  = \im(d^{i-1})$ be
the coboundaries and $Z^i (\pi;su(3)_\al) = \ker(d^i)$
   the cocycles in the reduced bar complex. Thus
   $$H^1(\pi;su(3)_\al) = Z^1 (\pi;su(3)_\al) / B^1 (\pi;su(3)_\al).$$
Since $H^0(\pi;su(3)_\al)=0,$
   $d^0$ is injective and $B^1(\pi;su(3)_\al)$ has dimension 8.
   So to compute $H^1(\pi;su(3)_\al)$ we only need to determine the
dimension of
   the space $Z^1(\pi;su(3)_\al)$ of  1-cocycles.

The Fox calculus  identifies
$ Z^1 (\pi;su(3)_\al) $ with the set of
4-tuples $(X,Y,Z,H)$ in $su(3)$ satisfying the equations one gets
by taking Fox derivatives of the relations in (\ref{pres}).
   For example, the relation $hx =xh$ gives the equation
    $$ H+ \al(h)X = X+ \al(x) H.$$ Since $\al(h) X = X$, this reduces
    to $H=\al(x)H.$ Similarly, we get the equations
    $H=\al(y) H$ and $H=\al(z)H.$
    Since $\al$ is irreducible with image generated
    by $\al(x),\al(y)$, and $\al(z)$, these
   three equations imply
$H=0.$

Setting $H=0$ in  the remaining equations, we obtain:
   \begin{eqnarray*}
   (1+\al(x) + \cdots + \al(x^{p-1})) X &=& 0, \\
   (1+\al(y) + \cdots + \al(y^{q-1})) Y &=& 0, \\
   (1+\al(z) + \cdots + \al(z^{r-1})) Z &=& 0, \\
   X + \al(x) Y + \al(xy) Z &=& 0.
   \end{eqnarray*}

\medskip \noindent
{\sc Case 1:}  $\al(x), \al(y)$ and $\al(z)$
all have three distinct eigenvalues.

Since $\al(x)^p=\al(h)^a$ acts as the identity on $su(3)$ via the adjoint
action, $su(3)$ decomposes as $T_x\oplus U_x$, where $T_x$ is the tangent space
to the maximal torus containing $\al(x)$ and
$U_x$ is the kernel of the map $1+\al(x)+\cdots+\al(x)^{p-1}\colon su(3)\to
su(3)$. Note that $T_x$ is 2-dimensional
(since $\al(x)$ has three distinct eigenvalues) and that $\al(x)$ acts trivially on
$T_x$.
It follows from the equations above that $X$ lies in $U_x$.
Similar statements hold for $\al(y)$ and $\al(z)$. The space of 1-cocycles is
therefore a subspace of $U_x\oplus U_y\oplus U_z$.

Since $\al$ is irreducible,
$T_x\cap T_y=0$. In fact, if $t\in T_x\cap T_y$, then  exp$(t)\in SU(3)$
stabilizes both $\al(x)$ and $\al(y)$ and hence stabilizes $\al$.
Thus $U_x\cap U_y$ is 4-dimensional, and therefore $U_x+U_y$ is 8-dimensional,
i.e. $U_x+U_y=su(3)$.

Since $\al(x)^{-1}$ preserves the decomposition $su(3)=T_x\oplus U_x$ and acts
as an isomorphism on each factor, the linear map
$$U_x\oplus U_y\to su(3), \ (X,Y)\mapsto \al(x)^{-1}X + Y$$
is onto. Thus the linear map
$$(X,Y,Z)\mapsto \al(x)\big(\al(x)^{-1}X+Y\big) + \al(xy)Z=X+\al(x)Y+\al(xy)Z$$
is also onto. Its kernel is just the space of 1-cocycles, and so
$\dim Z^1(\pi;su(3)_\al)=10$. Hence $\dim H^1(\pi;su(3)_\al)=10-8=2$.

\medskip \noindent
{\sc Case 2:}  One of $\al(x), \al(y)$ and $\al(z)$
has a double eigenvalue.

We first show that
at most one of $\al(x), \al(y),$ and $\al(z)$
   can have a double eigenvalue. For example, if both
   $\al(x)$ and $\al(y)$ had a double eigenvalue, then
   the intersection of the
   corresponding  eigenspaces would
determine a    linear subspace invariant under $\al(x)$, $\al(y)$, and
$\al(z)=\al(xy)^{-1}$,  contradicting the irreducibility of $\al.$

So assume  that $\al(x)$ has a double
eigenvalue and $\al(y)$ and $\al(z)$ have three distinct
eigenvalues (the proofs of the other cases work the same way). Under
the adjoint action of $\al(x)$,
$su(3)$ decomposes as $\CC^2 \oplus \RR^4$,
where $ \al(x)$ acts trivially on $\RR^4$
and by multiplication by a nontrivial $p$-th
root of unity on $\CC^2$.
Thus we see that $X$ now lies in a (real) 4-dimensional subspace
$\CC^2 \subset su(3)$.
Arguing as before, we
conclude that $\dim Z^1(\pi;su(3)_\al) = 4+6+6 -8=8$,
from which it follows that
$$\dim H^1(\pi;su(3)_\al)
= \dim Z^1(\pi;su(3)_\al) - \dim B^1(\pi;su(3)_\al)
= 8-8 =0.$$

These two cases complete the proof of (i) because
   irreducibility of $\al$ precludes any other possibility.
   To see this, suppose  one of $\al(x), \al(y)$ or $\al(z)$ were central,
say $\al(x)$,
   then the relation $xyz=1$
would imply that $\al(y)$ commutes with $\al(z)$, and
hence that $\al$ is abelian. This would imply $\al$ is trivial (and in particular reducible).

In proving (ii), we regard $\al$ as an $SU(2)$ representation.
Irreducibility of $\al$ (as an $SU(2)$ representation)
implies that $H^0(\Si;su(2)_\al)=0$
   and $\al(h)=\pm I$.
The fact that  $H^1(\Si;su(2)_\al)=0$ for Brieskorn spheres is
well-known
(see  \cite{FS}).
Nontriviality of $\al$ implies $H^0(\Si;\CC^2_\al)=0$,
and  that leaves $H^1(\Si;\CC^2_\al)$, which we determine
with another application of
   the Fox calculus.
   The only
difference is that we use the defining
representation instead of  the adjoint representation.
In particular, $-I \in SU(2)$ acts nontrivially.

Suppose then that $(X,Y,Z,H)$ is a 4-tuple of
vectors in $ \CC^2$ satisfying the
equations one gets by taking Fox derivatives of
the relations in (\ref{pres}).
There are two cases.

\medskip \noindent
{\sc Case 1:}  $\al(h) = I.$

As before,   $H=0$ and
   $\al(x)$ acts by  multiplication by a nontrivial $p$-th
root of unity in each of the two complex factors.
Consequently $ (1+\al(x) + \cdots + \al(x^{p-1}))X=0$
for all $X \in \CC^2.$
Similar statements hold for $Y$ and $Z$ and it follows that
$(X,Y,Z,H)$ is a 1-cocycle
provided
$H=0$ and $X + \al(x) Y + \al(xy) Z =0.$
One can check that the last equation imposes
 four independent conditions, hence
   $\dim Z^1(\pi;\CC^2_\al) = 4+4+4-4=8$,
and it follows that
$$\dim H^1(\pi; \CC^2_\al)
= \dim Z^1(\pi; \CC^2_\al) - \dim B^1(\pi; \CC^2_\al)
= 8-4 =4.$$

\medskip \noindent
{\sc Case 2:}  $\al(h) = - I.$

In this case, $\al(h)$ acts on $\CC^2$ by multiplication
by $-1$ and it
   is no longer true that $H=0.$
Instead, we find that  $H$ determines $X$
by the equation $(1-\al(h))X =(1-\al(x)) H$
and similarly for $Y$ and $Z$.
An easy check shows that all the remaining equations
are automatically satisfied, and since $H$ is an arbitrary
element in $\CC^2,$ it follows that $\dim Z^1(\pi;\CC^2_\al) = 4$
and
$$\dim H^1(\pi; \CC^2_\al)
= \dim Z^1(\pi; \CC^2_\al) - \dim B^1(\pi; \CC^2_\al) =4-4=0$$
as claimed.
\end{proof}

\subsection{The representation variety R(${\boldsymbol \Si}$,SU(3))}
In this subsection, we classify the different
path components of the representation variety $R(\Si,SU(3)).$
To start off, we show that every component contains at most
one conjugacy class of reducible representations.

\begin{prop} \label{onlyone} If $\al_t,\ t\in [0,1],$ is a continuous
path of $SU(3)$ representations
of $\pi_1(\Si)$ with $\al_0$ and $\al_1$ both reducible, then
$\al_0$ and $\al_1$ are conjugate. Consequently, every path component of
  $R(\Si,SU(3))$ contains at most one  conjugacy class
  of reducible representations.
\end{prop}
\begin{proof}  
For the trivial representation $\th$, $H^1 (\Si; su(3)_\th) = H^1 (\Si;
\RR^8)=0$, so $[\th]$ is isolated.  Thus we assume that $\al_t$   is
nontrivial for all $t$.
If $\al_0(h) \neq I$, then
Proposition \ref{coho} implies that   $[\al_0]$ is isolated.
So we can assume that
$\al_0(h) = I$.  The continuous function $t\mapsto \tr(\al_t(h))$ takes
values in the discrete set $\{3, -1, 3e^{2\pi i/3}, 3e^{4\pi i/3}\}$ by
Proposition \ref{coho}.
It follows that $\al_t(h)=I$ for all $t$.
The relations (\ref{pres}) then imply that
$\al_t(x), \  \al_t(y),$ and $\al_t(z)$ are  conjugate to fixed $p$-th,
$q$-th, and
$r$-th roots of unity in $SU(3)$ for all $t$.
(To see this, use continuity and the fact that
the trace map $\tr\colon SU(3) \to \CC$ distinguishes conjugacy classes
and sends the set
$\{ A \in SU(3) \mid A^k = I\}$ of all
$k$-th roots of unity into a discrete set.)

Since $\al_0$ and $\al_1$ are both reducible and $SU(3)$ is path
connected, we may assume that the path $\al_t$ is conjugated so that
$\al_0$ and $\al_1$   take values in  $SU(2)\times \{1\}$.
Thus $\al_0(x)$ and $\al_1(x)$ each have one eigenvalue equal to 1. But since
$\al_0(x)$ and $\al_1(x)$  are conjugate (in $SU(3)$), the other two
eigenvalues of $\al_0(x)$ and $\al_1(x)$ coincide.  The same argument
applies to $y$ and $z$.

    It is well-known  that
the conjugacy class $[\be]$ of a  representation
$\be\colon \pi_1(\Si) \to SU(2)$ of a Brieskorn sphere
is completely determined by
the eigenvalues  of
$\be(x), \be(y),$ and $\be(z)$ (see \cite{FS}).  Hence $\al_0$ and
$\al_1$ are conjugate as $SU(2)$ and hence also as $SU(3)$
representations.
\end{proof}

\begin{prop}\label{onlyone2} Every path component of $R(\Si,SU(3))$
is either an
isolated point, a smooth 2-sphere consisting
of conjugacy classes of irreducible representations, or a pointed
2-sphere, which is smooth except for exactly one singular point,
the  conjugacy class of a  reducible representation.
\end{prop}
\begin{proof} It is proved in \cite{B2, Ha} that  each
path component of
$R(\Si,SU(3))$ is either an isolated point or a topological 2-sphere.
In the case  of an isolated point,
there is nothing to prove, so assume the path component is a 2-sphere.
Any conjugacy class $[\al]$ of irreducible representations lying
on such a component must have nonzero Zariski
tangent space, and Proposition \ref{coho} then
implies $H^1(\Si;su(3)_\al) \cong \RR^2$
and we conclude that $[\al]$ is indeed a smooth point of $R(\Si,SU(3))$.
On the other hand,
Proposition \ref{onlyone} shows that every path
component of $R(\Si, SU(3))$ contains at most one
conjugacy class  of reducible representations.
For  a pointed 2-sphere component, the  conjugacy class
$[\be]$ of reducible
representations is never a smooth point,
since Proposition \ref{coho} shows its Zariski
tangent space is $H^1(\Si;su(3)_\be) \cong \RR^4$.
(Note that the hypothesis on $\be$ implies that
$H^1(\Si;su(3)_\be) \neq 0$, and then Proposition \ref{coho}
shows that $\be(h) = I$.)
\end{proof}

   The next proposition shows that the pointed 2-spheres are in
one-to-one
correspondence with the nontrivial  reducible representations sending
$h$ to
the identity.

\begin{prop} \label{TFAE}
Given a nontrivial reducible representation
$\al \colon  \pi_1(\Si)  \to SU(3),$ the following
are equivalent:
\begin{enumerate}
\item[(i)] $\al(h) = I$,
\item[(ii)] $H^1(\Si;\CC^2_\al) \neq 0$,
\item[(iii)] There exists a   family of
irreducible
$SU(3)$ representations limiting to $\al$.
\end{enumerate}
The collection of pointed 2-spheres in $R(\Si,SU(3))$ are
therefore in one-to-one correspondence with conjugacy classes of
nontrivial reducible representations $\al\colon \pi_1(\Si) \to SU(3)$
with $\al(h) = I.$ Further, $\tr \al(z)$ is constant along a pointed 2-sphere.
\end{prop}
\begin{proof}
The statement (i) $\Leftrightarrow$ (ii) follows from
Proposition \ref{coho}, (ii).
The implication (iii) $\Rightarrow$ (ii) follows
because the Kuranishi map locally  embeds $R(\Si, SU(3))$
near $[\al]$ into its Zariski
tangent space $H^1(\Si; su(3)_\al)$ modulo $\Ga_\al$, and the Zariski tangent space
equals $H^1(\Si; \CC^2_\al)$ by
Proposition \ref{coho}.

For the implication (i) $\Rightarrow$ (iii), notice that a
representation  $\al\colon \pi_1(\Si) \to SU(3)$ satisfying  $\al(h) =
I$
uniquely determines an $SU(3)$  representation of the  (free) group
$F= \langle x,y,z \mid xyz = 1\rangle$.
Fix three conjugacy classes $\ba,\bb,\bc$ in $SU(3)$
and consider the space
$\cM_{\ba\bb\bc}$  consisting  of
conjugacy classes of representations
$\al\colon F \to SU(3)$
with $\al(x) \in \ba, \ \al(y)\in \bb,$ and $\al(z) \in \bc.$
In \cite{Ha}, Hayashi gives necessary and sufficient  conditions on
$\ba,\bb,\bc$ for $\cM_{\ba\bb\bc}$
to be nonempty.
The resulting
  inequalities (18 in all)   determine a convex, 6-dimensional
polytope $P$ parameterizing  all triples $(\ba, \bb ,\bc)$
with $\cM_{\ba\bb\bc} \neq \varnothing.$
Hayashi observes further that $\cM_{\ba\bb\bc}$ is a 2-sphere
whenever $(\ba, \bb ,\bc)$ lies in the interior of $P$ and  is a point
whenever $(\ba, \bb ,\bc)$ lies on the boundary of $P$.
For more details, turn to Subsection \ref{sfuse}
and  read Theorem \ref{fused}.

The key to proving that (i) $\Rightarrow$ (iii) is  to show
that the triple $(\ba, \bb ,\bc)$ determined by $\al(x), \al(y), \al(z)$
lies in the interior of $P$. From this, it follows that $
\cM_{\ba\bb\bc}$, which is connected and contains
$[\al]$, is a 2-sphere. Assume to the contrary that
$(\ba, \bb ,\bc)$ is a boundary point of $P$.
  There are two possibilities, because there
are two kinds of boundary points.
The first kind occurs when one of the inequalities in equation
(\ref{fusion})
is an equality. This cannot happen for $(\ba, \bb ,\bc)$ because
$\al(x), \al(y), \al(z)$  are, respectively,
$p$-th, $q$-th, and $r$-th roots of unity in $SU(3)$
and $p,q,r$ are pairwise relatively prime. The other kind of boundary
point of $P$ occurs when one of
$\al(x), \al(y), \al(z)$ has a repeated eigenvalue.
If $\al(x)$  were to
  have a repeated eigenvalue, then since $\al$ has image in $SU(2)
\times \{1\}$
(up to conjugation),
  it follows that  $1$ is an eigenvalue of $\al(x)$, and so its other
eigenvalues
   are either both $+1$ or  both $-1.$
  In either case, it follows easily that $\al(x)$ commutes
  with $\al(y)$ and $\al(z)$, and the relation $xyz=1$
then shows that $\al$ has abelian image. Since $\Si$ is a homology
sphere, this implies $\al$ is trivial and gives the desired
contradiction.
\end{proof}

The computations of Propositions  \ref{coho}--\ref{TFAE} give a
decomposition
of $R(\Si,SU(3))$ into the various different types, summarized in
the following theorem.

\begin{thm}\label{sthm}
The path components of the representation space
$R(\Si, SU(3))$  come in the following  four types.
(Notice that in each case, the image of $h$ is constant
along the component and the conjugacy classes of the images
of $x,y,$ and $z$ are also constant.)
\begin{enumerate}
\item[(i)] The Type Ia  components  consist of one isolated  conjugacy
class
$[\al]$ of irreducible
representations with the property that exactly one of
$\al(x),\al(y),\al(z)$
has a repeated eigenvalue.
The representation $\al\colon \pi_1(\Si) \to SU(3)$ sends $h$ to a
central element
and has $H^1(\Si; su(3)_\al)=0$.

\item[(ii)] The Type IIa components  are smooth 2-spheres consisting of
conjugacy classes of irreducible
representations $\al$ with the property that $\al(x),\al(y),\al(z)$
each
have three distinct eigenvalues. For any conjugacy class
$[\al]$ in a Type IIa component,
the representation $\al\colon \pi_1(\Si) \to SU(3)$ sends
$h$ to a central element and has
$H^1(\Si; su(3)_\al)\cong \RR^2$.

\item[(iii)] The Type Ib components consist of isolated conjugacy
classes $[\be]$
of reducible representations.
The only isolated,  reducible conjugacy class $[\be]$
with $\be(h)=I$ is the conjugacy class of the trivial representation.
If $[\be]$ is isolated,  reducible,
and nontrivial, then
$\tr(\be(h))=-1$ (i.e.~ $\be(h)=-I$ as an $SU(2)$ element) and
$H^1(\Si; \CC^2_\be)= 0$.  
\item[(iv)] The Type IIb components are topological 2-spheres containing
exactly one
conjugacy class $[\be]$ of reducible representations with
$H^1(\Si;su(3)_\be)=H^1(\Si;\CC^2_\be)\cong \RR^4$.
Every other conjugacy class $[\al]$ in a Type IIb component is a smooth
point with $\al$ irreducible and satisfying $H^1(\Si;su(3)_\al)\cong \RR^2$. In particular, the reducible orbit is the only singular
point.
   Every conjugacy class of representations
in a Type IIb component sends $h$ to the identity  and
sends $x, y$ and $z$ to elements with
three distinct eigenvalues.
   \end{enumerate}
\end{thm}

The way in which a component type contributes to
the integer valued $SU(3)$ Casson invariant
is explained in Theorem \ref{cct}.

\begin{prop}
The representation variety
$R(\Si(p,q,r),SU(3))$ contains a  Type IIb  component if (and only if)
none of $p,q,r$ equal $2.$
\end{prop}

\begin{proof} Suppose first that $r=2$ and
$\al\colon \pi_1 \big(\Si(p,q,2)\big) \to SU(2)$ is
a  representation with $\al(h) = I.$
Then $\al(z)^2 = I,$ hence $\al(z) = \pm I$ is central.
Thus $\al(y)=\pm \al(x)^{-1}$, which implies
$\al$ is abelian and hence trivial. Thus, up to reordering, if one of
$p,q,r$ equals 2,
then $R(\Si(p,q,r),SU(3))$ does not contain a Type IIb
component.

On the other hand, if none of $p,q,r$ equals $2,$ the results of
\cite{FS} prove the existence of nontrivial
representations $\al \colon \pi_1 \big(\Si(p,q,r)\big) \to SU(2)$
with $\al(h) = I$.
Apply  Proposition
\ref{TFAE} to complete the proof.
\end{proof}

\section{SU(3) representation spaces of knot complements} \label{s3}
We next carry out a similar analysis of  the
$SU(3)$ representation
variety $R(Z,SU(3))$ of the knot complement $Z$ obtained by
removing a neighborhood of one of the singular fibers
of $\Si(p,q,r).$

We explain our  purpose first.
The inclusion $Z \hookrightarrow \Si$ induces
a surjective map $\pi_1(Z) \to \pi_1(\Si)$.
In terms of the presentation (\ref{pres-Z}), this map is
given by imposing the relation $\mu=1$.
Consequently the representation
variety $R(\Si,SU(3))$ can be viewed as
the subvariety of $R(Z,SU(3))$
cut out by the equation  determined
by the condition that ``the meridian is sent to the identity.''
By perturbing, we will replace this equation  by  a condition of the
form ``the
meridian and  longitude are related by the equation \ref{tnoia}.''
Hence, the perturbed flat moduli space can also be
identified as a subset of $R(Z,SU(3))$.
The results on the  local and global structure of the representation
variety
$R(Z,SU(3))$ that are developed in this section will therefore
be essential to our understanding of the behavior of the
moduli space under perturbation.

   \subsection{Cohomology calculations}
\label{CC-Z}
Let $Z$ be the complement
of the singular $r$-fiber in $\Si(p,q,r).$
In contrast to the homology sphere case, the abelianization of
$\pi_1(Z)$ is nontrivial.
Consequently, $\pi_1(Z)$ admits nontrivial abelian representations,
and reducible representations of $\pi_1(Z)$
do not always reduce  to
$SU(2)\times \{1\}$.  Given
a  representation $\al\colon \pi_1(Z) \to SU(3)$,
there are three possibilities:
\begin{itemize}
\item[(i)] $\al$ is irreducible,
\item[(ii)] $\al$ is nonabelian and reducible, or
\item[(iii)] $\al$ is abelian.
\end{itemize}

The first result is the analogue of Proposition \ref{coho}
for the knot complement $Z$.
\begin{prop} \label{coho-Z}
Suppose  $\al\colon \pi_1(Z) \to SU(3)$  is a nonabelian representation.
\begin{enumerate}
\item[(i)] If $\al$ is irreducible, then
   $\al(h) = e^{2 \pi i k/3} \cdot I $,
   $H^0(Z;su(3)_\al)  =0,$ and
$$H^1(Z;su(3)_\al)  =
\begin{cases} \RR^4 & \text{if $\al(x)$ and $\al(y)$
have three distinct eigenvalues,} \\
\RR^2 & \text{otherwise.}
\end{cases} $$
\item[(ii)] If $\al $ is reducible and has been
conjugated to take values in $S(U(2) \times U(1)),$
then
$$\al(h) = \left[ \begin{array}{ccc}
   e^{iv} &0&0 \\
    0    & e^{   iv}&0 \\
   0 &0&e^{- 2 iv}
   \end{array} \right].$$
With respect to the splitting $su(3)  = s(u(2)\times u(1)) \oplus \CC^2$
(see equation (\ref{dec2})), we have
$H^0(Z; s(u(2)\times u(1))_\al)   = \RR, \ H^0(Z;\CC^2_\al)=0, \
H^1(Z;s(u(2)\times u(1))_\al)  = \RR^2,$ and
   $$H^1(Z;\CC^2_\al) =
\begin{cases} \CC^2 & \text{if $\al(h)$ is central, i.e., if $e^{3iv}=1$,}
\\
0 & \text{otherwise.}
\end{cases}  $$
\end{enumerate}
\end{prop}

\begin{proof}
   This proof is similar to that given for
   Proposition \ref{coho}. We leave the details
   as an exercise for the reader. See also the proof of Lemma \ref{pl}.
\end{proof}

\begin{prop} \label{yum}
Suppose  $\al\colon \pi_1(Z) \to SU(3)$  is a nonabelian representation.
\begin{enumerate}
\item[(i)] If $\al$ is irreducible, then
   $$H^1(Z,\partial Z;su(3)_\al)  =
   \begin{cases} \RR^4 & \text{if $\al(x)$ and $\al(y)$
each have three distinct eigenvalues,} \\
\RR^2 & \text{otherwise.}
\end{cases} $$
\item[(ii)] If $\al $ is reducible and has been
conjugated to take values in $S(U(2) \times U(1)),$
then with respect to the
splitting  $su(3)  = s(u(2)\times u(1)) \oplus \CC^2$,
we have $H^1(Z, \partial Z;s(u(2)\times u(1))_\al)=\RR$
and
$$  H^1(Z,\partial Z;\CC^2_\al)  =
\begin{cases} \CC^2 & \text{if $\al(h)$ is central,} \\
0 & \text{otherwise.}
\end{cases}  $$
The map $H^1(Z,\partial Z;\CC^2_\al)\to H^1(Z;\CC^2_\al)$ induced by inclusion
is an isomorphism.
\end{enumerate}
\end{prop}

\begin{proof}
Associated to
$$(\partial Z, \varnothing)  \hookrightarrow (Z,\varnothing)
\hookrightarrow (Z,\partial Z)$$
is the
long exact sequence in cohomology
\begin{equation}\label{eqstar} \qquad \qquad  \cdots \to H^i(Z,\partial
Z) \to H^i(Z)
\to H^i(\partial Z) \to H^{i+1}(Z,\partial Z) \to \cdots. \qquad
\qquad\end{equation}
(We are temporarily omitting the
coefficients from the notation.)
To prove part (i), consider this sequence (\ref{eqstar}) with
coefficients
$su(3)_\al.$  The previous proposition shows that
$H^0(Z;su(3)_\al) =0$.
Additionally, since $\partial Z$ is
a 2-torus, Poincar\'{e} duality shows that
$H^0(\partial Z;su(3)_\al) =\RR^n = H^2(\partial Z;su(3)_\al) $
and $H^1(\partial Z;su(3)_\al) = \RR^{2n}.$ (In fact,
   $n = 4$ or $2$ depending on whether $\al(\mu)$ has
a double eigenvalue or not, but this has no bearing on the
rest of the argument.)

The nondegenerate pairing between relative and
absolute cohomologies gives
   $$\dim H^i(Z,\partial Z;su(3)_\al) = \dim H^{3-i}(Z;su(3)_\al).$$
   (E.g., $ H^3(Z;\partial Z;su(3)_\al)=0$.)
   The long exact sequence (\ref{eqstar}) with coefficients $su(3)_\al$
   has only seven nontrivial terms. Any long exact
   sequence has Euler characteristic zero, and
   so
   $$\dim H^1(Z,\partial Z;su(3)_\al) = \dim H^1(Z;su(3)_\al).$$

The proof of part (ii) is similar; in fact, for the coefficients
$\CC^2_\al$,
it is simplified by the observation that $H^i(\partial Z;\CC^2_\al)=0$,
and hence $ H^1(Z;\CC^2_\al) = H^1(Z,\partial Z; \CC^2_\al)$
as claimed. The long exact sequence
(\ref{eqstar}) with coefficients $s(u(2)\times u(1))_\al$
has nine nontrivial terms,
starting with $H^0(Z; s(u(2)\times u(1))_\al)$
which equals $\RR$ by the previous
proposition, and ending with $ H^3(Z, \partial Z; s(u(2)\times
u(1))_\al) $
which also equals $\RR$ by the nondegenerate pairing.
Arguing as before, it is not hard to see that
   $$\dim H^1(Z,\partial Z; s(u(2)\times u(1))_\al) =
   \dim H^1(Z; s(u(2)\times u(1))_\al)-1=2-1=1.$$
\end{proof}

We now turn our attention to the cohomology of the  abelian
representations of $\pi_1(Z).$
   To warm up, consider a nontrivial representation
   $\al\colon \pi_1(Z) \to U(1)$. The following lemma computes
    $H^0(Z;\CC_\al)$ and  $H^1(Z;\CC_\al)$.

\begin{lem} \label{pl}
Suppose $\al\colon \pi_1(Z) \to U(1)$ is a nontrivial representation.
Then
$H^0(Z;\CC_\al)=0$
and
$$H^1(Z;\CC_\al)=\begin{cases}
\CC& \text{if  $\al(\mu)^{pq}=1$,   $ \al(\mu)^{ap} \neq 1$ and $
\al(\mu)^{aq}  \neq 1$,}\\
   0&\text{otherwise.} \end{cases}
$$
\end{lem}

\begin{proof} We can compute the first two
cohomology groups of $Z$ using the cellular cohomology
of the 2-complex $Z_2$ determined by the presentation
of
$\pi_1(Z)$. The group presentation determines a
cellular structure  with one 0-cell, three 1-cells,
and four 2-cells.
The
   differentials in the cellular cochain
complex for the universal cover of $Z_2$ are
$$d^0=\left[ \begin{array}{c}x-1\cr y-1\cr h-1 \end{array}\right]
d^1=\left[ \begin{array}{ccc}1-h&0&x-1\cr 0&1-h&y-1\cr
1+x+\cdots + x^{p-1}&0&-x^p(h^{-1}+\cdots +h^{-a})\cr
0&1+y+\cdots +y^{q-1}&-y^q(h^{-1}+\cdots+h^{-a})
\cr\end{array}\right].$$

Taking the tensor product with $\CC$ over the
representation $\al$ has the effect of replacing $x,y,$ and $h$
in the matrices $d^0$ and $d^1$ by $\al(x), \al(y)$ and $\al(h)$.
We denote the resulting matrices by $d^0_\al$ and $d^1_\al$, and so
the cochain complex
$C^*(Z_2;\CC_\al)$ has the form
$$0\to\CC\,\mapright{d^0_\al} \;\CC^3
\; \mapright{d^1_\al} \; \CC^4 \to 0$$

Since $\al$ is nontrivial, it follows that
$d^0_\al$ is injective, hence $H^0(Z;\CC_\al)=0$
and $H^1(Z;\CC_\al)$ has dimension $3-$rank$(d^1_\al)-1$.

Notice that $[x]=aq  [\mu], \ y=ap [ \mu], $
and $h=pq [\mu] $ in homology.
If $\al(\mu)^{pq}\ne 1$, then $\al(h)-1\ne 0$, and so the
rank of $d^1_\al$ is at least two, and hence is exactly two. Thus
$H^1(Z;\CC_\al)=0$ if $\al(\mu)^{pq}\ne 1$.
This leaves the case when $\al(\mu)^{pq}=1$, i.e., when $\al(h)=1$.
In this case,  $\al(x)=\al(\mu)^{aq}$ is a $p$-th root of unity
and $\al(y)= \al(\mu)^{ap}$ is a $q$-th root of unity. Thus
$$d^1_\al= \left[ \begin{array}{ccc}0&0&\al(x)-1\\
0&0&\al(y)-1\\
1+\al(x)+\cdots + \al(x)^{p-1}&0&-a  \cr
0&1+\al(y)+\cdots +\al(y)^{q-1}&-a \cr\end{array}\right].$$
This matrix has rank $2$ unless $\al(x)$ is
a  nontrivial $p$-th root of unity and $\al(y)$ is a
nontrivial $q$-th root of unity, in which case it has rank 1. The
lemma follows.
\end{proof}

Now consider abelian representations $\al\colon \pi_1(Z) \to SU(3)$.
By conjugation, we can assume that $\al$ takes values
in the maximal torus $T \subset SU(3)$.
Under the adjoint action of $T$,  the Lie algebra
$su(3)$ decomposes as
\begin{equation}\label{lasplit}su(3) = \CC^3 \oplus
\RR^2.\end{equation}The $\CC^3$ corresponds to the off-diagonal entries
and $\RR^2$ to the diagonal entries.
Then $T$ acts trivially on $\RR^2$
and by rotations on each of the three complex factors.
More precisely, the action on $\CC^3$ is given by
$$ \left[ \begin{array}{ccc} \om_1 &0&0\\
0& \om_2 &0\\ 0&0& \bar{\om}_1 \bar{\om}_2 \end{array} \right]
\cdot \left[ \begin{array}{c} z_1\\ z_2 \\ z_3  \end{array} \right]
= \left[ \begin{array}{c} \om_1 \bar{\om}_2 \, z_1\\
\om_1^2 \om_2 \, z_2\\
\om_1 \om_2^2 \, z_3  \end{array} \right].$$

An abelian representation $\al\colon \pi_1(Z) \to SU(3)$
is completely determined
by $\al(\mu)$, since $H_1(Z;\ZZ)$ is generated by $[\mu]$.
Suppose in addition that $\al$
is the limit
of a sequence of $SU(2)\times \{1\}$ representations.
Then
we can arrange that
\begin{equation} \label{diagon2}
\al(\mu) =  \left[ \begin{array}{ccc} \om &0&0\\
0& \bar{\om} &0\\ 0&0& 1 \end{array} \right].
\end{equation}
In this case, there is a distinguished $\CC^2$ subbundle
of the adjoint bundle $Z \times su(3)$ on which $\al(\mu)$ acts by
$(z_1,z_2) \mapsto (\om z_1, \bar{\om} z_2)$
(namely the last two coordinates in $\CC^3$).
Suppose further that  $\al$ is nontrivial.
Then  $H^0(Z;\CC^2_\om) =0$.
Applying Lemma \ref{pl}
to $ \CC^2_\al = \CC_{\om} \oplus \CC_{\bar{\om}}$,
and noting that $H^*(X;\CC_\om) \cong H^*(X;\CC_{\bar{\om}}),$
we  see that
$$H^1(Z;\CC^2_\al)
= \begin{cases}
\CC^2 & \text{if $\om^{pq}=1$ and $\om^{ap} \neq 1 \neq \om^{aq}$,} \\
0 & \text{otherwise.}
\end{cases}$$
The next proposition extends these computations to
abelian representations in a neighborhood of $\al$.

\begin{prop} \label{abelians}
Let $\al\colon \pi_1(Z) \to SU(3)$ be a fixed nontrivial, abelian
representation with $\al(\mu)$
given by the diagonal matrix in equation (\ref{diagon2}).
Suppose further that $\om^{pq}=1$
and $\om^{ap} \neq 1 \neq \om^{aq}.$
(Thus $H^1(Z;\CC^2_\al)  = \CC^2$.)
Consider abelian representations $\be \colon \pi_1(Z) \to SU(3)$
near to but distinct from $\al$.
   Conjugating, we can arrange that
   $$\be(\mu)  =\left[ \begin{array}{ccc} \om_1 &0&0\\
0& \om_2 &0\\ 0&0& \bar{\om}_1 \bar{\om}_2 \end{array} \right]$$
with $\om_1$ close to $\om$ and $\om_2$ close to $\bar{\om}$
(so $\om_1 \om_2$ is close to 1).
Then, for $\be$ close enough to
   $\al$, we have
   $H^0(Z;\CC^2_\be)=0$ and
$$H^1(Z;\CC^2_\be)=H^1(Z,\partial Z;\CC^2_\be)=
   \begin{cases} \CC & \text{if $(\om^2_1  \om_2)^{pq}=1$ or if $(\om_1
\om^2_2)^{pq}=1$,} \\
0& \text{otherwise.}
\end{cases}$$
\end{prop}

\begin{proof}
That $H^0(Z;\CC^2_\be) = 0$
follows from  upper semicontinuity of $\dim H^0$ on the
representation variety.
The computation of $H^1 (Z;\CC^2_\be)$ follows
from Lemma \ref{pl}, keeping in mind that our hypotheses
exclude the possibility $\be = \al$.
 All that remains is to prove the claim
about relative cohomology.
   Set $T= \partial Z$. If $\ga\colon \pi_1(T)\to SU(2)$ is any {\it
nontrivial} representation, then
$H^*(T;\CC^2_\ga)=0$ (cf.~equation (3.4) of \cite{BHKK}).
Now using the long exact sequence in cohomology, it follows
that $H^1(Z;\CC^2_\be)=H^1(Z,\partial Z;\CC_\be)$
for $\be$
in a small enough neighborhood of $\al$.
\end{proof}

\subsection{The representation variety R(Z,SU(3))}
\label{dot}

In this subsection, we give a description of $R(Z,SU(3))$.
This space is the union of three different strata:
\begin{enumerate}
\item[(i)] $R^*(Z,SU(3)),$ the stratum of
irreducible representations.
\item[(ii)] $R^\red(Z,SU(3))$,
the stratum of reducible, nonabelian representations.
\item[(iii)] $R^\ab(Z,SU(3)),$ the stratum of
abelian representations.
\end{enumerate}

We will  describe each of these  strata
presently. For $R^*(Z,SU(3))$, this involves certain
  double coset spaces,
  and for $R^\red(Z,SU(3))$,
  this builds on the results in \cite{K}.
  Note that, given any  finitely presented group $\pi$, two nonabelian
representations $\al_0,\al_1\colon\pi\to S(U(2)\times U(1))$ are conjugate
in $SU(3)$ if and only if they are conjugate by  a matrix in
$SU(2)\times \{1\}$. In particular, the natural map
 $R^*(Z, S(U(2)\times U(1)))\to
R(Z, SU(3))$ is injective  and has image in  $R^\red(Z, SU(3))$.

We begin with the description of $R^\ab(Z, SU(3))$  because it
is the simplest.
Since the homology class of the meridian
$\mu$   generates $H_1(Z;\ZZ)$,
a conjugacy class $[\al]$ of abelian representations
is completely determined by the conjugacy class
of  $\al(\mu)$.
Thus,   $R^\ab(Z,SU(3))$
  is parameterized by
the quotient $SU(3)/\hbox{conj},$ which is just
the quotient $T/S_3$ of the maximal torus by  the
Weyl group. This is parameterized by the
standard 2-simplex $\De$,
see
equation (\ref{tri}) in Subsection \ref{sfuse}.

For the stratum $R^\red(Z,SU(3))$, note that
every reducible representation can be conjugated to
have image in $S(U(2) \times U(1)).$ We will see that
every $S(U(2) \times U(1))$
representation of $\pi_1(Z)$ is obtained by twisting an
$SU(2)$ representation, and we will combine this observation with
an explicit description of the $SU(2)$ representation varieties
of $\pi_1(Z)$ (essentially from Klassen's work \cite{K})
to prove that $R^\red(Z,SU(3))$
is a union of $(p-1)(q-1)/4$ open 2-dimensional cylinders
under the assumption that $p,q$ are both odd (see Proposition
\ref{seam}).

Let $\al\colon \pi_1(Z)\to SU(3)$ be a nontrivial reducible
representation sending $(xy)^r h^c$ to the identity.
Thus $\al$ extends over the solid torus and gives
a reducible representation $\pi_1(\Si) \to SU(3).$
In particular, $\al$ reduces to $SU(2)\times \{1\}$
and is nonabelian.

In Proposition \ref{coho-Z}, we computed that
$ H^1(Z;s(u(2)\times u(1))_\al) = \RR^2,$
hence  the reducible stratum $R^\red(Z,SU(3))$
has 2-dimensional Zariski tangent
space at $[\al]$.
In this subsection, we construct an explicit
2-parameter family
of reducible representations $\al_{s,t}\colon \pi_1(Z) \to SU(3)$ near
$\al$, showing  that all the Zariski
tangent vectors are integrable.
  From this, we will conclude that the reducible
stratum $R^\red(Z,SU(3))$ is smooth and 2-dimensional near $[\al].$

The 2-parameter family will be obtained by twisting
$SU(2)\times \{1\}$ representations of $\pi_1(Z)$
   to  representations with image in $S(U(2)\times U(1))$.
   To get started, we describe the
$SU(2)$
representation variety of $\pi_1(Z)$.
Note that
we have assumed that $p$ and $q$  are both odd
in the presentation (\ref{pres-Z}).
The following result is proved by techniques developed by  Klassen in
\cite{K}. The methods he uses to describe $SU(2)$
representation varieties of complements of torus knots
work equally well for the 3-manifolds $Z$ considered
here.
We view $SU(2)$ as the unit quaternions and write a
typical element
as $a + i b + j c + k d$ where $a,b,c,d \in \RR$ satisfy
$a^2+b^2+c^2+d^2=1.$

\begin{prop}  \label{pad}
$R^*(Z,SU(2))$ consists of  $(p-1)(q-1)/2$ open arcs of
irreducible representations.
These arcs are given  as follows. For each   $k\in\{1,\cdots, p-1\}$,
$\ell\in\{1,\cdots ,q-1\}, \varep \in\{0,1\}$
satisfying $k\equiv \ell\equiv a \varep \pmod{2}$,
the assignment to   $s\in [0,1]:$
  \begin{eqnarray*}
\be_s(x) &=& \cos(\pi k/p) + i \sin(\pi k/p),  \\
\be_s(y) &=&  \cos(\pi \ell/q) + \sin(\pi \ell/q)(i \cos(\pi s) + j
\sin(\pi s))\\
  \be_s(h)&=&(-1)^{\varep}
   \end{eqnarray*}
  defines a path of $SU(2)$ representations which are irreducible for
$s\in (0,1)$.
Moreover,  for $s\in(0,1)$,
$$  H^1(Z;\CC^2_{\be_s}) =
\begin{cases} \CC^2 & \text{if $\varep =0$, i.e., if $\be_s(h)=1$,}\\
0 & \text{if $\varep=1$, i.e., if $\be_s(h)=- 1$.} \end{cases}$$

The two limit points of each open arc, $\be_0$ and $\be_1$,
are   abelian representations sending $\mu$
to $(-1)^k e^{\pi i\left(\frac{r(kq+\ell p)}{pq}\right)}$
and $ (-1)^k e^{\pi i\left(\frac{r(kq-\ell p)}{pq}\right)}$.\qed
\end{prop}
(The cohomology calculation in Proposition \ref{pad} follows from
Proposition \ref{coho-Z}.)

To summarize, the subspace of $R(Z,SU(3))$ consisting of conjugacy
classes of nonabelian $SU(2)\times \{1\}$
representations of $\pi_1(Z)$ is a union of $(p-1)(q-1)/2$ open arcs
with ends that limit to points in the abelian stratum.
The intersection of the subspace $R(\Sigma, SU(3))\subset
R(Z,SU(3))$ with such an arc of reducible representations
  consists of either   reducible representations on   pointed
2-spheres or   isolated reducible representations (i.e.~ Type
Ib representations),   depending on whether or not $h$ is sent to $I$.

In defining these 1-parameter families of representations,
we arranged  that $x$ was sent to a
diagonal matrix.  For future applications, it is
convenient to arrange (by conjugation) that
$xy$ is sent to a diagonal matrix, because then it follows
from  equations (\ref{mu}) that the
meridian and longitude will also be diagonal.

   Fix a connected component   of $R^*(Z,SU(2))$ determined by the
triple $(k,\ell,\varep)$  with $k\equiv\ell\equiv a\varep\pmod{2}$ as above,
and denote by $\al_s$ the corresponding arc of $SU(2)\times \{1\}$
representations sending $xy$ to a diagonal matrix.
A short calculation shows that
$$\al_s(xy) = \left[ \begin{array}{ccc} e^{iu} & 0 &0 \\
0 & e^{-iu} & 0 \\
0 & 0 & 1\end{array}\right]$$
where $u$ satisfies the equation
\begin{equation} \label{coseq}
\cos(u) =\cos(\pi k/p)\cos(\pi \ell/q) -
\sin(\pi k/p) \sin(\pi \ell/q) \cos (\pi s).
\end{equation}
   We next show that   the arc $[ \al_s]$ of
$SU(2)\times\{1\}$-representations
   is  a codimension one subset of
   $R^\red(Z, SU(3))$.
The other degree of freedom
comes from {\em twisting} a representation out of $SU(2)\times \{1\}$,
keeping it in $S(U(2)\times U(1))$.

First, given
$$A=\left[\begin{array}{cc} a & b  \\
-\bar{b} & \bar{a}   \end{array} \right] \in SU(2),$$
the {\em twist} of $A$ by $e^{i\th} \in U(1)$ is
the $S(U(2) \times U(1))$ matrix
$$
\left[\begin{array}{ccc} e^{i\th}   &0 &0\\
0  \ & e^{i\th}   & 0 \\
0&0& e^{-2 i\th}\end{array} \right]
\left[\begin{array}{ccc}  a &  b &0\\
 - \bar{b}  &  \bar{a} & 0 \\
0&0& 1\end{array} \right]
=\left[\begin{array}{ccc} e^{i\th} a & e^{i\th} b &0\\
- e^{i\th} \bar{b} \;\; \ & e^{i\th} \bar{a} & 0 \\
0&0& e^{-2 i\th}\end{array} \right].
 $$
The map  $SU(2) \times U(1) \to S(U(2) \times U(1))$
defined by twisting is a 2-to-1 map.
   In terms of $U(2)$, this is simply the description
   $U(2) = SU(2) \times_{\ZZ_2} U(1),$
   and twisting is just scalar
multiplication by $e^{i \th}$.
Notice that the matrix $\Phi(u,v)$ appearing in equation (\ref{defnPhi})
is  the twist of the diagonal $SU(2)$ matrix $A$ with entries $e^{iu}, e^{-iu}$ by $e^{iv}$.

Suppose $\chi\colon \pi_1(Z)\to U(1)$ is a character, i.e.~ a
homomorphism into the abelian group $U(1)$, and let $\be\colon
\pi_1(Z)\to SU(2)$ be a representation.
The reducible $SU(3)$
representation obtained by twisting $\be$ by $\chi$ is
defined to be representation $\pi_1(Z) \to S(U(2)\times U(1))$
taking an element $w \in
\pi_1(Z)$ to the twist of $\be(w)$ by $\chi(w)$.
Notice that, since $H_1(Z;\ZZ) \cong \ZZ$ is   generated by the
meridian $\mu,$ any character $\chi$ is completely
determined by the element $\chi(\mu)\in U(1),$ which can be arbitrary.
If $\chi(\mu) = -1,$ then
the twist of $\be$ by $\chi$ is again an $SU(2)$ representation,
and twisting by this central character
defines an  involution on the $SU(2)$
representation variety of knot complements.

We give a more explicit description of the
stratum $R^\red(Z,SU(3))$ of reducible $SU(3)$ representations in terms
of
twisting  the arcs $\al_s$ described above.
\begin{defn}\label{defnofalsth}
Fix $e^{i \th} \in U(1)$ and let $\chi_\th$ be the character
sending $\mu$ to $e^{i \th}$. Let $\al_s$ be representation described
in     Proposition \ref{pad}  corresponding to a triple $(k,\ell,\varep)$
and  $s\in(0,1)$.
Define the reducible  $SU(3)$ representation  $ \al_{s,\th}\colon
\pi_1(Z) \to S(U(2)\times U(1))\subset SU(3)$ to be  the twist of $
\al_{s}$ by $\chi_\th$.\end{defn}

\begin{prop}\label{rolm} Fix $(k,\ell,\varep)$ with $k\equiv\ell\equiv
a\varep\pmod{2}$ as in Proposition  \ref{pad}  and let $\al_{s,\th}$ be
the 2-parameter family of $S(U(2)\times U(1))$ representations
corresponding to twisting $\al_s$ by $\th$. Then
the representation $\al_{s,\th}$  sends
$x$ to the twist of $ \al_{s}(x)$ by $e^{i a q \th},$
$y$ to the twist of $ \al_{s}(y)$  by $e^{i a p \th},$
and $h$ to the twist of $ \al_{s}(h)$  by $e^{i pq \th}.$
Moreover,
$$\al_{s, \th}(\mu)=
   \left[ \begin{array}{ccc} (-1)^{kc} e^{i(\th+ru)} & 0 &0 \\
0 & (-1)^{kc} e^{i(\th-ru)} & 0 \\
0 & 0 & e^{-2i\th}\end{array}\right]$$
and
$$\al_{s, \th}(\la)= \al_s(\la) =
   \left[ \begin{array}{ccc} (-1)^{ka(p+q)} e^{ipqu} & 0 &0 \\
0 & (-1)^{ka(p+q)} e^{-ipqu} & 0 \\
0 & 0 & 1\end{array}\right],$$
where $u$ satisfies equation (\ref{coseq}).
The representation $\al_{s,\th}$ is conjugate to an $SU(2)\times\{1\}$
representation only for $\th\in\pi\ZZ$, and the arcs $\al_{s,0}$ and
$\al_{s,\pi}$   are different components of $R(Z,SU(2))$. The map
$(s,\th)\mapsto  \al_{s,\th}$ defines a smooth
2-dimensional subvariety of $R(Z,SU(3))$ contained in $R^\red(Z,SU(3))$
and homeomorphic to $(0,1)\times
S^1$.
\end{prop}

\begin{proof}
The first few assertions  follow immediately from the definitions and
equations (\ref{mu}) and (\ref{homolog}).

By taking the determinant of $e^{i\th}\al_s$, it is easy to check that
$\al_{s,\th}$ is an $SU(2)\times\{1\}$ representation if and only if
$\th\in \pi \ZZ$. The representation $\al_{s,0}$ takes $h$ to the
diagonal matrix with
entries $(-1)^\varep,(-1)^\varep, 1$ and $\al_{s,\pi}$ takes $h$ to the
diagonal matrix with
entries $(-1)^{pq+\varep},(-1)^{pq+\varep}, 1$.
Since $p$ and $q$ are both odd,   $\al_{s,0}$ and $\al_{s,\pi}$ are
different arcs. The map $(s,\th)\mapsto [\al_{s,\th}]\in R(Z,SU(3))$ is
injective, and since
$H^1(Z;s(u(2)\times u(1)))_{\al_{s,\th}})=\RR^2$ by Proposition
\ref{coho-Z}, this parameterizes a smooth subvariety.
\end{proof}

Every representation $\al$ in $R^\red(Z,SU(3))$ is conjugate to some
$\al_{s,\th}$ for some choice of $(k,\ell,\varep)$ and $(s,\th)$. The
reason for this is that one can first conjugate  $\al$ into
$S(U(2)\times U(1))$,
and then if the $(3,3)$ entry of $\al(\mu)$ is $e^{2i\th}$, $\al$ must
be the $\th$-twist of some $SU(2)$ representation $\al_s$.

By Proposition \ref{pad}, it follows that  $R^\red(Z,SU(3))$
  has exactly $(p-1)(q-1)/4$ components, each of which is a  smooth open
  cylinder with two seams of $SU(2)\times\{1\}$
  representations  (see
Figure 2).

\begin{figure}[t]  
\begin{center}
\leavevmode\hbox{}
\includegraphics[width=3in]{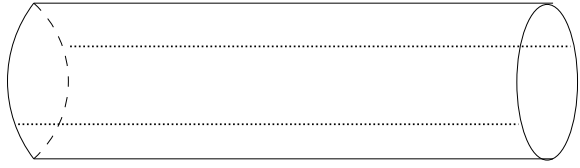}
\caption{An open cylinder of reducible $SU(3)$ representations with two
  seams  (the dotted lines) representing the $SU(2)\times \{1\}$
representations.}
\end{center} 
\end{figure}

\begin{remark} It is not hard to extend Proposition \ref{rolm} to the
case when either $p$ or $q$ is even. In that case it is possible for $\al_{s,0}$
and $\al_{s,\pi}$ to represent the same arc, with  a reversal of
orientation, and hence there are  components of $R^\red(Z, SU(3))$
homeomorphic to an open M\"obius band.  For example, the space of
nonabelian reducible $SU(3)$ representations of the Trefoil  knot
complement is  a single M\"obius band.\end{remark}

The following theorem  summarizes our discussion.

  \begin{thm} \label{seam}
  Suppose $\Si(p,q,r)$ is a Brieskorn sphere and reorder
  $p,q,r$ so that $p$ and $q$ are both odd. Let $Z$ be the complement
  of the singular $r$-fiber of $\Si(p,q,r).$ Then the
  stratum $R^\red(Z,SU(3))$ of conjugacy classes
  of nonabelian reducible representations
  is a smooth, open, 2-dimensional manifold consisting of
  $(p-1)(q-1)/4$ path components, each of which is diffeomorphic to
the open cylinder $(0,1)\times S^1$. 
The closure of such a component  in $R(Z, SU(3))$ contains
two boundary circles, which are circles immersed in the abelian stratum $R^\ab(Z,SU(3))$
with isolated double points.   \end{thm}

Fix $(k,\ell,\varep)$ with $k\equiv\ell\equiv a\varep\pmod{2}$ as in
Proposition  \ref{pad} and let $\al_{s,\th}\colon \pi_1(Z)\to S(U(2)\times
U(1))$ denote the corresponding 2-parameter family of representations.
Suppose  for some $s$, $\al_{s,0}$ extends to a reducible representation
on $\pi_1(\Si).$ This is the case if and only if
   $ \al _{s,0}(\mu) = I$, namely if $  \al_{s,0}(xy)$
is an $r$-th root of $ \al_{s,0}(h^c).$

Since $H^1(\Si; su(2)_{\al_{s,0}})=0$ and $\Si$ is a homology sphere,
none  of the nearby
representations in the
2-parameter family
$\al_{s, \th}$ of $\pi_1(Z)$ extend to representations of $\pi_1(\Si).$

If $[ \al_{s,0}]$ lies on a  2-sphere component of
$R(\Si,SU(3)),$ then  $H^1(\Si;\CC^2_\al) \neq 0$ and
$ \al_{s,0}(h)=I$ (i.e. $\varep=0$).
Hence $\al_{s,0}(xy)$ is an $r$-th root of $I$
and $s$ satisfies the equation
$$ \cos\left(\tfrac{2 \pi m}{r}\right) = \cos(\pi k/p)\cos(\pi \ell/q) -
\sin(\pi k/p) \sin(\pi \ell/q) \cos (\pi s)$$
for some $0<m<r$. In particular,
\begin{equation} \label{rroott}
\al_{s,0}(xy) =\left[ \begin{array}{ccc} e^{2 \pi i m/r} & 0 &0 \\
0 & e^{-2 \pi i m/r} & 0 \\
0 & 0 & 1\end{array}\right].
\end{equation}

We now consider irreducible representations
$\al\colon \pi_1(Z) \to SU(3)$ and
  give a description of the
closure of $R^*(Z,SU(3))$.
We begin with a simple observation. If
$\al\colon \pi_1(Z)\to  SU(3)$ is an irreducible representation, then
$\al(h)$ lies in the center of $SU(3)$ and it follows
from the presentation
(\ref{pres-Z}) that
$\al(x)^{p}=\al(y)^{q}=\al(h)^a$. Conversely, suppose
we are given matrices $A,B, H\in SU(3)$ with $H$ central such that
\begin{equation} \label{ABHeq}A^{p}=B^{q}=H^a,
\end{equation}
then setting
$\al(x)= A, \al(y) = B,$ and $\al(h)=H$ uniquely determines a
representation $\al\colon \pi_1(Z) \to SU(3)$. This representation is
reducible  if and only if $A$ and $B$ share an eigenspace.

For $A,B,H$ diagonal $SU(3)$ matrices with $H$ central and satisfying equation
(\ref{ABHeq}),
we can write $H = e^{2 \pi i \ell /3} I $ for a unique $\ell \in
\{0,1,2\}$ and
we denote by
$\cC_{AB}^\ell \subset R(Z,SU(3))$
the subset of conjugacy classes $[\al]$
of representations with
$\al(x)$ conjugate to $A$, $\al(y)$ conjugate to $B$, and $\al(h) =e^{2
\pi i \ell/3}I$.
There is a map $\Psi\colon  SU(3) \to \cC_{AB}^\ell$ where $\Psi(g) =
[\psi_g]$
is the conjugacy class of the representation  $\psi_g$ with
  $\psi_g(x) = A$ and $\psi_g(y) = g B g^{-1}.$
Let $\Ga_A$ and $\Ga_B$ denote
the stabilizer subgroups of $A$ and $B$.
If $\ga \in \Ga_B$, then
$\psi_{g\ga} = \psi_g$ for all $g \in SU(3)$.
Likewise, if $\ga \in \Ga_A,$ then
$\psi_{\ga g} = \ga \psi_{g} \ga^{-1}$ for all $g \in SU(3)$. Thus,
$\Psi$ factors
through left
multiplication by $\Ga_A$
and right multiplication by $\Ga_B$ and determines a map from
the double coset space
$$\Psi\colon \Ga_A\backslash SU(3)/ \Ga_B \to \cC_{AB}^\ell$$
which is a homeomorphism which is smooth on the stratum
of principal orbits.

Elementary dimension counting gives that $\cC_{AB}^\ell$
has dimension four if both $A$ and $B$ have three distinct
eigenvalues and dimension two if exactly one of $A$ or $B$ has a
2-dimensional eigenspace. In all other cases, $\cC_{AB}^\ell$
does not contain any irreducibles. For example, if both
$A$ and $B$ have double eigenspaces, then the eigenspaces
intersect nontrivially in an invariant linear subspace, giving a reduction. Similarly,
if either $A$ or $B$ has an eigenvalue of multiplicity three, then
the corresponding representation is necessarily abelian.

Observe further that the set
$\cC_{AB}^\ell$ depends only on $\ell \in \{0,1,2\}$ and the conjugacy
classes of the matrices $A$ and $B$. Thus, we can assume
without loss of generality that $A$ and $B$ are both
diagonal.
\begin{thm} \label{Zsthm}
The closure of the stratum
$R^*(Z,SU(3))$ of irreducible representations is a union
$\bigcup \cC_{AB}^\ell$, where the  union is over
pairs $([A],[B]) \in \left( SU(3)/\hbox{conj}\right) ^2$
 and $\ell \in \{0,1,2\}$
satisfying the conditions:
\begin{enumerate}
\item[(i)] $A^{p}=B^{q} = H^a$, where $H= e^{2 \pi i \ell/3} I$,
\item[(ii)] neither $A$ nor $B$ is central, and
\item[(iii)] one of $A$ or $B$ has three distinct eigenvalues.
\end{enumerate}
In particular
\begin{enumerate}\item[$\bullet$] If one of $A$ or $B$ has a repeated
eigenvalue, then
$\cC_{AB}^\ell$ is 2-dimensional and is called a Type I component of
$R(Z,SU(3))$.
\item[$\bullet$] If both $A$ and $B$ have three distinct eigenvalues,
then
$\cC_{AB}^\ell$ is 4-dimensional and is called a Type II component of
$R(Z,SU(3))$.\end{enumerate}
\end{thm}

Given a nonabelian reducible representation $\al\colon \pi_1(Z) \to
SU(3)$,
we would like to know when there exists a  1-parameter family of
irreducible
representations limiting to $\al.$ If there is, then
Proposition \ref{coho-Z} implies
that $\al(h)$ is central.
The following proposition is a
partial converse.

\begin{prop}   \label{TFAE-Z}
If $\al\colon \pi_1(Z) \to SU(3)$ is a nonabelian reducible
representation
satisfying:
\begin{enumerate}
\item[(i)]  $\al(h)$ is central, and
\item[(ii)]  one of $\al(x)$ or $\al(y)$ has three distinct eigenvalues,
\end{enumerate}
then there exists a  1-parameter family of irreducible
$SU(3)$ representations limiting to $\al$.
\end{prop}
  \begin{remark}
  Notice that the condition $H^1(Z;\CC^2_\al) \neq 0,$
  which is equivalent to (i), is not enough to guarantee that there be a
  family of irreducible representations limiting to $\al$. There are
  nonabelian reducible representations with $\al(h)$ central such that
  $\al(x)$ and $\al(y)$ both have repeated eigenvalues.
  Such representations are not in the closure of
  $R^*(Z,SU(3))$ even though $H^1(Z;\CC^2_\al) \neq 0$.
  \end{remark}
\begin{proof}
Set $A = \al(x)$ and $B = \al(y).$
Notice that the assumption that $\al$ is nonabelian
implies that neither $A$ nor $B$ is central.
Obviously $[\al] \in \cC_{AB}^\ell.$
The subspace $\cC ^{\ell, \red}_{AB}$ of conjugacy
classes of reducible representations
has codimension greater than or equal to one,
and this completes the proof.
   \end{proof}
It is not hard to show that $\cC_{AB}^{\ell,\red}$ has dimension
one. We leave this as an exercise for the reader.  Note that 
$\cC^{\ell,\red} _{AB}$ is also a codimension one subset of $R^\red(Z,SU(3))$.
 The next lemma is a slight reformulation of
\cite[Lemma 2.4]{Ha}.
We include the proof
for the sake of completeness.
\begin{lem} \label{sbmrs}
Suppose $A, B \in SU(3)$ are diagonal matrices
and consider the
map $\varphi\colon SU(3) \lto \CC$ defined by setting $\varphi(g)=
\tr(AgBg^{-1}).$
Then, for fixed $g\in SU(3)$,
the differential $d \varphi_g$ is surjective provided
\begin{enumerate}
\item[(i)] $A$ and $gBg^{-1}$ have no common eigenvectors, and
\item[(ii)] the product $A gBg^{-1}$ has three distinct eigenvalues.
\end{enumerate}
Equivalently, $d\phi_g$ is surjective if $\psi_g \colon \pi_1 Z \to SU(3)$ is
irreducible and $\psi_g (xy)$ has 3 distinct eigenvalues.  
\end{lem}
\begin{proof}
Since condition (i) cannot hold when  $A$ and $B$ both have double
eigenspaces,
we assume (by switching the roles of $A$ and $B$,
if necessary) that the eigenvalues of $B$ are distinct.

There is an element $h \in SU(3)$
so that  $h (B g^{-1} A g) h^{-1}$
is diagonal. Let
$$C= h (B g^{-1} A g) h^{-1}=\left[  \begin{array}{ccc}
c_1 & 0 & 0 \\ 0 &c_2 &0 \\ 0&0& c_3
\end{array}\right]$$
be the resulting matrix.
Condition (ii) implies that $c_1,c_2,$ and $c_3$
are all distinct.

If $g_t \in SU(3)$ is a path passing through $g$ at $t=0,$
then
$g_t = g(I + tX   + O(t^2))$ for some $X \in su(3)$
and  $g_t^{-1} =(I - tX  + O(t^2))g^{-1}$.
Since $tr(A g B g^{-1})=tr(B g^{-1} A g)$, we have
\begin{eqnarray*}
d\varphi_g(X)&=& \left.  \tfrac{d}{dt} \tr(B g_t^{-1} A g_t )
\right|_{t=0} \\
&=& tr(B g^{-1} A g X - g^{-1}AgBX)\\
&=& tr(B g^{-1} A g \cdot (X-BXB^{-1})) \\
&=& tr(C \cdot h (X-BXB^{-1}) h^{-1}) \\
&=& (c_1 -c_3) (ir_1) + (c_2-c_3)(ir_2),
\end{eqnarray*}
where
$ir_1$ and $ir_2$ are the $(1,1)$ and $(2,2)$  entries of $h
(X-BXB^{-1}) h^{-1}$.
(Here,  $r_j \in \RR$ since $h (X-BXB^{-1}) h^{-1} \in su(3).$)

Write  $h=(h_{ij})$ and  let $X$ have the form
\begin{equation} \label{Xdefn}
X= \left[ \begin{array}{ccc}
0 & u & v \\ -\bar{u} &0 &w \\ -\bar{v} &-\bar{w}& 0
\end{array} \right],
\end{equation}
we compute that
\begin{eqnarray*}
r_1 &=& 2 \Im \left( h_{11} \bar{h}_{12} (1-b_1 b_2^{-1}) u +
h_{11} \bar{h}_{13} (1-b_1 b_3^{-1}) v +
h_{12} \bar{h}_{13} (1-b_2 b_3^{-1}) w\right), \\
r_2 &=& 2 \Im \left( h_{21} \bar{h}_{22} (1-b_1 b_2^{-1}) u +
h_{21} \bar{h}_{23} (1-b_1 b_3^{-1}) v +
h_{22} \bar{h}_{23} (1-b_2 b_3^{-1}) w\right),
\end{eqnarray*}
where $\Im(x+iy)=y$ is the imaginary part of a complex number.

Suppose that two of the entries of $(h_{ij})$
vanish.
Orthogonality of the rows and columns of $h$ then
implies that two of the other entries of $(h_{ij})$ also vanish.
Thus $h$ must send one of the standard basis vectors $e_j$ to
another (possibly different) standard basis vector, perhaps multiplied
by a unit complex number.
Therefore $e_j$ is
an eigenvector for both $B$ and $h^{-1} C h$,
and this contradicts condition (i).

Thus at most one of the entries of $(h_{ij})$
equals zero.
This implies that $r_1 \neq 0$  for some choice
of $u,v,$ and $w.$
For if $r_1=0$ for all $u,v,w \in \CC$, then
two of $\{ h_{11},h_{12}, h_{13}\}$  must vanish
(because $\{b_1, b_2, b_3 \}$ are all distinct).
Similarly   $r_2 \neq 0$ for some choice
of $u,v,$ and $w.$

Since at most one of the entries of $(h_{ij})$ vanishes,
one of the two cases holds:

\smallskip\noindent
{\sc Case 1:} each $\{h_{11},h_{12}, h_{13}\}$ is nonzero, or \\
{\sc Case 2:} each $\{h_{21},h_{22}, h_{23}\}$ is nonzero.

\smallskip

The proofs for the two cases are similar,
and we supply the details for Case 1 only, leaving the
rest of the argument to the reader.

Notice that the set of matrices of the form (\ref{Xdefn})
form a real vector space of dimension six.
We will apply
$d \varphi_g$ to a basis $\{X_1,\ldots, X_6\}$
and show that the image spans $\CC$ as a real vector
space.   It is useful to make the
simplifying substitutions:   $$u=\frac{u'}{h_{11} \bar{h}_{12}
(1-b_1b_2^{-1}) },
\quad v =\frac{v'}{h_{11} \bar{h}_{13}(1-b_1b_3^{-1})},
\quad \hbox{ and } w=\frac{w'}{h_{12} \bar{h}_{13} (1-b_2b_3^{-1})}.$$
Define six distinct matrices $X_1,\ldots, X_6 \in su(3)$ as in equation
(\ref{Xdefn}) as follows: for $X_1$ and $X_2$, take
$u' \in \{1,i\}$ and $v'=0=w'$; for $X_3$ and $X_4$, take
$v'\in \{1,i\}$  and $u'=0=w'$; and for $X_5$ and $X_6$, take
$w' \in \{1,i\}$ and $u'=0=v'$.
We claim the set
$$S= \{ d \varphi_g(X_1), \ldots , d \varphi_g(X_6)\}$$
spans $\CC$ as a 2-dimensional real vector space.
Suppose otherwise, namely suppose $S$ does not span $\CC$.
Condition (ii) implies that $c_1, c_2, c_3$
are all distinct, from which it follows that
$\frac{c_1-c_3}{c_2 -c_3} \not\in \RR$.
The only way   $S$ could be linearly dependent is if
$$\frac{ h_{21} \bar{h}_{22}}{h_{11} \bar{h}_{12}} =
\frac{ h_{21} \bar{h}_{23}}{h_{11} \bar{h}_{13}} =
\frac{{h}_{22} \bar{h}_{23}}{h_{12} \bar{h}_{13}} \in \RR.$$
Taken one at a time, we obtain the three equations:
$$h_{11} h_{22} = h_{21} h_{12}, \quad
h_{11} h_{23} = h_{21} h_{13}, \quad
h_{12} h_{23} = h_{22} h_{13}.$$
Expanding along the bottom row of $(h_{ij})$,
these equations imply that  $\det(h)=0$,
which contradicts the fact that $h \in SU(3)$
and completes the proof.
\end{proof}

Now suppose $A,B,\ell$ satisfy the hypotheses
of Theorem \ref{Zsthm}. If we define $\phi\colon \cC_{AB}^\ell \to \CC$
by setting  $\phi([\al]) = \tr(\al(x) \al(y))$, then
the following triangle commutes:

\[\begin{diagram}\dgmag{1000}\dgsquash
\node{SU(3)}\arrow{s,l}{\Psi}\arrow{se,t}{\varphi}\\
\node{\cC_{AB}^\ell}\arrow{e,b}{\phi}\node{\CC}\end{diagram}
\]

Define $\De=SU(3)/\text{conjugation} = \text{maximal torus}/\text{Weyl group}.$  This quotient
space is a topological 2-simplex, described in Section 6 in more detail.
The edges contain conjugacy classes of matrices with double eigenvalues,
and the vertices are the conjugacy classes of the central elements.

The map $\phi:\cC^\ell_{AB} \to \CC$ clearly factors through the map 
$\xi\colon \cC^\ell_{AB} \to \De$ sending  $\al \mapsto [\al(xy)], $
and the map $ \tr \colon \De \to \CC,$ which is smooth on the interior of the
simplex.
In Section 6 (following Hayashi \cite{Ha}) we identify the image
 $\xi(\cC^\ell_{AB})\subset \De$ (which we denote by $Q^\ell_{AB}$) as a convex polygon.  
Indeed, $ Q^\ell_{AB}$ is a hexagon if $\cC^\ell_{AB}$ is a Type I component (i.e. if
one of $A$ or $B$ has a repeated eigenvalue) and $ Q^\ell_{AB}$  is a nonagon if $\cC^\ell_{AB}$ is a Type II component (i.e. if  $A$  and $B$ each have three distinct eigenvalues).
 If $\cC^\ell_{AB}$ is a Type II component, then
$\xi^{-1}(p)$ is homeomorphic to a 2-sphere for all $p$ in the interior $ Q^\ell_{AB}$.

\begin{cor}  Set $\cC^{\ell,*}_{AB} = \cC^\ell _{AB} \cap R^*(Z, SU(3)).$ 
Then $\xi |_{\cC^{\ell, *}_{AB}} : \cC^{\ell, *}_{AB} \to \De$   
is a submersion except on the preimages of the intersection $ Q^\ell_{AB} \cap \partial \De$. 
\end{cor}

\begin{proof}  
Lemma \ref{sbmrs} effectively states that the differential of 
the composition 
$\tr \circ \xi |_{\cC^{\ell, *}_{AB}} \colon \cC^{\ell, *}_{AB} \to \CC$ 
has rank 2 except on $\xi^{-1}(\partial \De)$.  
By the chain rule, the same must apply to 
$\xi |_{\cC^{\ell,*}_{AB}}.$   
\end{proof}

When $\cC^\ell_{AB}$ is 4-dimensional, the structure of the fiber $\xi^{-1}(p)$ 
is described by Theorem \ref{fused}.  We summarize this information below.

 \begin{thm} \label{fiber}
Suppose $\cC^\ell_{AB}$ is a Type II component (i.e. 4-dimensional), and 
set $Q^{\ell,\red}_{AB} =\xi(\cC^{\ell, \red}_{AB})$.
Then $Q^{\ell,\red}_{AB}$ is
1-dimensional and the fiber of $\xi\colon \cC^\ell_{AB} \to \De$ over $p \in Q^\ell_{AB}$  is:
\begin{enumerate}
\item[(i)] A point if $p \in \partial Q^\ell_{AB},$
\item[(ii-a)] A smooth 2-sphere if $p \in {\sl Int} \, Q^\ell_{AB}$ and $p \not \in Q^{\ell, \red}_{AB},$
\item[(ii-b)] A pointed 2-sphere if $p \in {\sl Int} \, Q^\ell_{AB}$ and $ p   \in Q^{\ell, \red}_{AB}$.
\end{enumerate}
By a pointed 2-sphere, we mean a 2-sphere which is smooth  away from one point.

If $\al \colon \pi_1 Z \to SU(3)$ is a representation with $[\al] \in  \cC^\ell_{AB}$
such that
$\al(\la)$ does not have $1$ as an eigenvalue, then 
$p = \xi([\al]) \not\in Q^{\ell,\red}_{AB}$. If, in addition, $p\in {\sl Int} \, Q^{\ell}_{AB}$,
then it follows that $\xi^{-1}(p)$ is a smooth 2-sphere.
\end{thm}

\begin{proof}
The subset $\cC^{\ell, \red}_{AB}$ of reducible representations can be identified 
with the image under 
$\Psi: \Ga_A \backslash SU(3) / \Ga_B \to \cC^{\ell}_{AB}$ of the subset
$$\{g=(g_{ij}) \in SU(3) \mid g_{12}=g_{13}=0 \text{ or } g_{13}=g_{23}=0 \text{ or }  g_{12}=g_{23}=0\}
\subset SU(3).$$
This subset is
4-dimensional, and the principal orbits under the $\Ga_A \times \Ga_B$ action
are 3-dimensional (because their isotropy
subgroup of $\Ga_A \times \Ga_B$, which is 4-dimensional). 
Thus its image in  $\Ga_A \backslash SU(3) / \Ga_B$, and hence 
in $\cC^{\ell, \red}_{AB},$ is 1-dimensional.  

Suppose that $p=\xi([\al])\in Q^{\ell,\red}_{AB}$. Then we have
a reducible representation $\be\colon \pi_1 Z \to SU(3)$ with $[\be] \in \xi^{-1}(p)$.
Clearly $\be(xy)$ and $\al(xy)$ are conjugate in $SU(3).$ 
Since $\be$ is reducible and $\la$ lies in the commutator subgroup of $\pi_1(Z)$, it follows that
 $\be(\la)$ has (at least) one eigenvalue equal to $1$.  Because $\la=(xy)^{pq}h^{-(p+q)a}$ and $\al$ and $\be$ send $h$ to the same central element, it follows that $\al(\la)$ and $\be(\la)$ are conjugate, and hence $\al(\la)$ must also have $1$ as an eigenvalue.

The rest of the statement follows from Theorem \ref{fused}, and we 
explain the relationship between the different notations here and there. 
Suppose $A,B,C$ are diagonal $SU(3)$ matrices with  eigenvalues   
$\{e^{2 \pi ia_1}, e^{2 \pi i a_2}, e^{2 \pi i a_3}\}, \{e^{2 \pi i b_1},e^{2 \pi i b_2},e^{2 \pi i b_3}  \}$ and 
$\{e^{2 \pi i c_1},e^{2 \pi i c_2},e^{2 \pi i c_3}  \}$, respectively.  Then  
$\xi^{-1}([C])$, the preimage in $\cC^\ell_{AB}$ of the conjugacy class of $C$,
can be  identified with the moduli space $\cM_{\ba\bb\bc}$ 
described in Theorem \ref{fused}.

\end{proof}

\section{Perturbations} \label{s4}
The representation varieties for $\Si$ and $Z$  discussed in the
previous sections can be identified with the moduli spaces  
of flat $SU(3)$  connections on $\Si \times SU(3)$ and $Z\times
SU(3)$.  The principal advantage of this perspective is that
flat moduli space is the critical set of a function on the space of
all connections, modulo gauge, and this gives a framework to perturb
  for transversality purposes.  In particular, we deform the function
of which the flat moduli space is the critical set, and consider the
critical set of the deformed function to be the ``perturbed moduli
space.''

After introducing some notation,
we will define the twisting perturbations
and analyze their effect on the moduli space.
Of central importance is the behavior of pointed 2-spheres
under  twisting perturbations.
In Subsection \ref{TEOT},
we show that  under a twisting perturbation, every pointed
2-sphere resolves into two pieces:
an isolated reducible orbit and a smooth, nondegenerate
2-sphere.
\subsection{Gauge theory preliminaries and the odd signature operator}
Fix a 3-manifold $X$ with Riemannian metric.
Let:
$\cA(X)$ be the space of
  $SU(3)$ connections over $X$,
completed in the $L^2_1$ topology,
  $\cG(X)$  be the group of $SU(3)$ gauge transformations,
completed in the $L^2_2$ topology,
$\cB(X)$  be the quotient $\cA(X)/\cG(X)$ , and
$\cM(X)$ be the moduli space of gauge orbits of flat connections.

\medskip

When the manifold is clear from context, we will drop it from
the notation and simply write $\cA, \cG, \cB$ and $\cM$.

The spaces $\cA, \cB,$ and $\cM$ are
stratified by levels of reducibility, and we
adopt a notation  consistent with that used for
the representation varieties. In particular:
\begin{enumerate}
\item[(i)] $\cM^*$ is the moduli space of irreducible flat
$SU(3)$ connections.
\item[(ii)] $\cM^\red$ is the moduli space of reducible, nonabelian flat
$SU(3)$ connections.
\item[(iii)] $\cM^\ab$ is the moduli space of abelian, flat $SU(3)$
connections.
\end{enumerate}

Given
an $SU(3)$ connection $A$, covariant differentiation defines
a map $$d_A\colon \Om^0(X;su(3)) \to \Om^1(X;su(3)).$$
If $A$ is flat, we obtain the twisted de Rham complex
\begin{equation}\label{deco}
\Om^0 (X;su(3)) \stackrel{d_A}{\lto} \Om^1 (X;su(3)) \stackrel{
d_{A}}\lto
\Om^2 (X;su(3)) \stackrel{d_A}\lto \Om^3 (X;su(3))
\end{equation}
with cohomology groups
$H^0_A(X;su(3)) = \ker(d_A\colon \Om^0 \to \Om^1)$, the Lie algebra of the
stabilizer of $A$, and
$H^1_A(X;su(3)) = \ker(d_A\colon \Om^1 \to \Om^2)/ \im(d_A\colon \Om^0 \to \Om^1),$
the Zariski tangent space of $\cM$ at $[A]$. When $X$ is closed, the
Hodge star  isomorphism $\star\colon \Om^i (X;su(3))\to  \Om^{3-i}
(X;su(3))$ induces isomorphisms $H^i_A(X;su(3))\cong
H_A^{3-i}(X;su(3))$.
The  de Rham theorem  for twisted cohomology
gives isomorphisms $H^i_A(X;su(3)) \cong H^i(X;su(3)_\al)$
where $\al\colon \pi_1(X) \to SU(3)$ is the holonomy representation of
the flat connection $A.$

In this section, we will consider perturbations
of the moduli space.
The perturbations we use are of Floer
type, meaning that we perturb the flatness equations
in a neighborhood of a finite collection of loops in $X$ (see
Definitions \ref{tpd} and \ref{pfeq} below.)
Given an admissible perturbation $h$,
a connection $A \in \cA$ is called $h$-perturbed flat  if
$F_A = \star 4 \pi^2 \nabla h(A).$
We denote the moduli space of $h$-perturbed flat
$SU(3)$ connections by $\cM_h$.
For more details on perturbations in the $SU(3)$
context, see Section 2.1 in \cite{BH1}.

If $A$ is $h$-perturbed flat, we define
$d_{A,h}= d_A - \star 4 \pi^2 \Hess h(A)$ and
the perturbed deformation complex
\begin{equation}\label{pdeco}
\Om^0 (X;su(3)) \stackrel{d_A}{\lto} \Om^1 (X;su(3)) \stackrel{
d_{A,h}}\lto
\Om^2 (X;su(3)) \stackrel{d_A }\lto \Om^3 (X;su(3)).
\end{equation}
The argument that (\ref{pdeco})
is Fredholm when $X$ is closed can be found in \cite{T} or \cite{H}.
The
first cohomology  of this Fredholm  complex is denoted
$H^1_{A,h}(X;su(3)) = \ker(  d_{A,h})/\im(d_A)$
and is the Zariski tangent space of $\cM_h$ at $[A]$. The
cohomology of this complex is independent of the choice of
Riemannian metric on $X$ since $d_A$ and $d_{A,h}$ are.

\begin{defn} The {\em odd signature operator}  twisted by a connection
$A$ is the linear elliptic differential operator
\begin{eqnarray*}
&D_A\colon  \Om^{0+1}(X;su(3))  \lto \Om^{0+1}(X;su(3))& \\
&D_A(\si,\tau) = (d_A^* \tau, d_A \si + \star d_A \tau).&
\end{eqnarray*}
 It can obtained by folding up the complex (\ref{deco}), although $D_A$
 is defined whether or not $A$ is flat. It is a generalized Dirac
 operator (in the sense of \cite{BW}).

The
{\em perturbed odd signature operator} is defined similarly using the
complex (\ref{pdeco}),  for a connection $A$ and a perturbation $h$ to
be
\[ D_{A,h}\colon  \Om^{0+1}(X;su(3))  \lto \Om^{0+1}(X;su(3))  \]
\begin{eqnarray*}
  D_{A,h}(\si,\tau) &=& (d_A^* \tau, \ d_A \si + \star d_{A,h} \tau) \\
& = & \left(d_A^* \tau, \ d_A \si + \star d_A \tau - 4\pi^2 \Hess h(A)(\tau)\right)\\ &=&
D_{A}(\si,\tau)+ (0,-4\pi^2 \Hess h(A)(\tau)).
\end{eqnarray*}
Here we use the metric to view
  $\Hess h(A)(\tau)$ as a 1-form with $su(3)$ coefficients.
\end{defn}

The Hessian is bounded as a map from $L^2$ to $L^2$
  (\cite{T}, \cite{BH1}, \cite{HKL}).  Thus
  the composite of the compact inclusion of $L^2 _1 \to L^2$ with the
  bounded Hessian $L^2\to L^2$ is a compact map $L^2 _1 \to L^2$, and
  the addition of the Hessian to the signature operator is a compact
  perturbation.
Since $D_{A,h}$ differs from $D_A$ by a compact perturbation, it is again Fredholm  when $X$ is
closed.

 The usual Hodge theory argument
shows that if $X$ is closed, the kernel of $D_{A,h}$ is isomorphic to $H^0_{A,h}(X;su(3))\oplus
H^1_{A,h}(X;su(3))$.  When $X$ is not closed, then $D_{A,h}$ is not Fredholm. The operators $D_{A}$
and $D_{A,h}$ are symmetric: $\langle D_{A,h}(\phi_1),\phi_2\rangle= \langle\phi_1,
D_{A,h}(\phi_2)\rangle$ if $\phi_1$ and $\phi_2$ are supported on the interior of $X$. Thus if $X$
is closed $D_{A}$ and $D_{A,h}$ are self-adjoint.

The operator $D_{A,h}$ is not local. It is not a differential nor pseudodifferential operator.
However, $(D_{A,h}-D_A)(\phi)$ depends only on the restriction of $A$ and $\phi$ to the compact
domain in $X$ along which the perturbation is supported   and moreover $(D_{A,h}-D_A)(\phi)$
vanishes outside of this domain (in the case considered in this article, the compact domain is a
neighborhood of the $r$-singular fiber). The proof of this fact can be found e.g.~ in \cite{HKL},
and follows in the present context straightforwardly from Proposition \ref{grad}.

Basic for us will be  the splitting \begin{equation}\label{splitSi}\Si(p,q,r) = Y \cup_T Z.
\end{equation}
Here, $$T = S^1\times S^1=\{(e^{ix},e^{iy})\}$$ is the 2-torus with the product metric and
orientation so that $dx dy$ is a positive multiple of the volume form. Its fundamental group
$\pi_1(T)$ is generated by the loops $\mu=\{(e^{ix},1)\}$ and $ \la=\{(1,e^{iy})\}.$

The 3--manifold $Y$ is the solid torus
$$Y =D^2\times S^1 = \{ (re^{i x},e^{iy}) \mid 0 \leq r \leq 1\}$$
  oriented so that
$dr dx dy $ is a positive multiple of the volume form; it is a neighborhood of the $r$-singular
fiber in $\Si(p,q,r)$. Choose a metric on $Y$ so that a collar neighborhood of the
  boundary is isometrically identified with $[-1,0]\times T$.
  As oriented manifolds,
  $\partial Y=\{0\}\times T$.
  The fundamental group $\pi_1(Y)$ is infinite cyclic generated
     by the longitude $\la$. (The meridian $\mu$
     bounds the disc $D^2\times\{1\}$ and so is trivial in $\pi_1(Y)$.)

The 3--manifold $Z$ is the complement of an open tubular
   neighborhood of the $r$-singular fiber in $\Si(p,q,r)$.
Choose
a metric on $Z$ so that a collar neighborhood of the boundary
$\partial Z$ is isometrically identified with $[0,1]\times T$,
and $\la$ is
null-homologous in $Z$.
As oriented manifolds,  $\partial Z=-\{0\}\times T$.

The metrics on $Y$ and $Z$ induce one on $\Si$ with the property
that a bicollared neighborhood of  $T \subset \Si$ is
isometric to $[-1,1] \times T$. We call  $[-1,1] \times T$ the neck.
Every connection $A$ on $\Si$ which is flat on the neck is gauge equivalent
to one in {\it cylindrical form}, meaning that
its restriction $A|_{[-1,1] \times T}$ to the neck
is the pullback of a
connection on the torus under the projection
$[-1,1] \times T \to T.$
There are similar results for
$Y$  using the collar $[-1,0] \times T \subset Y$
and for $Z$ using $[0,1] \times T \subset Z$.
A connection  in cylindrical form and which is
flat on the neck is gauge equivalent to one
whose
meridinal and longitudinal
holonomies  are diagonal.

\subsection{The twisting perturbation on the solid torus}\label{ttp}
In this subsection, we define the   twisting
perturbation and  study the perturbed flatness
equations on the solid torus.
The crucial issue is to determine
which flat connections on the boundary
extend as perturbed flat
connections over the solid torus.

We begin with some notation.
For a complex number $\zeta,$ let
$\Re(\zeta) $ be its real part and $\Im(\zeta)$
its imaginary part.
Recall the parameterization $\Phi\colon  \RR^2 \to SU(3)$
   of the maximal torus $T \subset SU(3)$ given by
   equation (\ref{defnPhi}).

We use $x=(x_1,x_2)$ for the  coordinates on the
2-disk $D^2$
and  $\th$  for the circle $S^1$.
Suppose $\eta\colon  D^2 \to \RR$ is a radially symmetric nonnegative
function
supported in a small neighborhood of $x=0$ with
$\int_{D^2} \eta(x) \ dx = 1$.

Fix a basepoint $\th_0 \in S^1$.
For a connection $A$ on the solid
torus $  D^2 \times S^1$, let
  $\hol_x (A)$ denote its holonomy
around $\{ x\} \times S^1$
starting and ending at $(x,\th_0).$ Although  $\hol_x (A)$
depends on the choice of basepoint,
  its trace
$\tr \hol_x(A)$  is independent
of this choice.

\begin{defn}\label{tpd}
Define the {\it twisting perturbation} function
$f\colon \cA(D^2 \times S^1) \to \RR$ by setting
\begin{equation}
f(A) = - \tfrac{1}{4\pi^2}
  \int _{D^2} \Im (\tr \hol _{x }(A)) \eta(x) \ dx
\label{tpde}
\end{equation}
\end{defn}

Let $M_3(\CC)$ be the vector space of $3\times 3$ complex matrices and regard $su(3)$ as a subspace
of $M_3(\CC)$. Define $\Pi_{su(3)}\colon M_3(\CC)\to su(3)$ to be orthogonal projection with
respect to the standard inner product on $M_3(\CC)$.

\begin{prop} \label{grad}
The gradient of the perturbation of  (\ref{tpde}) is given by $$\nabla f(A)=-\tfrac{1}{4\pi^2}
\Pi_{su(3)}\left( i\hol_x(A) \right) \eta(x) \, d \th.$$
\end{prop}
\begin{proof}
If $A$ is a connection on $S^1$ and $\al$ is an $su(3)$-valued 1-form
on $S^1$,
then Proposition 2.6, \cite{BH1} gives the
differentiation formula
$$
\left. \frac{d}{ds}\Im \tr \hol_{x}(A+s\al) \right|_{s=0} =
\Im \tr (\hol_{x}(A) {\textstyle \int_{S^1}\al}),$$
where $\int_{S^1} \al$ is interpreted as in Section 6 of
   \cite{BH1}.

  From  equation (\ref{tpde}), $f(A+s\al)$ is clearly
independent of all components
of $\al$ except the $d\th$ component.  We can find its derivative by
integrating
the formula in the circle case:
\begin{equation}\label{der1}
\left. \frac{d}{ds} f(A+s\al)\right|_{s=0} = -\tfrac{1}{4\pi^2}
\int_{D^2} \Im \tr \left(\hol_{x}(A) ~ {\textstyle \int_{S^1}\al}
\right)
\eta(x) \ dx
\end{equation}
Since $\int_{S^1} \al$ is $su(3)$-valued, we have
\begin{eqnarray}
\Im \tr ( \hol_{x}(A) ~ {\textstyle \int_{S^1}\al} )&=&
-\Re \tr\left( i~ \hol_{x}(A) {\textstyle \int_{S^1}\al} \right)\cr
&=& \left\langle   \Pi_{su(3)} (i\hol_{x}(A) ),  {\textstyle
\int_{S^1}\al}
    \right\rangle_{su(3)},
\end{eqnarray}
where we identify $-\tr (AB)$ with the standard inner product
$\langle \cdot, \cdot \rangle _{su(3)}$
on $su(3)$.
Therefore equation (\ref{der1}) can be rewritten as
$$
\left. \frac{d}{ds} f(A+s\al)\right| _{s=0}=\left\langle
   -\tfrac{1}{4\pi^2}  \Pi_{su(3)} (i\hol_{x}(A) ) \eta(x) \, d\th ,
    \al  \right\rangle_{L^2(D^2 \times S^1)}.  $$
Here, $\hol_x(A)$ is interpreted
as a section of the bundle $\End (E)$ of endomorphisms
of the rank three bundle $E \to D^2 \times S^1.$
The section $\hol_x(A)$
is covariantly constant around the
circle fibers with respect to the induced connection
on $\End(E)$.
\end{proof}

\begin{defn}\label{pfeq} Given $t\ge 0$, a connection $A$ on the solid torus
is called {\em $(tf)$-perturbed flat}
if it satisfies the equation
$$F_A=\star 4\pi^2 t \nabla f(A)$$
where $F_A$ denotes the curvature of $A$.
Since $\eta$ is supported on a small neighborhood of $0\in D^2$,
a $(tf)$-perturbed flat connection is flat
near the boundary torus (see Proposition \ref{foc} below).
\end{defn}

The next two propositions are well-known.
The first  was initially observed by Floer in
\cite{F}. Its proof is based on the previous observation
that a perturbed flat connection
has curvature only in the $dx_1 dx_2$ direction.
\begin{prop} \label{lhi}
Suppose $A$ is a  connection
on the solid torus.
If $A$ is
$(tf)$-perturbed flat, then $\hol_x (A)$ is
independent of $x \in D^2$.  \end{prop}
\begin{proof}
On the disk $D^2 \times \{\th_0\},$ trivialize the
$SU(3)$ bundle using radial parallel translation
starting at the center $(0,\th_0)$.
For each $x \in D^2,$
take the line segment
$\overline{0x}$ and consider the annulus $\overline{0x} \times S^1$.
Since $\star F_A = i \, d\th,$ the restriction of $A$ to this annulus
is flat.
But parallel translation along the line segment $\overline{0x}$
is trivial, and so $\hol_x(A) = \hol_0(A)$ and is independent of $x \in
D^2$.
\end{proof}

Proposition \ref{lhi} shows that for a perturbed flat connection $A$ on the solid torus, we can denote
$\hol_x(A)\in SU(3)$ unambiguously by $\hol_{\la}(A)$.  We call this
the {\em longitudinal holonomy} of $A$. The holonomy of $A$ along
the meridian $\del D^2\times\{\theta_0\}$ is called the {\em
meridinal holonomy}.

The next result states that perturbed flat connections are
flat outside a neighborhood of the perturbation curves.
\begin{prop}  \label{foc}
If $A$ is perturbed flat with respect to a
perturbation $h$ supported on
   a single thickened
curve $\ga\colon  D^2\times S^1 \to \Si$, then $A$ is flat
on the complement $\Si-\ga(D^2 \times S^1)$.
\end{prop}
\begin{proof}  Under the hypothesis, the equation for perturbed
flatness is
$\star F_A =4\pi^2 \nabla h(A)$, but $\nabla h(A) =0$ outside  the
image $\ga ( D^2 \times S^1)$.
\end{proof}

The twisting perturbation is well-defined as a function $$f\colon\cA(\Si(p,q,r))\to \RR,$$ once one
fixes a framing on the solid torus $Y$ in the decomposition (\ref{splitSi}).
 We use the framing $Y\cong D^2\times S^1$ in
which the longitude $\la$ is homotopic to $\{x\} \times S^1$ in the complement of $K$ for all
nonzero $x\in D^2.$ We assume further that the bump function $\eta(x)$ is supported in a small
enough neighborhood that it vanishes on the neck $[-1,1] \times T$. Proposition \ref{foc} then
implies that every $(tf)$-perturbed flat connection $A$ on $\Si$ restricts to  a flat connection on
$([-1,0] \times T) \cup Z$. The definition of $f$ and Proposition \ref{grad} show that $f(A)$,
$\nabla f(A)$, and $\Hess f(A)$ depend only of the restriction of $A$ to the interior of $Y$.

\medskip

The last result in this subsection determines an equation
on meridinal and longitudinal holonomies that a
connection $A$ must satisfy in order for it
to be  $(tf)$-perturbed flat.

\begin{prop} \label{pfhc}
    Suppose that $A$ is a  connection
on $D^2 \times S^1$ which is perturbed flat with
respect to the twisting perturbation $tf$.  Then there is a smooth
gauge representative for $[A]$.  Furthermore, if
  $\hol_{\la}(A)=\Phi(u,v)$, then
   the meridinal holonomy is given by
\begin{equation}
\label{tnoia}
\hol_{\mu}(A) =\Phi\left(- t\sin u \sin v,
\tfrac{ t}{3}(  \cos u \cos v - 2\cos^2 v +1)\right).\end{equation}
\end{prop}

\begin{proof} The smoothness property holds for all holonomy type
perturbations, not just the twisting perturbation we have defined
here.  This is claim (1) of Lemma 8.3 in \cite{T}.

The second claim is a generalization to $SU(3)$ (and imaginary part of
trace)
   of a well-known fact for $SU(2)$
perturbed flat connections, going back to Floer.
Note first that $\nabla (tf)=t\nabla f$.  Let $A$ be a smooth
$tf$-perturbed
flat connection, gauge transformed so that $\hol_\la (A)$ is
diagonal.

Since the curvature $F_A=
\star 4\pi^2 t \nabla f(A)$ takes only diagonal matrix
values, we can find the meridinal holonomy by integrating $F_A$ over a
disk
that the meridian bounds, namely
\begin{eqnarray*}
\hol_\mu (A) 
&=& \exp \left( - {\textstyle \int_{\partial D^2 }} A \right) \\
&=&  \exp \left( -{\textstyle \int _{D^2}} dA \right) \\
&=& \exp \left( -{\textstyle \int_{D^2}} F(A) \right) \\
&=&  \exp \left(- {\textstyle \int_{D^2}} 4\pi^2 \star \nabla (tf)(A) \right)  \\
&=& \exp \left(  t{\textstyle \int_{D^2}} \Pi_{su(3)} (i\hol_\la (A) )
\, \eta(x) \, dx_1 \wedge dx_2    \right) \\
&=&  \exp \left(  t \Pi_{su(3)}  \left( i \hol_\la (A)
\right)  \right).
\end{eqnarray*}
 The projection of a diagonal matrix $B$ onto $su(3)$ is given by
taking the imaginary part of $B - \frac 1 3 \tr (B) I$.
Applying this to $ti\hol_{\la}(A)$
shows that
\begin{eqnarray*}
 \Pi_{su(3)} \left( i \hol_\la (A) \right)
&=& \Pi_{su(3)}  i \Phi(u,v) \\
&=&\Im \left[ i\Phi(u,v) - \tfrac{i}{3} \tr \Phi(u,v)I \right] \\
 &=& \left[ \begin{array}{ccc}i a_1 & 0&0\\ 0&ia_2&0 \\ 0&0&
ia_3
\end{array} \right],\end{eqnarray*}
where \begin{eqnarray*}
a_1&=& \tfrac{1}{3} \left( 2 \cos(u+v) -\cos(-u+v) -\cos(2v)\right) ,\\
a_2&=& \tfrac{1}{3} \left(  -\cos(u+v) +
2\cos(-u +v) - \cos(2v) \right), \\
  a_3 &=& \tfrac{1}{3} \left(  -\cos(u+v) -\cos(-u+v) +
2\cos(2v)\right) .\end{eqnarray*}

Setting $\tilde u = \tfrac{a_1 -a_2}{2}$ and $\tilde v = \tfrac{a_1 + a_2}{2}$ and applying the angle addition formulas, we see
that $\tilde u =  -\sin u \sin v$ and $\tilde v= \tfrac{1
}{3} \left( \cos u \cos v - 2 \cos^2 v + 1 \right).$
These substitutions simplify the formula for $\hol_\mu(A)$ to give
\begin{eqnarray*}
\hol_\mu (A)
 &=& \Phi\left(  -t\sin u \sin v,
\tfrac{t}{3} \left( \cos u \cos v - 2 \cos^2 v + 1 \right)
\right).\end{eqnarray*}
\end{proof}
 \begin{remark}\label{simplerem}
Notice that if $\hol_\la(A)= \Phi(u,0)$ in the above proposition
(namely if $v=0$), then the conclusion is
that $\hol_\mu(A) = \Phi(0,\tfrac{ t}{3} (\cos u-1)).$
  \end{remark}

\subsection{The effect of the twisting perturbation on a pointed
2-sphere}
\label{TEOT}
We now consider
  twisting perturbations on $\Si = Y \cup_T Z$
supported on the solid torus $Y$. In the
last subsection we showed  that any
perturbed flat connection $A$ on $\Si$ is indeed flat
on $Z$ (Proposition \ref{foc}) and we obtained an
equation that the meridinal and longitudinal
holonomies must satisfy to
extend as a perturbed flat connection on $Y$
(Proposition \ref{pfhc}).
In this subsection, we use this equation   to analyze the topology of
the perturbed flat moduli space.
We are  particularly interested in the effect
of  the twisting perturbation on the
  pointed 2-spheres
in $\cM$. We show that the perturbed flat moduli space
near a pointed 2-sphere resolves into two pieces:
an isolated gauge orbit of reducible connections
and a smooth, nondegenerate 2-sphere of gauge orbits
of irreducible connections.

We identify the perturbed flat moduli space
$\cM_{tf}(\Si)$ as the subset of the flat moduli space
$\cM(Z)$ of gauge orbits which extend as perturbed
flat connections over the solid torus.
We explain the geometric picture before going into
details.

The moduli space
  $\cM(T)$  is the quotient
of the product of two copies of the
maximal torus of $SU(3)$ modulo the diagonal action
of Weyl group $S_3$, the group
of symmetries on three letters.
Thus $\cM(T)$ is 4-dimensional.

With respect to the splitting
 $\Si = Y \cup_T Z,$
we have restriction maps
$$\zeta_Y\colon \cM(Y) \to \cM(T), \quad \zeta_Z ^\red \colon \cM^\red
(Z)\to \cM(T) \quad \hbox{and} \quad
\zeta_Z ^*\colon \cM^* (Z) \to \cM(T)$$ defined by sending $[A]$
to $[A|_T]$. Denote the images of these maps by
$\cZ_Y = \im(\zeta_Y)$, $\cZ_Z ^\red = \im(\zeta_Z ^\red)$ and $\cZ_Z ^* =
\im(\zeta_Z ^*)$.
We also have restriction maps $r_Z\colon \cM(\Si)\to \cM(Z)$ and $r_Y \colon \cM(\Si)\to \cM(Y)$,  a commutative diagram
\begin{equation}\label{restrictiondiagram}
\begin{diagram}
\node[2]{\cM(Z)}\arrow{se,t}{\zeta_Z}\\
\node{\cM(\Si)}\arrow{ne,t}{r_Z}\arrow{se,b}{r_Y}\node[2]{\cM(T)}\\
\node[2]{\cM(Y)}\arrow{ne,b}{\zeta_Y}
\end{diagram}
\end{equation}
and similar diagrams for the reducible and irreducible moduli spaces.

All three of  $\cZ_Y$, $\cZ_Z ^\red$ and $\cZ_Z ^*$ are codimension two
submanifolds of $\cM(T)$.
The map  $r_Z$ is injective. This is just the statement
that the flat connections on $\Si$ can be identified with
those flat connections on $Z$ which
extend flatly over
the solid torus $Y$. The crux of the matter is that
the flat extension to  the solid
torus is uniquely determined by $A|_T$ up to gauge transformation.

Thus   the moduli space $\cM ^\red(\Si)$ can be identified with
 $$r_Z^\red(\cM^\red(\Si))=(\zeta_Z^\red)^{-1}(\cZ_Y)=\{ [A] \in \cM^\red (Z) \mid  [A|_T] \in \cZ_Y \cap \cZ_Z ^\red\},$$
 and  likewise we can identify $\cM ^* (\Si)$
as the subset of $\cM^* (Z)$ given by
$$r_Z^*(\cM^*(\Si))=(\zeta_Z^*)^{-1}(\cZ_Y) = \{ [A] \in \cM^* (Z) \mid [A|_T] \in \cZ_Y \cap \cZ_Z ^* \}.$$

If  $[A_0]$ lies on a pointed 2-sphere, then $\zeta_Z^\red$ and $\zeta_Z^*$ are individually
transverse to $\cZ_Y$ at $[A_0|_T]$. But
  $\cZ_Y$ intersects both $\cZ_Z ^\red$ and $\cZ_Z ^*$ at
$[A_0|_T]$, causing difficulties.
The reducible part $(\zeta ^\red _Z) ^{-1} ([A_0 |_T])$  is simply $[A_0]$,
 while the irreducible part $(\zeta ^*_Z ) ^{-1} ([A_0|_T])$
  is the complement of $[A_0]$ in the pointed
2-sphere (and in particular is not compact).

To make $\cM(\Si)$ non-degenerate, we apply a   twisting
perturbation which moves $\cZ_Y$ slightly.  As with the flat moduli space, we have a restriction
map $\zeta_{Y,tf}\colon \cM_{tf}(Y) \to \cM(T)$
defined by sending  $[A] \in \cM_{tf}(Y)$
to $[A|_T]$. (Recall that $A|_T$ is necessarily flat.)
Denote the image
of this map by $\cZ_{Y,tf} = \im(\zeta_{Y,tf})$.
As before, we can  identify
the strata of reducible and irreducible gauge orbits
in the perturbed flat moduli space $\cM_{tf}$
as the subsets of $\cM(Z)$ given by
$$ \cM_{tf} ^\red (\Si) = (\zeta_Z^\red)^{-1}(\cZ_{Y, tf})
=\{ [A] \in \cM^\red(Z) \mid [A|_T] \in \cZ_{Y,tf} \cap \cZ_Z ^\red\},$$
and
$$ \cM_{tf} ^*(\Si)= (\zeta_Z^*)^{-1}(\cZ_{Y, tf})
= \{ [A] \in \cM^*(Z) \mid [A|_T] \in \cZ_{Y,tf} \cap \cZ_Z ^*\}.$$

We will show that  for small $t>0$
$\cZ_{Y,tf}$ intersects $\cZ_Z ^\red$ and
$\cZ_Z^*$ at points with nondegenerate preimages in the following sense:
$(\zeta_Z^\red)^{-1}(\cZ_{Y,tf})$ is  an isolated reducible connection $[A]$ with
  $H^1 _A(Z;\CC^2 )=0$, and $(\zeta_Z^*)^{-1}(\cZ_{Y,tf})$  a smooth 2-sphere.

We will show this to be the case by determining,  to first order in
$t$,  where this intersection
point lies.
  The idea is to pin down their meridinal and longitudinal
holonomies.

  Throughout the remainder of this section, $A_0$ will be a fixed
  reducible flat connection whose gauge orbit $[A_0]$ lies
on a 2-sphere component.
Identify $\cM^\red(\Si)$ with  $\cM^*_{S(U(2)\times U(1))}(\Si)$
and note that $[A_0]$ is a regular point of this latter
moduli space. This follows because $[A_0]$ can be represented by
  an $SU(2)\times \{1\}$ connection $A_0$
and Proposition \ref{coho} implies that
$$H^1_{A_0}(\Si; s(u(2)\times u(1)))
=H^1_{A_0}(\Si; su(2))\oplus H^1_{A_0}(\Si;u(1))=0.$$

Regularity of $\cM^\red(\Si)$  near $[A_0]$ implies that,
for $0\leq t \leq \ep$, there is a family
of reducible $(tf)$-perturbed flat  connections $A_t$ which are
deformations of
$A_0$.
Our first goal is to
show that $[A_t]$ is an isolated point  in the perturbed
flat moduli space $\cM_{tf}$.

\begin{prop} \label{redps}  Assume $[A_0]$ is a gauge orbit
of reducible flat connections on
$\Si$ that lies on a 2-sphere component.
Choose a representative $A_0$ in cylindrical form
whose holonomy
on the torus $T$ is diagonal.
Equation (\ref{rroott}) gives that
  $$\hol_{xy}(A_0) = \Phi\left(\tfrac{2 \pi m}{r},0\right)$$ for some
integer
$m$ with $0 < m < r.$
For $0 \leq t \leq \ep$, let $[A_t]$ be the family of
gauge orbits of
  reducible $(tf)$-perturbed flat connections near $[A_0]$.
  As before, choose representatives  in cylindrical form.
Since each  $A_t$ restricts to a flat connection on $Z$,
we can also arrange
that $A_t$ has diagonal holonomy on the torus $T$.
Then   the
holonomies satisfy:
\begin{eqnarray*}
\hol_{\mu}(A_t) &=& \Phi\left(2 \pi m, \tfrac t 3  \left(\cos \left(
\tfrac{2 \pi m}{r}\right) -1\right)\right),\\
\hol_{\la}(A_t) &=& \Phi\left(\tfrac{2 \pi m p q}{r}, 0\right).
\end{eqnarray*}
\end{prop}
\begin{proof} The Implicit Function Theorem
 implies the path
$[A_t]$ is smooth.  As with single connections, the path of gauge
representatives for $[A_t]$ can be chosen to be smooth, in cylindrical
form, and with the property that $\hol_{xy}(A_t)$ and $\hol_h (A_t)$
are diagonal.  Note that by Proposition 4.6 these connections are flat
on $Z$.

Equation (\ref{rroott}) and the discussion immediately preceding it
imply that
$$\hol_{xy}(A_0) = \Phi(\tfrac{2\pi m}{r}, 0)
\quad \text{and} \quad  \hol_h (A_0) = I.$$
Therefore, $$\hol_{xy}(A_t) = \Phi(u_t, v_t) \quad \text{and} \quad \hol_h (A_t) =
\Phi(0,w_t)$$ for some functions $u_t$, $v_t, w_t$ satisfying
$u_0=\tfrac{2\pi m}{r}$, $v_0 =0=w_0$.
Here we know that $\hol_h (A_t)$ has the form stated because it
commutes with the nonabelian representation $\hol (A_t)\colon \pi_1 Z \to
S(U(2)\times U(1))$, so it is in the center of $S(U(2)\times U(1))$.

It follows from equation (\ref{mu}) that
$$
\hol_\mu (A_t) = \Phi (ru_t, rv_t + cw_t)\quad \text{and} \quad \hol_\la (A_t) =
\Phi ( pqu_t , pq v_t -(p+q)aw_t).$$
Proposition \ref{rolm} shows that the second argument in $\hol_\la
(A_t)$, namely $pqv_t - (p+q)aw_t $, must equal zero.
Proposition \ref{tnoia} (see Remark \ref{simplerem}) now implies that
\[ \Phi(ru_t, rv_t + cw_t)= \Phi(0, \tfrac t 3
(\cos(pqu_t ) - 1) ).\]

 From this it follows that $u_t = \frac{2\pi m}{r}$, independent of
$t$, and that $rv_t + cw_t = \frac t 3 \left(
  \cos ( \frac{2\pi pqm}r) - 1 \right)$.
\end{proof}

\begin{cor}  \label{Z cohom goes away} For small enough $t>0$, the representation
$\al_t\colon \pi_1 (Z) \to SU(3)$ induced by the reducible flat
connection $A_t$ is twisted (i.e. takes values in $S(U(2)\times U(1))$
but not in $SU(2)\times \{1\}$) and satisfies $H^1(Z; \CC ^2
_{\al_t}) = 0$.  \end{cor}
\begin{proof}  Proposition \ref{redps} shows $\hol_\mu (A_t)$ is
twisted, and therefore $\al_t$ is twisted.  The cohomology claim
then follows from Proposition \ref{abelians}.
\end{proof}

Corollary \ref{Z cohom goes away} will be used in Section \ref{s5} to
show that, for small $t$, the orbit $[A_t]$  of reducible perturbed
flat connections near $[A_0]$ is isolated in $\cM _{tf} (\Si)$.

We now turn our attention to  understanding the effect
of the twisting perturbation on the stratum of irreducible
connections.
We continue to assume that $A_0$ is a reducible flat connection, in
cylindrical form, with $\hol_{xy}(A_0)$ diagonal, and that $[A_0]$
lies on a pointed 2-sphere.
As pointed out in the proof of Proposition \ref{redps},
  there is  an integer $m$ with $0 < m < r$ such that
$\hol_{xy}(A_0)=\Phi(\frac{2 \pi m}{r},0)$ and
$\hol_\la(A_0) = \Phi(\frac{2 \pi pq m}{r},0).$
 
Now consider an irreducible $(tf)$-perturbed flat
  connection $A$  near   $A_0.$
We assume $A$ is in cylindrical form on the neck and
that the meridinal and longitudinal holonomies of $A$ are diagonal.
Since $\hol_\ga(A)$ is close to $\hol_\ga(A_0)$ for
all $\ga \in \pi_1(Z),$ we can write
\begin{equation} \label{uvwz}
\hol_\la(A)=\Phi(u, v) \quad \hbox{and} \quad \hol_\mu(A)=\Phi(w,z)
\end{equation}
for $(u,v)$ near $(\frac{2 \pi pq m}{r},0)$
and $(w,z)$ near $(0,0)$.
Because
the restriction of $A$ to $Z$ is irreducible and flat,
$\hol_h(A) =I$. (To see this, note that
$h \in \pi_1(Z)$ is central and $\hol_h(A)$ is
a priori near $\hol_h(A_0)=I$.)
Equation (\ref{mu}) now implies that
$$\left(\hol_\mu(A) \right)^{pq} = \left(\hol_\la(A)\right)^r,$$
and plugging this into equation (\ref{uvwz}) gives that
\begin{equation} \label{isl}
w = \frac{ru}{pq}  - {2\pi m}
\quad \hbox{and} \quad
z = \frac{rv}{pq}.
\end{equation}
On the other hand, if  $A$
extends as a $(tf)$-perturbed flat connection
over $Y,$ equation (\ref{tnoia}) implies that
\begin{equation} \label{pert uv}
w = -t \sin u \sin v \quad \hbox{and} \quad   z= t(\cos u  \cos v
-2\cos^2 v  + 1).\end{equation}
Combining equations (\ref{isl}) and (\ref{pert uv}), we obtain a
pair of equations (depending on the parameter $t$) which determine
$u$ and $v$.

We now solve for $u$ and $v$ to first order in $t$.
  To facilitate the argument, define
the function $P\colon \RR^3 \to \RR^2$ given by

$$P(t,u,v)=\left( \left(\tfrac{ru}{pq}\right)  t\sin u \sin v, \ \tfrac{rv}{pq} -
\tfrac{t}{3}\left( \cos u \cos v - 2\cos ^2 v + 1 \right) \right).$$

The map $(u,v) \mapsto P(0,u,v)$ is clearly a submersion,
and  the
Implicit Function Theorem
provides smooth functions $u(t)$ and $v(t)$
  near $t=0$ such that $(t, u(t), v(t))$ parameterizes
the solutions of the equation
$P(t,u,v)=0$ near $\left(0, \frac{2 \pi pqm}{r}, 0\right)$.
Differentiating the
equation $P(t,u(t), v(t))=0$ with respect to $t$ at $t=0$ yields
$$
u'(0)=0 \quad \hbox{and} \quad
v'(0)= \tfrac{pq}{3r} \left( \cos\left(\tfrac{2\pi pqm}{r}\right)
-1\right).  $$
Thus any irreducible
$(tf)$-perturbed flat connection $A$ near $A_0$
satisfies:
\begin{equation} \label{holy}
\begin{split}
\hol_{\la}(A) &= \Phi
\left(\tfrac{2\pi pqm}{r}, \tfrac{tpq}{3r} \left(
\cos\left(\tfrac{2\pi pqm}{r}\right)  - 1 \right) \right) +O(t^2),  \\
\hol_{\mu}(A) & = \Phi\left( 0, \tfrac{t}{3}  \cos\left( \tfrac{2\pi pqm}{r}\right)
-\tfrac{t}{3} \right) +O(t^2).
\end{split}
\end{equation}

This characterization
of the longitudinal and meridinal holonomies of
the perturbed flat irreducible connections near $A_0$ allow us to
prove the following theorem, which describes the perturbed flat moduli
space of $\Si$ in a neighborhood of the pointed 2-sphere.
\begin{thm} Let $S\subset \cM(\Si)$ be a pointed 2-sphere, and let
$[A_0]\in S$ be the gauge orbit of reducible connections.
   For a sufficiently small neighborhood $\cU \subset
\cB(\Si)$ of $S$, and for sufficiently small $t>0$,
$\cU \cap \cM_{tf}(\Si)$ consists of two components.  The first is an
isolated gauge orbit of reducible connections, and the second is a
smooth 2-sphere of gauge orbits of irreducible connections.
\end{thm}
\begin{remark}  In this theorem we do not claim that the reducible
connection $[A_t]\in \cM_{tf}(\Si)$ near $[A_0]$ satisfies the
nondegeneracy condition $H^1 _{A_t, tf} (\Si; su(3))=0$.  This will be
proved in Proposition  \ref{noH1forsmallt}.
\end{remark}

\begin{proof}
Choose a neighborhood $\cU_Z$ of $[A_0|_Z]$ in $\cB_Z$ with the
following properties:
\begin{itemize}
\item[(i)] $\cU_Z \cap \cM^* (Z) \subset \cC,$ where $\cC$ is
the 4-dimensional Type II component of $\overline{\cM ^* (Z)}$ containing
$[A_0 |_Z]$, as in Theorem \ref{Zsthm}.
\item[(ii)] $\cU_Z \cap \cM^\red (Z) \subset \cC^\red,$ where $\cC^\red$ is  the 2-dimensional
component of $\cM^\red (Z)$ containing $[A_0|_Z]$, as in Theorem
\ref{seam}.
\item[(iii)] $ r^{-1} (\cU_Z) \cap \cM (\Si) = S,$ where $r\colon \cB(\Si) \to \cB(Z)$
is the restriction map.
\item[(iv)] The  restriction of $\zeta^\red_Z$ to
$\cC^\red \cap \cU_Z  $  is injective.
\end{itemize}
Set $\cU = r^{-1}(\cU_Z).$
The intersection $\cM _{tf} ^\red (\Si) \cap \cU$ is identified with
$$
\{ [A] \in \cC^\red\cap \cU_Z \mid [A|_T ] \in \cZ_Z ^\red \cap Z_{Y,tf}
\}.$$
This intersection is a single point, identified in Proposition
\ref{redps} and the restriction map $\zeta^\red_Z$ maps
$\cC^\red\cap \cU_Z $ injectively into $\cM(T)$.
Thus $\cM_{tf} ^\red (\Si) \cap \cU$ is a single point.

Now consider $\cM _{tf} ^* (\Si) \cap \cU$, which is identified with
$$
\{ [A] \in \cC \cap \cU_Z \mid [A|_T]\in \cZ_Z ^* \cap Z_{Y, tf}
\}.$$
In equations (\ref{holy}) we have identified the unique point in
$\zeta_Z (\cC \cap \cU_Z  ) \cap Z_{Y,tf}.$  This point has a 2-sphere
preimage in $\cC \cap \cU_Z$ for small $t$, because $\cC$ is topologically
a 2-sphere bundle.  
This can be seen by observing that the map $\zeta_Z : \cC \cap \cU_Z 
\to \cM(T^2) $ factors through $\xi\colon \cC \cap \cU_Z \to \De,$ which
sends $\al$ to $[\al(xy)],$ because $\al(h)=e^{2\pi i \ell/3}I, 
\la=(xy)^{pq} h^{-(p+q)a}$  and $\mu=(xy)^{r}h^c.$ 
Again by equations (\ref{holy}), the longitudinal
holonomy does not have 1 as an eigenvalue, and hence the 2-sphere fiber
does not contain any reducibles, so by Theorem \ref{fiber}
it is a smooth 2-sphere of gauge orbits of irreducible connections.
\end{proof}

\section{Spectral flow arguments} \label{s5}
In this section, we perform computations of the spectral
flow of the odd signature operator.
These are necessary to calculate the contribution of
the pointed 2-spheres to the invariant $\tau_{SU(3)}(\Si)$.
The main result here is that, given a path
$A_t$ of reducible $(tf)$-perturbed connections
on $\Si$ where $[A_0]$ is flat and lies on a 2-sphere,
the $\CC^2$ spectral flow of the perturbed odd signature
operator equals $SF_{\CC^2}(A_t;\Si) = -2.$
This is proved by splitting the spectral flow
according to the manifold decomposition $\Si = Y \cup_T Z$
(Theorem \ref{SFST}),
and then computing the spectral flow on $Z$
(Theorem \ref{KCSF}).

\subsection{The odd signature operator, spectral flow, and splittings} \label{ap}
As in Section \ref{s4} we assume that $\Si=\Si(p,q,r)$ is endowed with a metric isometric to the
product metric on a bicollared neighborhood $[-1,1] \times T$, where $\Si = Y \cup_T Z$.

The  operator $D_A$ is a self-adjoint  Dirac-type operator. Thus on the closed manifold
$\Si(p,q,r)$, $D_A$ has a compact resolvent and hence the spectrum of $D_A$ is unbounded but
discrete, and each of its eigenspaces is finite dimensional. Although $D_{A,h}$ is not a Dirac-type
operator, it is a compact perturbation of  $D_A$ and also has  a compact resolvent.

Given a suitably continuous path $D_t, \ 0 \leq t \leq 1,$ of self-adjoint operators with discrete,
real spectrum
  each of whose eigenspaces is finite dimensional,
one can define the spectral flow $SF(D_t)\in\ZZ$ to be the algebraic
intersection
in $[0,1] \times \RR$ of
the track of the spectrum
$$\{(t,\la) \mid t \in [0,1], \ \la \in \Spec(D_t)\}$$
with the line segment from $(0,-\varep)$ to $(1,-\varep)$,
where $\varep>0$ is chosen smaller than the modulus of the
largest negative eigenvalue of $D_0$ and of $D_1$ (this is called the
$(-\varep,-\varep)$ convention).

If $A_t$ is a continuous path of $SU(3)$ connections on the closed 3-manifold $X$  and $h_t$ a
continuous path of perturbations, we denote by $SF(D_{A_t,h_t};X)$ or $SF(A_t, h_t;X)$ the spectral
flow of the family of odd signature operators $D_{A_t, h_t}$ on $\Om^{0+1}(X;su(3))$.  (A proof
that the family $D_{A_t}$ is suitably continuous and a careful definition of the spectral flow can
be found in \cite{BLP} and \cite{HKL}.) The spectral flow is an invariant of homotopy rel
endpoints, and to emphasize this point we will occasionally write $SF(A_0,A_1;X)$ instead of
$SF(D_{A_t, h_t};X)$ when the path of perturbations is understood (the parameter space of pairs
$(A,h)$ is contractible).

If $A$ is an $S(U(2) \times U(1))$ connection on $X$, then $D_{A,h}$ respects the decomposition on
forms induced by the splitting of coefficients $su(3) = s(u(2)\times
u(1)) \oplus \CC^2.$
 In particular, for a path $A_t$ of $S(U(2) \times U(1))$
connections and path of perturbations $h_t$, we denote by $SF_{\CC^2}(A_t, h_t;X)$  the spectral
flow of the restriction of the path $D_{A_t,h_t}$ to $\Om^{0+1}(X;\CC^2).$ Similar
notation applies to the other summand in this decomposition of $su(3)$.

In computing the $\CC^2$ spectral flow, we count eigenvalues with their real multiplicity, thus
$SF_{\CC^2}(A_t, h_t;X)$ is always a multiple of two and we have $$SF_{su(3)}(A_t, h_t;X) =
SF_{s(u(2)\times u(1))}(A_t, h_t;X) + SF_{\CC^2}(A_t, h_t;X).$$

 \medskip

   When $X$ is compact
but has nonempty boundary $\partial X = W$ the constructions must be refined in order to obtain
suitable families of operators for which one can define the spectral flow. We must draw on deeper
results from the Calder\'on-Seeley theory of boundary-value problems for Dirac operators.

Assume the metric on $X$  is isometric to the product
  metric on  a
  collar $W \times (-1,0]$ of the boundary $\partial X = W \times \{0\}$.
  We work with connections $A$ on $X$ that are in  {\it cylindrical
  form}, namely we assume that the restriction of $A$ to
  the collar $W \times (-1,0]$ is the pullback  of a connection $a$ on
$W$
  under the natural projection $W \times (-1,0] \to W$.

  Given an $su(3)$ connection $a$ on $W,$
  define the de Rham operator
\begin{eqnarray*}
&S_a\colon \Om^{0+1+2}(W;su(3)) \lto \Om^{0+1+2}(W;su(3))& \\
&S_a(\al,\be,\ga)
  = (* d_a\be, - * d_a \al - d_a * \ga, d_a * \be).&
  \end{eqnarray*}
  Here, $*\colon  \Om^i(W;su(3)) \to \Om^{2-i}(W;su(3))$
denotes the Hodge star operator on $W$.
Define $P^\pm_a$ to be the positive and negative
eigenspans of this operator on the space of $L^2$ forms
$L^2 ( \Om^{0+1+2}(W;su(3))).$

If $a$ is a flat connection on $W,$
then the Hodge and de Rham theorems identify
  the kernel of $S_a$ with the cohomology groups
  $H^{0+1+2}_a (W;su(3))$ with coefficients in the local system $su(3)$
  twisted by $a$.
Define the operator
\begin{eqnarray*}
& J\colon \Om^{0+1+2}(W;su(3)) \lto \Om^{0+1+2}(W;su(3)) \\
&J(\al,\be,\ga) = (- * \ga, * \, \be, * \, \al).&
\end{eqnarray*}
Notice that $J^2 = -1$. Setting
$w(x,y) = \langle x,J y\rangle_{L^2}$ defines
a symplectic structure on
the Hilbert space $L^2(   \Om^{0+1+2}(W;su(3)))$
of $L^2$ forms.   By
restricting this also gives a symplectic structure to $\ker S_a$.

If $A$ is an $SU(3)$ connection on $X$ in cylindrical form,
and $a$ is its restriction to the boundary $\partial X = W,$
then along the collar
$W \times [-1,0],$ we have
\begin{equation} \label{neck}
D_A =  J\left(S_a + \tfrac{\partial}{\partial s}\right),
\end{equation}
where $s $ denotes the collar coordinate. (See Lemma 2.4 of \cite{BHKK}.) This holds more generally
for $D_{A,h}$ provided the    perturbation is supported away from  the collar. Given a Lagrangian
subspace $L\subset \ker S_a$, the operator $D_{A,h}$ taken with domain those $L^2_1$ sections
$\phi\in \Om^{0+1}(X;su(3))$  satisfying the APS boundary condition $$\phi |_{W}\in L\oplus P^+_a$$
 is self-adjoint with compact resolvent and hence discrete
spectrum. Given a family $(A_t,h_t)$ and a   choice of  Lagrangian subspaces $L_t\subset \ker
S_{a_t}$ so that $L_t\oplus P^+_{a_t}$ is continuous, the spectral flow $SF(D_{A,h},P^+_a)\in \ZZ$
is well defined (see e.g.~\cite{BLP}). In our context below we will have $\ker S_{a_t}=0$ for all
$t$ and $P^+_{a_t}$ continuous.

Given a  connection $A$ on $X$ in cylindrical form,and $h$ a
perturbation of the type we described above
we define an (infinite-dimensional) Lagrangian subspace
$$\La_{X,A,h}\subset L^2\left( \Om^{0+1+2}(W;su(3))\right)$$ as
follows.The main result of  \cite{HKL} implies that
  there is a well-defined injective  map
$$r \colon  \ker \left(D_{A,h}\colon  L^2_{1/2}\left(
\Om^{0+1}(X;su(3)\right) \to
  L^2_{-1/2}\left( \Om^{0+1}(X;su(3)\right) \right)
  \to L^2\left(   \Om^{0+1+2}(W;su(3))\right)$$
  given by restriction whose image is a closed,
  infinite dimensional Lagrangian subspace  called the
{\it Cauchy data space} of the operator $D_{A,h}$ on $X$
and is denoted $\La_{X,A,h}.$ Since the restriction map
$r$ is injective,  kernel of $D_{A,h}$ with $P^+_a$  (i.e.~APS)
boundary conditions is isomorphic to  $\La_{X,A,h}\cap P^+_a$.
  When the context is clear, we will abbreviate $\La_{X,A,h}$ to
$\La_{X,A}$ or even $\La_X$.

  The space $\La_{X,A,h}$ varies continuously (in the graph topology
on closed subspaces) with respect to $A$, $h$, and the metric on $X$.
This result is well known in the case of Dirac-type operators (such as
$D_A$), see e.g.~\cite{BW}. The theorems of the article \cite{HKL}
  extend   these standard results to the more  general setting of
  small
perturbations of Dirac operators such as $D_{A,h}$ (which is not
a differential or even a pseudodifferential operator).

\begin{remark}
The previous remarks change slightly when the collar of $\partial
X$ is parameterized as $[0,1)\times W$ with $\del X=\{0\}\times
W$. The significant difference is that the positive eigenspan
$P^+_a$ of $S_a$ is replaced by the negative eigenspan $P^-_a$.
\end{remark}

We will apply these observations to the decomposition
$\Sigma=Y\cup_T Z$. Parameterize a collar of the separating torus $T$
as $(-1,1)\times T$ in $\Sigma$, with $(-1,0]\times T$ a
collar of the boundary of the solid torus $Y$ and $[0,1)\times T$
a collar of the boundary of $Z$.

 The fact that the operator $D_{A,h}$ on $\Si$ is Fredholm is
equivalent to the fact that the pair $(\La_{Y,A,h},\La_{Z,A,h})$
form a Fredholm pair of (Lagrangian) subspaces, and hence if $(A_t,
h_t)_{t\in[0,1]}$ is a path, the Maslov index
$\Mas(\La_{Y,A,h},\La_{Z,A,h})$ is well defined. Similarly the
restriction of
$D_{A,h}$ to $Y$ with $P^+_a$ boundary conditions is Fredholm
because the pair of subspaces $(\La_{Y,A,h},P^+_a)$ is Fredholm,
and
the
restriction of
$D_{A,h}$ to $Z$ with $P^-_a$ boundary conditions is Fredholm
because the pair of subspaces $(P^-_a,\La_{Z,A,h})$ is Fredholm.
Proofs of these facts can be found e.g.~in \cite{Nico} or
\cite{KL}.

\subsection{Some vanishing results}
\label{SVR}
This subsection consists of an interlude to prove some needed vanishing
results
for the perturbed flat cohomology groups.
To begin with, we note the following property of perturbed flat
cohomology.  The proof  is the same as the standard proof of the exactness of the
Mayer-Vietoris sequence and is left as
an exercise. Note that the restriction of $A$ to $Z$ is
flat and so $H^*_{A,tf}(Z;\CC^2)=H^*_{A}(Z;\CC^2)$ and similarly
for $T$.
\begin{lem} \label{MV}
If $A$ is a $(tf)$-perturbed flat connection on $\Si$, the
Mayer-Vietoris sequence
\[ \cdots \to H^0_A(T;\CC^2)\to H^1_{A,tf}(\Si;\CC^2)\to
H^1_{A,tf}(Y;\CC^2)\oplus H^1_{A}(Z;\CC^2) \to
H^1_{A}(T;\CC^2)\to\cdots\]
is exact.
 \end{lem}

To use the Mayer-Vietoris sequence in the present context, we need to
know the perturbed flat cohomology of the perturbed flat connections
on $Y$.  This information is provided by the following lemma.

\begin{lem}  \label{pfcv}
For $0 < \de <  \frac \pi 2$, define the open rectangle
$$R_\de = \{ (u,v) \mid  \de <u< 2\pi-\de, \ -\de/4< v < \de/4 \}.$$
Given $0 < \de <  \frac \pi 2$, there
exists an $\ep > 0$
such that, if $-\ep<t<\ep$ 
then $H^0 _A(Y;\CC^2)=0$ and $H^1 _{A,tf}(Y;\CC^2)=0$
for every
$(tf)$-perturbed flat connection $A$ on $Y$
with
$\hol_\la(A)=\Phi(u,v)$  for $(u,v) \in R_\de$.
\end{lem}

\begin{proof}
Fix $0 < \de < \frac \pi 2$ and
consider the subset $\cM_\de(Y)$ of
$\cM(Y)$ consisting of gauge orbits of flat connections
  $A$  with
$\hol_\la(A)$ conjugate to $\Phi(u,v)$
for $(u,v)$ in the closure of $R_\de.$
Obviously $\cM_\de$ is a compact subset of $\cM$.
Moreover,  the conditions on $(u,v)$ guarantee that
$\hol_\la(A)$ acts nontrivially on $\CC^2$ for all
$[A] \in \cM_\de.$  From this, it
follows that $H^0_A(Y;\CC^2)  = 0$
for all $[A] \in \cM_\de.$
Poincar\'e duality on the circle (a retract of $Y$) then
gives $H^1_A(Y;\CC^2)=0$ as well.

On the  closed manifold $\Si$, if $A$ is a flat connection, then
  one may identify the cohomology $H^p _A(\Si; su(3))$ with  the kernel
of the
  operator $d_A \oplus d_A^*\colon  L^2 _1 \Om^p(\Si;su(3)) \to L^2
  \Om^{p+1}(\Si; su(3)) \oplus L^2 \Om^{p-1}(\Si; su(3))$, which
  is elliptic and hence Fredholm.  On the manifold $Y$, with non-empty
  boundary, one must impose Neumann boundary conditions for this to be
  an elliptic operator, namely replace the domain by
$$
L^2 _1 \Om^p _\nu (Y; su(3)) = L^2 _1 \{ \al \in \Om^p (\Si;
su(3))\mid \star \al |_{T} =0 \}.$$

The map $D_A$ is equivalent to the sum of the de Rham operator and its
adjoint from odd forms to even forms, except that we have used the
Hodge star operator to replace 3-forms by 0-forms and 2-forms by
1-forms.  Hence the appropriate Dirichlet/Neumann-type boundary conditions for
$D_A$ are to restrict the domain to
$$L^2 _1 \Om^{0+1} _{\tau, \nu} =
\{ (\al , \be) \in L^2 _1 \Om^{0+1}(Y;su(3)) \mid \al|_T = 0,
\star \be|_T=0 \}.$$  If $A$ is not flat, then this operator $D_A$
differs from that of a flat connection (for example, the trivial
connection) by a compact operator (see \cite{T}).
As pointed out above,  the  operator
$D_{A,tf}$  also  differs from
$D_A$ by a compact operator and hence, with these boundary conditions,
  is still Fredholm.   Again the (perturbed) cohomology $H^0 _A
(Y;su(3)) \oplus H^1 _{A,tf} (Y;su(3))$ of a $tf$-perturbed flat
connection is identified with
the kernel of this operator with the restricted domain.

For flat connections $A$, we have $[A] \in \cM_\de$,
  the kernel of $D_A$ restricted to $L^2 _1 \Om^{0+1} _{\tau,\nu}
  (Y;su(3))$ equals $H^{0+1}(Y;\CC^2),$
which vanishes for $[A] \in \cM_\de$   by the previous argument.
Using upper semicontinuity of the dimension of the kernel of a
  continuous family of Fredholm operators,
  the family  $D_{A,tf}$, with the same
  boundary conditions, must have trivial kernel
  neighborhood
of  $([A_0],0)$ for fixed $[A_0] \in \cM_\de$.
Using compactness of $\cM_\de,$ we obtain an $\ep
0$ such that if $A$ is $(tf)$-perturbed flat for $-\ep< t < \ep$ and if $\hol_\la(A) = \Phi(u,v)$
for $(u,v) \in R_\de,$ then  $H^0_A(Y;\CC^2)$ and $H^1_{A,tf}(Y;\CC^2)$ vanish.
\end{proof}

As in Section \ref{s4},
suppose $A_0$ is a reducible flat connection on $\Si$ whose gauge orbit
$[A_0]$ lies on a 2-sphere component.
For $0 \leq  t \leq \ep,$ let $A_t$ be the family constructed
in Subsection \ref{TEOT} of reducible
$(tf)$-perturbed flat
connections on $\Si$ limiting to $[A_0]$ as $t\to 0$.

\begin{prop} \label{noH1forsmallt}
If $t>0$ is sufficiently small,   then $H^1 _{A_t, tf}(\Si;su(3))=0$.
\end{prop}
\begin{proof}
We split the coefficients according
to the decomposition $su(3) = s(u(2) \times u(1)) \oplus \CC^2$
and argue the two cases separately.  The fact that $H^1 _{A_0}(\Si;
s(u(2)\times u(1))) = 0 $ implies that the same holds true for the
perturbed cohomology for small $t$.  As far as the $\CC^2$ cohomology
goes, we cannot make the same argument since $H^1 _{A_0} (\Si; \CC^2)
= \CC^2$.  Instead, we combine
Corollary  \ref{Z cohom goes away} and Lemma \ref{pfcv}, using the
Mayer-Vietoris sequence, to obtain the desired conclusion.
\end{proof}

\subsection{The spectral flow to the reducible perturbed flat
connection}
  We turn now to an analysis of the spectral flow   from the
reducible flat connection whose orbit lies on a pointed 2-sphere to the nearby reducible perturbed
flat  connection. The set-up is as follows. We have a path  $A_t$  of reducible $(tf)$-perturbed
flat connections on $\Si$ such that $A_0$ is a flat connection whose gauge orbit lies on a 2-sphere
component. In Theorem \ref{KCSF} we compute the spectral flow $$SF_{\CC^2}(A_t,tf; \, \Si; \, 0
\leq t \leq \ep)$$ of the perturbed odd signature operators $D_{A_t,tf}\colon  \Om^{0+1}(\Si;
\CC^2) \to \Om^{0+1}(\Si; \CC^2)$ from $t=0$ to $t=\ep$. The strategy is to use the machinery of
Cauchy data spaces to prove a splitting result for spectral flow. This is accomplished in Theorem
\ref{SFST} which which shows that the spectral flow is concentrated on $Z$. The path $A_t$
restricts to a path of flat connections on $Z$ which allows us to compute the the resulting
spectral flow by topological methods. (For the remainder of this subsection, we restrict
$D_{A_t,tf}$ to $\CC^2$ valued forms and write $SF$ for $SF_{\CC^2}$ without further reference.)

As before we let $a_t$ denote the path of flat connections on the separating torus $T$ in the
decomposition (\ref{splitSi}) and let $S_{a_t}$ be the corresponding path of of twisted de Rham
operators on $\Om^{0+1+2}(T; \CC^2)$.
  Since the twisting perturbation is supported on the interior
  of the solid torus and vanishes on the neck,  it follows that
  the operators $D_{A_t,tf}$ and $D_{A_t}$
  coincide on  $([-1,0]\times T) \cup Z.$
Thus on the neck, equation (\ref{neck}) gives that
  \begin{equation}\label{neck2}D_{A_t,tf} =J(S_{a_t}+\tfrac{\partial}{\partial s}).
  \end{equation}

Let $P^\pm_t$ denote the positive and negative
eigenspans  of the operator
$S_{a_t}$.
  Denote by $\La_Y(t)\subset L^2(\Om^{0+1+2}(T;\CC^2))$ the Cauchy data
space of
the operator $D_{A_t,tf}$ on $Y$ and by $\La_Z(t)$ the
Cauchy data space of
  $D_{A_t,tf}$ on $Z$.
Thus the
kernel of
$D_{A_t,tf}$
is isomorphic to the intersection
$\La_Y(t)\cap \La_Z(t)$.

Let $Y^R$ be $Y$ with a collar of length $R$ attached, namely $$Y^R = Y \cup ([0,R] \times T).$$
Any connection  $A \in \cA(Y)$ in cylindrical form extends in the obvious way  to give a connection
on $Y^R$ in cylindrical form. Thus the family $D_{A_t,tf}$ of perturbed odd signature operators on
$Y$ extends (using (\ref{neck2})) to give a family of operators on $Y^R$. Let $\La^R_{Y}(t)$ denote
the Cauchy data space of the operator $D_{A_t, t f}$ on $\Om^{0+1}(Y^R;\CC^2)$. Similarly, set
$Z^R=([-R,0] \times T) \cup  Z$ and denote by $\La^R_{Z}(t)$  the Cauchy data space of the operator
$D_{A_t, t f}$ on $\Om^{0+1}(Z^R;\CC^2)$.

\begin{lem}\label{geah}
There exists an $\ep>0$ such that  $0\leq t \leq \ep$ implies
\begin{enumerate} \item[(i)]  $\ker S_{a_t}=0$.
\item[(ii)] $\La_Y^R(t)\cap P^+_t=0$ for all $R>0$.
\item[(iii)]  $ {\displaystyle \lim_{R\to \infty}} \La_Y^R(t)=P^-_t$.
\item[(iv)] $\La_Z^R(\ep) \cap P^-_\ep=0$ for all $R>0$.
\item[(v)] ${\displaystyle \lim_{R\to \infty}} \La_Z ^R (\ep)= P^+_\ep.$
\end{enumerate}
\end{lem}
\begin{proof}  As in Subsection \ref{TEOT},
the reducible flat connection $A_0$ has longitudinal
holonomy $\hol_{\la} (A_0) = \Phi(\tfrac{2 \pi pqk}{r},0)$
for some $0< k < r.$ The matrix $\Phi(\tfrac{2 \pi pqk}{r},0)$
acts nontrivially on $\CC^2$, and it follows
that  $H^0_{a_0} (T; \CC^2)=0$.
Poincar\'e duality implies $H^2_{a_0} (T; \CC^2)=0,$
and Euler characteristic
considerations give that $H^1_{a_0}(T;\CC^2)=0$ as well.
Hence
$$\ker S_{a_0} = H^{0+1+2}_{a_0} (T;\CC^2) = 0 .$$
By upper semicontinuity, $\ker S_{a_t}=0$ for small $t$.
This proves (i).

Proposition 2.10 of \cite{BHKK} states that if $A$ is a flat
   $S(U(2)\times U(1))$
connection on a 3-manifold $X$ with boundary, and if
$a=A|_{\partial X}$, then $\La_{X,A} \cap P^+ _a$ is isomorphic to the
image of the relative cohomology in the absolute
$$\Image \left(H^1_A (X,\partial X; \CC^2) \to H^1_A (X; \CC^2)
\right).$$
The proof involves identifying the intersection with the space of
$L^2$ harmonic forms on the infinite cylinder and  applying
Proposition 4.9 of \cite{APS}.    If $H^{0+1+2}_a(\partial X;\CC^2)=0$,
then the image of the relative cohomology  in the
absolute is exactly $H^1_A (X; \CC^2)$.

Apply this result to the case $A=A_0$ and $X=Y^R$.
Since $H^1_{A_0} (Y^R;\CC^2)=0$ (by Lemma \ref{pfcv}), we conclude that
$\La _Y ^R (0 )\cap P^+_0=0$ for all $R$.
This generalizes to perturbed flat connections as
follows.  The proof of \cite{APS} that
the space  of $L^2$ harmonic forms injects into $H^1_A (X;\CC^2)$
works just as easily to show that the space of
$L^2$ solutions
to $D_{A_t,tf}(\si,\tau)=0$ on $Y^\infty$  injects into
$H^1 _{A_t, tf} (Y;\CC^2)$.  But Lemma \ref{pfcv}
shows that $H^1 _{A_t, tf} (Y;\CC^2) = 0$. This
proves (ii).

Assertion (iv) follows by applying the same argument to the case
   $A=A_\ep$ and $X=Z$. Note that Proposition
\ref{redps}  implies  that, for $\ep>0$ small enough,
$H^1_{A_\ep} (Z;\CC^2) =0$.

Assertion (iii)  follows from (i) and (ii)
and a theorem of Nicolaescu
(\cite[Corollary 4.11]{Nico}; see Theorem 2.7 of \cite{BHKK} for the
result in the present context).
Similarly, Assertion (v) follows from (iv).
   \end{proof}

The restriction of the operator $D_{A_t,tf}$  to $Z$ coincides with $D_{A_t}$ on $Z$. The operator
$D_{A,tf}$ restricted to those $L^2_1$ sections whose restriction to the boundary lie in  $P^-_t$
is a well-posed  elliptic boundary value problem which, furthermore, is self-adjoint since
 $\ker
S_{a_t}=0$. This implies that the spectral flow $SF(D_{A_t, tf}; \, Z; \, P^-_t)$ is well-defined.
(These are well-known facts, originating in \cite{APS}, whose proofs can be found in many places,
e.g.~\cite{KL}.)

The next result is a splitting theorem which uses the vanishing of cohomology on the solid torus
$Y$ to localize the spectral flow on the knot complement $Z$.

\begin{thm}\label{SFST}
For small $\ep>0$, $$SF(D_{A_t,tf}; \,  \Si; \,  0 \leq t \leq \ep)= SF(D_{A_t}; \, Z; \, P^-_t; \,
0 \leq t \leq \ep).$$
\end{thm}
\begin{proof}
By part (i) of Lemma \ref{geah}, we have
$H^{0+1+2}_{a_t} (T;\CC^2)=0$  for $0 \leq t\leq\ep$.
A theorem of Nicolaescu (\cite{Nico}; see also \cite{KL}) states that
\begin{equation}\label{maseq0}SF(D_{A_t, tf};
\Si)=\Mas(\La_Y(t),\La_Z(t)). \end{equation}
   As in
\cite{BHKK} and \cite{DK},
we use homotopy invariance and additivity of the Maslov index
to complete the argument.
(For a precise definition of the Maslov index in this context,  see
\cite{Nico,KL} and \cite[Definition 2.13]{BHKK}).

Consider the 2-parameter family
$$L(s,t)= \begin{cases}
\La^{1/(1-s)}_Y(t) & \text{for $0 \leq s<1$} \\
P^-_t& \text{if $s=1$} \\
\end{cases}$$ for $0\leq s\leq 1, 0\leq t\leq \epsilon$.  Lemma
\ref{geah}(iii)  and the appendix to \cite{DK} shows that for each
fixed $t$ this is a continuous path. What we need is uniform continuity
in the $t$ parameter.
Such families are not always continuous (see \cite{BHKK} for a
discontinuous example) but in this case the family is continuous by
\cite[Corollary 4.12]{Nico}.  The  required {\em nonresonance}
hypothesis  is exactly what Lemma \ref{geah} (ii) asserts.

Since $L(0,t)=\La_Y(t)$ and $L(1,t)=P^-_t$, additivity and homotopy
invariance of the Maslov index implies that
\begin{equation}\label{maseq1}
\Mas(\La_Y(t),\La_Z(t))=\Mas(L(s,0),\La_Z(0))+\Mas(P^-_t,\La_Z(t))
-\Mas(L(s,\epsilon),\La_Z(\epsilon)).
\end{equation}

Since $A_0$ is flat, Proposition \ref{coho} shows that, for $0 \leq
s<1$,
$$\dim (L{(s,0)}\cap \La_Z(0))
=\dim \ker D_{A_0} =\dim(H^{0+1}_{A_0}(\Si;\CC^2))=4.$$
(Note, all dimensions computed here are real.)
For $t=1$,
$$L{(1,0)}\cap \La_Z(0)=P^-_0\cap \La_Z(0)  \cong
\Image \left(H^1_{A_0}(Z,T;\CC^2)\to H^1_{A_0}(Z;\CC^2)\right).$$
  Since
$H^{0+1+2}_{a_0}(T;\CC^2)=0$,  the image of the
relative cohomology in the absolute is all of $H^1_{A_0}(Z;\CC^2)$
which has complex dimension
$2$
by Proposition \ref{coho-Z}.
Thus
$\dim(L{(s,0)}\cap \La_Z(0))$ is constant in $t$
and it follows that
\begin{equation}\label{maseq2}\Mas(L(s,0),\La_Z(0))=0.\end{equation}

By Part (iv) of Lemma \ref{geah}, $L(1,\ep)\cap \La_Z(\ep)=0$.
For $0\leq s<1$, we have
\begin{equation*}
\begin{split}
\dim \left( L(s,\ep)\cap \La_Z(\ep) \right) &= \dim \ker \left( D_{A_\ep,\ep f}\colon
\Om^{0+1}(\Si^R; \CC^2) \to \Om^{0+1}(\Si^R;\CC^2)\right)\\ & = \dim  H^{0+1}_{A_\ep, \ep f} (\Si;
\CC^2) =0.
\end{split}
\end{equation*}

Here, $\Si^R = Y^R \cup_T Z$ is the result of adding a collar
of length $R$ to the neck.
The computation that
$H^{0+1}_{A_\ep, \ep f} (\Si; \CC^2)=0$ follows by
a Mayer-Vietoris argument,
using Lemma \ref{pfcv} and Proposition
\ref{TFAE-Z}.  Therefore
\begin{equation}\label{maseq3}\Mas(L(s,\epsilon),\La_Z(\epsilon))=0.
\end{equation}

Next,
\begin{equation}\label{maseq4}\Mas(P^-_t, \La_Z(t)) = SF({A_t}; \, Z;
\,  P^-_t; \, 0\leq t \leq \ep)\end{equation}
(This result is also due to Nicolaescu; see  \cite{KL} and
\cite[Theorem 2.18]{BHKK} for    proofs  in the present context).
  Combining (\ref{maseq0}), (\ref{maseq1}), (\ref{maseq2}),
(\ref{maseq3}),  and (\ref{maseq4}) with the observation that $D_{A_t,tf}$ and $D_{A_t}$ agree on
$Z$ completes the argument.
\end{proof}

Theorem \ref{SFST} reduces the problem of computing   $SF_{\CC^2}(D_{A_t,tf};\Si)$ from the flat
irreducible connection $A_0$  and zero perturbation to the $\ep f$-perturbed flat reducible
connection $A_\ep$ and perturbation $\ep f$ to  the problem of computing
  the spectral flow on the knot complement,
namely,  $SF_{\CC^2}(D_{A_t}; \, Z; \, P^-_t)$. This is a much easier problem for the following
reason. The  path of perturbed-flat connections $A_t$ restricts to a path of {\em flat} connections
on $Z$, and the kernel of $D_{A_t}$ acting on $\CC^2$-valued forms with boundary conditions $P^-$
is isomorphic to the image of $H^1(Z,T;\CC^2_{\al_t})\to H^1(Z;\CC^2_{\al_t})$ (see the proof of
Lemma \ref{geah}). Corollary \ref{Z cohom goes away} then implies that this kernel is $0$ for
$t>0$, and Proposition \ref{yum} shows that the kernel is $\CC^2\cong \RR^4$ for $t=0$.  We will
prove that two zero modes become positive and two become negative, so that the spectral flow equals
$-2$ (with our conventions). The homotopy will be a disk in the cylinder $S^1\times \RR$ of Theorem
\ref{seam}.

\begin{thm} \label{KCSF}
With $\ep>0$ as in Theorem \ref{SFST}, we have
$$SF(D_{A_t}; \, Z; \, P^-_t; \, 0 \leq t \leq \ep)=-2.$$\end{thm}
\begin{proof}
As mentioned above, for $t=0,$  the kernel of $D_{A_0}$ with $P^-$ boundary conditions has real
dimension 4, but for $t>0$, the kernel is trivial.

In Subsection \ref{dot}, we constructed  2-parameter families of reducible $SU(3)$ representations
on $Z$. These results give 2-parameter families of based gauge orbits of flat connections on $Z$.
The based gauge group is the subgroup of $\cG(Z)$ consisting of those gauge transformations in the
path component of the identity. The point is that spectral flow is a well defined concept for
connections modulo based gauge transformations, so we can use the parameterization from Subsection
\ref{dot} to compute spectral flow.

If needed,  gauge transform the path $A_t$ so that its path of holonomy representations
$\gamma_t\colon \pi_1(Z)\to SU(3)$ takes values in $S(U(2)\times U(1))$ and so that $xy$ is sent to a
diagonal matrix. Notice that $\gamma_0$ takes values in  $SU(2)\times \{1\}$  since $A_t$ is the
restriction of flat connection on $\Si(p,q,r)$. Thus $\gamma_0$ lies on an arc $\al$  (see
Definition \ref{defnofalsth}) for some $k,\ell,$ and $ \ep$. The precise values of $k,\ell,$ and
$\ep$ are not needed for our argument.

Suppose that $\gamma_0=\al_{s_0}$ for some $s_0\in(0,1)$.  Proposition \ref{pfhc} (in  particular
(\ref{tnoia})) shows that $\gamma_t$ lies off the seam of $SU(2)\times 1$ representations for
$t>0$, and in particular $\gamma_t$ is a $S(U(2)\times U(1)$ representation but not an
$SU(2)\times\{1\}$ representation for $t>0$. Hence Theorem \ref{seam} implies that $\gamma_t$ is of
the form $\al_{s_t,\th_t}$ for paths $s_t\in(0,1)$ and $\th_t\in [0,\pi]$. (We assume $\ep$ is
small so that $\th_t$ is also small.)

Now the construction of Definition \ref{defnofalsth} gives a 2-parameter family of representations:
namely the disk in the cylinder bounded by union of the 4 curves (see Figure 3):
\begin{enumerate}
\item[(i)] $\gamma_t=\al_{s_t,\th_t}, t\in [0,\ep]$,
\item[(ii)] $\al_{s_\ep,(1- u)\th_\ep}, u\in [0,1]$,
\item[(iii)] $\al_{(1-u)s_\ep,0}, u\in [0,1]$,
\item[(iv)] $\al_{u, 0}, u\in [0, s_0]$.
\end{enumerate}

  \begin{figure}[hbt]
\begin{center}
\leavevmode\hbox{}
\includegraphics[width=2.8in]{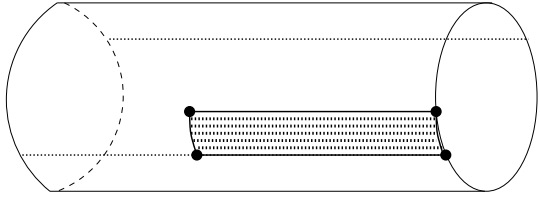}
\caption{The family $A_{s,t}$ is drawn as the shaded region, a subset of  one of the cylinders in
$R^\red(Z,SU(3))$. The seams are indicated by the two dotted lines which span the cylinder lengthwise.
Taken clockwise from upper left hand corner, the four vertices are the flat
connections $A_{0,\ep}, A_{1,\ep}, A_{1,0}$ and $A_{0,0}$.}
\end{center} 
\end{figure}

This disk determines a 2-parameter family of reducible flat connections $$\{A_{s,t} \mid 0 \leq s
\leq 1, 0 \leq t \leq \ep\}$$ such that:
\begin{enumerate}
\item[(i)] $A_{0,t} = A_t$ for $0 \leq t \leq \ep.$
\item[(ii)] $A_{s,\ep}$ is a flat $S(U(2)\times U(1))$ connection
with $H^1_{A_{s,\ep}}(Z; \CC^2) =0$ for $0\leq s < 1$ (see Lemma
\ref{pfcv}).
\item[(iii)] $A_{1,t}$ is a flat abelian connection for $0 \leq t \leq
\ep$ and $H^1_{A_{1,t}}(Z; \CC^2) =0$ for $0 < t \leq \ep$.

\item[(iv)] $A_{s,0}$ is a flat $SU(2)\times\{1\}$ connection
with $H^1_{A_{s,0}}(Z; \CC^2) = \CC^2$ for $0\leq s \leq 1.$
\end{enumerate}
The parameterization in $s$ and $t$ may be chosen so that, when $s$ is near 1, the $t$ parameter is
simply twisting, $\hol_\gamma (A_{s,t})$ equals the twist of $\hol_\gamma (A_{s,0})$ by  the
character sending $\mu$ to $e^{it}$. This family parameterizes a thin strip on
  the cylinder $S^1 \times \RR$ with the edge corresponding to (iii) in the abelian
  flat connections.
We assume that $A_{s,t}$ is in cylindrical form and has diagonal holonomy on the boundary.

Let $a_{s,t}$ denote the restriction of $A_{s,t}$ to the torus, and let $P^\pm_{s,t}$ be the
positive and negative eigenspans of $S_{a_{s,t}}.$ Since $$\ker(S_{a_{s,t}})
=H^{0+1+2}_{a_{s,t}}(T;\CC^2)=0$$ for $0 \leq s \leq 1$ and $0 \leq t\leq \ep,$ the Lagrangian
spaces $P^-_{s,t}$ vary continuously. Thus the odd signature operator $D_{A_{s,t}}$ acting on
sections over $Z$ with $P^-_{s,t}$ boundary conditions  is a continuous 2-parameter family of
self-adjoint operators. This 2-parameter family
 gives a homotopy from  the path
$D_{A_{0,t}}, \ 0 \leq t \leq \ep,$ to the composition of the three paths
\begin{enumerate}
\item[(i)] $D_{A_{s,0}}, \ 0 \leq s \leq 1.$
\item[(ii)] $D_{A_{1,t}}, \ 0 \leq t \leq \ep.$
\item[(iii)] $D_{A_{1-s,\ep}}, \ 0 \leq s \leq 1$
\end{enumerate}
and hence
$$SF(D_{A_t})=
SF(D_{A_{s,0}})_{s\in[0,1]}+SF(D_{A_{1,t}})_{t\in[0,\ep]}+SF(D_{A_{1-s,\ep}})_{s\in[0,1]}.
$$

The flat connections $A_{s,t}$ act nontrivially on $\CC^2,$ so it follows that
$H^0_{A_{s,t}}(Z;\CC^2)=0$ for all $s,t$. The path $A_{s,0}, 0\leq s \leq 1$ runs along the seam of
the cylinder and Propositions \ref{yum} and \ref{pad} show that $H^1_{A_{s,0}}(Z;\CC^2) = \CC^2$
for $0\leq s \leq 1$. By choosing $\ep$ sufficiently small, we can arrange that
$H^1_{A_{s,t}}(Z;\CC^2) =0$ for $0 \leq s \leq 1$ and $0<t\leq \ep.$ (For this deduction, notice
that $A_{s,t}$ has been twisted out of the $SU(2)\times \{1\}$ stratum for $t>0.$)

Since
  the kernel of $D_{A_{s,t}}$ with $P^-$ boundary conditions
 is isomorphic to the image
with $H^{0+1}_{A_{s,t}}(Z,\partial Z;\CC^2)\to H^{0+1}_{A_{s,t}}(Z;\CC^2)$ (see the paragraph
preceding the statement of Theorem \ref{KCSF}), and this restriction map is surjective by
Proposition \ref{yum}, it follows that along the first path $D_{A_{s,0}}$ the kernel is constant
(and 4-dimensional) and along the third path the kernel is trivial. Hence the spectral flow along
the first and third paths vanishes.
 Thus $$SF(D_{A_{0,t}}; \, Z; \, P^-_{0,t}; \, 0 \leq t \leq \ep)
= SF(D_{A_{1,t}}; \,  Z; \, P^-_{1,t}; \, 0 \leq t \leq \ep).$$

 We have now reduced the proof to computing $SF(D_{A_{1,t}}; \,  Z; \, P^-_{1,t}; \, 0 \leq t \leq \ep)$, along the
 path $A_{t,1}$ of abelian flat connections. We will show that the 4 zero modes bifurcate
 into two positive and two negative eigenvalues.  The idea of the argument is simple but the execution is a bit technical,
 so we outline the argument first.  We will embed the path $A_{t,1},t\in[0,\ep]$ in a 2-parameter family
  $B_{u,v}, (u,v)\in \RR^2$ so that $A_{t,1}$ corresponds to a short path starting at the origin moving along the positive
  $v$-axis. The operator $D_{B_{u,v}}$ with $P^-$ boundary conditions will be seen to have kernel of dimension 2
  along the two lines $v=u/3$ and $v=-u/3$ (and hence 4-dimensional kernel at the origin). The spectral flow along the
 $u$ axis through the origin (i.e. $SF(D_{B_{u,0}}, P^-, -\ep\leq u\leq \ep)$) equals 4 or $-4$. Thus in the 4 cone
 shaped regions complementary to the two lines, the two regions containing the positive and
  negative $v$ axis must correspond to two of the zero modes becoming positive and two
  becoming negative.

Since $A_{1,t}$ is an abelian flat connection on $Z$, it is completely
determined by its meridinal holonomy.
Suppose $\hol_\mu(A_{1,0}) = \Phi(u_0, 0)$
and let $B_{s,t}$ be a 2-parameter family of abelian flat connections
with  $B_{0,t}=A_{1,t}$ and $\hol_\mu(B_{s,t}) = \Phi(u_0+s,t)$.
Notice that each $B_{s,0}$ is an
$SU(2)\times \{1\}$ connection.

By \cite{APS},  the kernel of $D_{B_{s,t}}$ with $P^-$ boundary
conditions is isomorphic to  the image of the
relative cohomology in the absolute
$$\Image\left( H^1_{B_{s,t}}(Z,T;\CC^2)
\to H^1_{B_{s,t}}(Z;\CC^2) \right).$$  For $s$ and $t$ small, $H^*
_{B_{s,t}}(T;\CC^2) =0$, so the latter image is simply
$H^1
_{B_{s,t}}(Z;\CC^2)$, which is computed in  Proposition \ref{abelians}.
In the present context, this proposition implies that, for small $s$
and $t$,  the kernel of
$D_{B_{s,t}}$ with $P^-$ boundary conditions is
$$H^1 _{B_{s,t}} (Z; \CC^2) =
\begin{cases}
\CC^2 & \text{if $s=t=0$,}\\
\CC & \text{if $t=\pm \frac s 3 \neq 0$,} \\
0 & \text{otherwise.}
\end{cases} $$

For paths of $SU(2) \times \{1\}$ connections, the odd signature
operator respects the quaternionic structure on $\CC^2$, and
for this reason, the spectral flow
$$SF(B_{s,0}; \, Z; \, P^-; \,  -\ep \leq s \leq \ep) =\pm 4$$
(cf.~Theorem 6.12 in \cite{BHKK}).
We assume this spectral flow equals $+4$. The
argument in the other case is similar and is left to the reader.
Because  there are only four zero modes, all at $s=0$,
we see that the spectral flow along the first half of this path
$\{(s,0) \mid  -\ep \leq s \leq 0\}$ must also equal  $+4$ (by our
spectral flow conventions).

The straight line
$\{(s,0) \mid -\ep \leq s \leq \ep\}$
is
homotopic to  the semicircle
$\{(-\ep\cos\th, \ep \sin\th) \mid 0 \leq \th \leq \pi\}.$
The semicircle  passes through the two diagonal lines
through $(u_0,0)$ exactly once. Each time it crosses
a diagonal line $t=\pm \frac s 3$,
exactly one eigenvalue (of multiplicity two)
of $D_{B_{s,t}}$
crosses zero from negative to positive (since
the total spectral flow is $+ 4$).
Thus, the spectral flow  along the quarter circle
$\{(-\ep\cos\th,\ep\sin\th) \mid 0 \leq \th \leq \pi/2\}$
must equal $+2.$ Of course, the quarter circle is homotopic
to the composition of the two straight lines
$\{(s,0) \mid -\ep \leq s \leq 0\}$ and $\{(0,t) \mid 0 \leq t \leq
\ep\}.$
We already concluded that the spectral flow along
the first line  equals $+4$, hence the spectral flow along
the second must equal $-2.$ Thus
$$SF(B_{0,t};\, Z; \, P^-; \, 0 \leq t \leq \ep) = -2.$$
In other words, the behavior of the
four zero modes of $D_{A_{1,t}}$ as $t$ increases from $t=0$ is
that two go up, the other two go down.  This completes the proof.
\end{proof}

\section{Applications} \label{s6}
In this section, we present computations of
  the integer valued  $SU(3)$ Casson invariant
$\tau_{SU(3)}$ for Brieskorn spheres $\Si(p,q,r).$
As we know from Theorem \ref{sthm}, there are
exactly four types of path components, so our first
task is to explain how each type contributes
to $\tau_{SU(3)}$. This reduces the problem of computing
$\tau_{SU(3)}(\Si(p,q,r))$ to an enumeration problem,
which we then phrase and solve in terms of
counting lattice points in rational polytopes. From this, we deduce
that $\tau_{SU(3)}$ is a quadratic polynomial
in $n$ for $1/n$-Dehn surgery on a $(p,q)$ torus knot,and more
generally for the families $\Si_n=\Si(p,q,pqn+m)$
for $p,q,m>0$ fixed, relatively prime integers with  $m < pq.$

\subsection{The integer valued SU(3) Casson invariant}
In this subsection, we
  review how the different component types contribute to
the integer valued $SU(3)$ Casson invariant defined in \cite{BHK}. For
Brieskorn spheres $\Si$, let $h$ be a
small perturbation so  that $\cM_h$ is regular. The $SU(3)$
Casson invariant is given by
\begin{equation} \label{taudef}
\begin{split}
\tau_{SU(3)}(\Si) = & \sum_{[A] \in \cM^*_h} (-1)^{SF(\Th,A;\Si)} \\
& + \tfrac{1}{4}\sum_{[A] \in \cM^{\hbox{\tiny \sl red}}_h}
(-1)^{SF(\Th,A;\Si)}
\left( 2SF_{\CC^2}(\widehat{A},A;\Si)+  \dim
H^1_{\widehat{A}}(\Si;\CC^2) \right).
\end{split}
\end{equation}
In this formula, $\dim H^1$ refers to the real dimension
and $\widehat{A}$ is a
reducible flat connection chosen close to a fixed representative
  $A$ of the gauge orbit $[A] \in \cM^\red_h$.

\begin{remark}
The general definition of
$\tau_{SU(3)}$ in \cite{BHK} is more complicated; it involves choosing two basepoints
$[\widehat{A}_+]$ and $[\widehat{A}_-]$ for each path component
of $\cM^\red$. But  for Brieskorn spheres $\Si$,
$\cM^\red(\Si)$ is discrete, so we take
$\widehat{A}_+=\widehat{A}_-$, and the definition of \cite{BHK}
reduces to (\ref{taudef}).
\end{remark}

Using standard results from Morse theory,
one can show that each Type Ia or Type IIa path component
contributes  $\pm 1$ times
its Euler characteristic to $\tau_{SU(3)}.$
It is a general fact (cf.~the proof of Lemma 7 in \cite{BHK})
that isolated reducible orbits
with vanishing normal cohomology do not contribute
to $\tau_{SU(3)}$. Thus the Type Ib components do not contribute
to $\tau_{SU(3)}.$ For the Type IIb components, the pointed
2-spheres, we
apply the twisting perturbation to resolve the
singularity and then use the spectral flow computations
of Section \ref{s5} to calculate the contribution.
These results are summarized in the following theorem.

\begin{thm} \label{cct}
Suppose $\Si$ is a Brieskorn sphere.
The contribution of a given path component of $R(\Si,SU(3))$
to  the integer valued $SU(3)$ Casson invariant
$\tau_{SU(3)}$ depends only on the component type and
is as follows.
\begin{enumerate}
\item[(i)] Type Ia components, which are isolated  points of
conjugacy class of irreducible $SU(3)$ representations,
contribute $+1$ to $\tau_{SU(3)}(\Si)$.
\item[(ii)] Type IIa components, which are a smooth
2-spheres of
conjugacy classes of irreducible $SU(3)$ representations,
contribute $+2$ to $\tau_{SU(3)}(\Si)$.
\item[(iii)] Type Ib components, which are isolated points of
conjugacy classes of reducible $SU(3)$ representations,  do
not contribute  to $\tau_{SU(3)}(\Si)$.
\item[(iv)] Type IIb components, which are pointed
2-sphere containing one
conjugacy class of reducible $SU(3)$ representations,
contribute $+2$ to $\tau_{SU(3)}(\Si)$.
\end{enumerate}
\end{thm}

\begin{proof}
This theorem uses Proposition 5.1  from \cite{B2}, which states
that for any irreducible flat $SU(3)$ connection $A$ on $\Si,$
the adjoint $su(3)$ spectral flow $SF(\Th,A)$ is even.
Given a
nondegenerate component $\cC \subset R^*(\Si,SU(3))$
and $[A] \in \cC,$ Proposition 8 of \cite{BH2}
states that $\cC$ contributes $(-1)^{SF(\Th,A)} \chi(\cC)$ to
$\la_{SU(3)}$.
But the only difference between
the invariants $\tau_{SU(3)}$ and $\la_{SU(3)}$ is in
their correction terms. In other words,
on the level of the irreducible stratum, these two invariants coincide.
Thus, since components of Types I and II are nondegenerate,
we conclude that
components of Type Ia
contribute $+1$
and components of Type IIa contribute $+2$
to $\tau_{SU(3)}(\Si)$.

Next, consider a component $\cC$  of Type Ib.
Thus $\cC=\{[A_0]\}$ for an isolated reducible orbit $[A_0] \in
\cM^\red$.
Proposition \ref{coho} implies $H^1_{A_0}(\Si;\CC^2)=0$.
Given a generic path $h_t$ of small perturbations, the
path $A_t$ of nearby reducible $h_t$-perturbed flat connections
also have $H^1_{A_t,h_t}(\Si;\CC^2) =0.$ As a result,
$SF_{\CC^2}(A_t,h_t;\Si) =0$ and we conclude
that
components of Type Ib do not
contribute to $\tau_{SU(3)}$.

Finally, consider a component $\cC$ of Type IIb.
So $\cC$ is a pointed 2-sphere and has two strata:
  $\cC = \cC^* \cup \cC^\red$.
Let $tf$ be
the path of twisting perturbations  on $\Si$ as in Section \ref{s4}.
Denote by $\cC_t \subset \cM_{ t f}$
the part of the $(t f)$-perturbed flat moduli space of
$\Si$ near $\cC$. As we have shown, for $t$ small, $\cC_t$ is a disjoint union of two
components
$$\cC_t = \cC^*_t \cup \cC^\red_t.$$

Choose $\ep>0$ as in Theorem \ref{SFST} and
suppose $[B_t] \in \cC^*_t$ is a path
of gauge orbits of irreducible $(tf)$-perturbed flat
connections on $\Si$ for $0 \leq t \leq \ep$.
Then $H^1_{B_t,tf}(\Si;su(3)) = \RR^2$
for $0 \leq t \leq \ep,$ and hence
$$SF(B_t,tf; \Si; 0 \leq t \leq \ep) =0.$$
Since $SF(\Th,B_0;\Si)$ is even, another application of
Proposition 8 of \cite{BH2}, together with the
fact that $\cC^*_\ep$ is a nondegenerate 2-sphere,
shows that $\cC^*_\ep$  contributes $+2$ to
$\tau_{SU(3)}(\Si).$

Now suppose $[A_t] \in \cC^\red_t$
is a path of gauge orbits of reducible $(tf)$-perturbed flat
connections on $\Si$.
Corollary \ref{noH1forsmallt}
implies that $\{ [A_t]\}$ is isolated for $0 < t \leq \ep,$
and Theorems \ref{SFST} and \ref{KCSF} imply that
$SF_{\CC^2}(A_0,A_\ep;\Si) = -2.$
In addition, Proposition \ref{coho} tells us
that $H^1_{A_0}(\Si;\CC^2) = \CC^2.$
Thus
$$2 SF_{\CC^2}(A_0,A_\ep;\Si) + \dim H^1_{A_0}(\Si;\CC^2) = -4+4=0,$$
and
the contribution of
$\cC^\red_\ep$ to $\tau_{SU(3)}(\Si)$ is $0$.
Consequently, each component of Type IIb contributes
$+2$ to $\tau_{SU(3)}(\Si)$,
and this completes the proof.
\end{proof}

\subsection{SU(3) fusion rules} \label{sfuse}
The set of $SU(3)$ matrices modulo conjugation is parameterized by
the 2-simplex
\begin{equation} \label{tri}
\De := \{(a_1,a_2, a_3) \in \RR^3 \mid a_1 \leq a_2 \leq a_3
\leq a_1+1
\hbox{ and }a_1+a_2+a_3 =0\}.\end{equation}
Suppose $\Ga$ is a discrete group
and $\al \in R(\Ga,SU(3))$. Define
the map $\la_\al\colon \Ga \to \De$ by
sending $\ga \in \Ga$ to the
unique $(a_1,a_2, a_3) \in \De$ such
that  $\al(\ga)$ has eigenvalues $e^{2 \pi i a_1}, e^{2 \pi i a_2},e^{2
\pi i a_3}.$

The fundamental group of a thrice-punctured 2-sphere has the
presentation
$G = \langle x,y,z \mid xyz=1\rangle,$
where $x,y,z$ are represented by loops around the three punctures.
(Of course $G$ is  a free group on 2 generators.)
Given any representation
$\al\colon G \to SU(3)$,
the assignment $\al \mapsto \left(\la_\al(x),
\la_\al(y),\la_\al(z)\right)$ defines  a map
$$\Psi\colon R(G,SU(3)) \lto \De
\times\De \times\De.$$
The following theorem, due to Hayashi (see Theorems 3.3 and 3.4 of
\cite{Ha}), describes the image of this map as a convex 6-dimensional
polytope $\cP$ in
$\De \times\De \times\De.$

Given $\ba,\bb,\bc \in \De$, let $\cM_{\ba \bb \bc}$ be the
moduli space of flat connections
on a thrice-punctured 2-sphere with monodromies around the three punctures
specified by $\ba,\bb,\bc.$  Clearly $\cM_{\ba\bb\bc}$ can be
identified with
  the fiber of the map $\Psi$ over $(\ba,\bb,\bc).$

\begin{thm} \label{fused}
  The moduli space  $\cM_{\ba\bb\bc}$
  is nonempty if and only if $\ba=(a_1,a_2,a_3), \bb=(b_1,b_2,b_3)$  and
  $\bc=(c_1,c_2,c_3)$   satisfy the 18 inequalities:
\begin{equation} \label{fusion}
\vspace*{0.1in} \left\{ \begin{split}
\hspace*{0.25in}
&a_1+b_2+c_2 \leq 0 \hspace*{0.5in} a_1+b_3+c_3 \geq 0\hspace*{0.5in}
a_2+b_3+c_3 \leq 1  \\
&a_2+b_1+c_2 \leq 0 \hspace*{0.5in} a_3+b_1+c_3 \geq 0\hspace*{0.5in}
a_3+b_2+c_3 \leq 1 \\
&a_2+b_2+c_1 \leq 0 \hspace*{0.5in} a_3+b_3+c_1 \geq 0\hspace*{0.5in}
a_3+b_3+c_2 \leq 1 \\
\\
&a_2+b_2+c_3 \geq 0 \hspace*{0.5in} a_1+b_1+c_3 \leq 0 \hspace*{0.5in}
a_1+b_1+c_2 \geq -1 \\
&a_2+b_3+c_2 \geq 0 \hspace*{0.5in} a_1+b_3+c_1 \leq 0 \hspace*{0.5in}
a_1+b_2+c_1 \geq -1\\
&a_3+b_2+c_2 \geq 0 \hspace*{0.5in} a_3+b_1+c_1 \leq 0 \hspace*{0.5in}
a_2+b_1+c_1 \geq -1.\\
\end{split} \right.
\vspace*{0.1in} \end{equation}
Let
  $\cP=\{(\ba,\bb,\bc) \mid \text{all 18 of the inequalities
(\ref{fusion}) are satisfied}\}.$
Then
$\cP= \im(\Psi)$ is convex and 6-dimensional.
Moreover, $\cM_{\ba\bb\bc}$ is homeomorphic to
  a 2-sphere if  $(\ba,\bb,\bc)$ lies in the interior of $\cP$and a
point
  if $(\ba,\bb,\bc)$ lies
  on the boundary of $\cP$.
\end{thm}

These equations can be used to describe the irreducible stratum
$R^*(Z,SU(3))$ of
the  representation variety of $\pi_1 Z$ as follows.
Fix  $A,B \in SU(3)$ and $\ell \in \{0,1,2\}$ as in Theorem \ref{Zsthm}
and let $\ba, \bb \in \De$ be the conjugacy classes of
$A,B$, respectively.
Recall the presentation (\ref{pres-Z}) for $\pi_1 Z$ and
denote by $\cC_{\ba \bb}^\ell \subset R(Z,SU(3))$ the subset
consisting of conjugacy classes of representations $\al\colon \pi_1 Z \to
SU(3)$
such that $\la_\al(x) = \ba, \, \la_\al(y) = \bb,$ and $\al(h)= e^{2
\pi i \ell/3}I.$
(This set was denoted $\cC_{AB}^\ell$ in Theorem \ref{Zsthm}.)

The assignment $\al \mapsto \la_\al((xy)^{-1})$ defines a map
$$\psi_{\ba\bb}^\ell \colon \cC_{\ba\bb}^\ell \lto \De.$$
Let $Q^\ell_{\ba\bb}=im(\psi_{\ba\bb}^\ell)$ be the image of this
map, so $Q^\ell_{\ba\bb}$ is the intersection
of $\cP$ with the 2-dimensional slice obtained by fixing $\ba$ and
$\bb.$
Solving equations (\ref{fusion}) for $c_1,c_2,c_3$, we see that
$$Q^\ell_{\ba\bb} \subset \De = \{(c_1,c_2,c_3) \in \RR^3 \mid
c_1 \leq c_2 \leq c_3 \leq c_1+1
\hbox{ and } c_1+c_2+c_3 =0\}$$
consists of triples $(c_1,c_2,c_3)$ satisfying the six inequalities:
\begin{eqnarray*}
\Xmin \; \leq  &c_1 &\leq  \; \Xmax, \\
\Ymin \; \leq  &c_2 &\leq  \; \Ymax, \\
\Zmin \; \leq  &c_3 &\leq  \; \Zmax,
\end{eqnarray*}
where \begin{eqnarray*}
\Xmin &= &\max\{ -1-a_1-b_2, -1-a_2-b_1, -a_3-b_3\},\\
\Xmax &= &\min \{ -a_1 - b_3, -a_3-b_1, -a_2-b_2\}, \\
\Ymin &= &\max\{ -1-a_1-b_1, -a_2-b_3, -a_3-b_2\}, \\
\Ymax &= &\min \{ -a_1 - b_2, -a_2-b_1, 1-a_3-b_3\}, \\
\Zmin &= &\max\{ -a_1-b_3, -a_3-b_1, -a_2-b_2\}, \\
\Zmax &= &\min \{ -a_1 - b_1, 1-a_2-b_3, 1-a_3-b_2\}. \\
\end{eqnarray*}
Using these equations, one can determine
that $Q^\ell_{\ba\bb}$ is either a hexagon
or a nonagon, depending on whether $\cC^\ell_{\ba\bb}$ is a Type I
or II component, respectively (recall the definition of Type I and II
in Theorem \ref{Zsthm}). With a little more work,
one sees that  the vertices of  $Q^\ell_{\ba\bb}$ are given
by
\begin{equation} \label{TypeIvert}
   \begin{split}
V_1 &= (\Xmax, -\Xmax-\Zmin, \Zmin), \quad
V_2 = (-\Ymax-\Zmin, \Ymax, \Zmin), \quad
V_3 = (\Xmin, \Ymax, -\Xmin-\Ymax), \\
V_4 &= (\Xmin, -\Xmin-\Zmax, \Zmax),  \quad
V_5 = (-\Ymin-\Zmax, \Ymin, \Zmax),   \quad
V_6 = (\Xmax, \Ymin, -\Xmax-\Ymin), \\
\end{split}
\end{equation}
in the hexagonal case (i.e. when $\cC^\ell_{\ba\bb}$ is Type I), and by
  \begin{equation} \label{TypeIIvert}
   \begin{split}
V_1 &= (\Xmax, -\Xmax-\Zmin, \Zmin),\;
V_2 = (-2\Zmin, \Zmin, \Zmin),   \qquad \qquad \,
V_3 = (-2\Ymax, \Ymax, \Ymax),  \\
V_4 &= (\Xmin, \Ymax, -\Xmin-\Ymax), \;\;
V_5 = (\Xmin, -1-2\Xmin, 1+\Xmin),  \;
V_6 = (\Zmax-1, 1-2\Zmax, \Zmax),   \\
V_7 &= (-\Ymin-\Zmax, \Ymin, \Zmax), \;\;\;\,
V_8 = (\Ymin, \Ymin, -2\Ymin),  \qquad \qquad \; \,
V_9 = (\Xmax, \Xmax, -2\Xmax),
\end{split}
\end{equation}
in the nonagonal case (i.e. when $\cC^\ell_{\ba\bb}$ is Type II).

\subsection{Lattice points in rational polytopes} \label{lattice}
In this subsection, we use Ehrhart's theorems on enumerating
lattice points in rational polytopes to establish two
results. The first, Theorem
\ref{torusconj}, is essential for the
computations in Subsection \ref{conclud}.
It shows that the integer valued $SU(3)$
Casson invariant on homology 3-spheres
obtained by $1/n$ surgery on a torus knot (or torus-like knot)
is a quadratic polynomial in the surgery coefficient $n$.
The second result, Proposition \ref{Numb},
enumerates Type I and II components
in the $SU(3)$ representation variety
of knot complements $Z$ obtained
by removing one of the singular fibers of $\Si(p,q,r).$

To begin, suppose $\Si=\Si(p,q,r)$ is a  Brieskorn sphere and $Z$ is
the complement of a regular neighborhood of its singular $r$-fiber.
Recall the presentations (\ref{pres}) and (\ref{pres-Z})
for the fundamental groups $\pi_1 \Si$ and $\pi_1 Z.$
Restriction from $\Si$ to $Z$ defines a natural inclusion map
$R(\Si,SU(3)) \hookrightarrow R(Z,SU(3)),$
under which
\begin{equation} \label{identt}
R(\Si,SU(3)) = \{ \al\colon\pi_1 Z \to SU(3) \mid
\al((xy)^r h^c) = I\} / {\rm conj} \subset R(Z,SU(3)).
\end{equation}
Any irreducible representation  $\al\colon \pi_1 Z \to SU(3)$
must send $h$ to a central element, thus
$\al(h) = e^{2 \pi i \ell/3}I$ for some $\ell \in \{0,1,2\}.$
Hence  $\al(x)$ and $\al(y)$
are $p$-th and $q$-th roots of the central element $\al(h)^a = e^{2 \pi
i \ell a/3}I,$
and the results in Subsection \ref{dot} imply that
  $R^*(Z,SU(3))$
is a union of components $\cC_{\ba \bb}^\ell$ over all
$\ba,\bb \in \De$ and $\ell \in \{0,1,2\},$ of the
form
\begin{equation}\label{abform}
\ba=\left(\tfrac{i_1}{3p}, \tfrac{i_2}{3p},
\tfrac{-i_1-i_2}{3p}\right), \quad
\bb=\left(\tfrac{j_1}{3q}, \tfrac{j_2}{3q}, \tfrac{-j_1-j_2}{3q}\right),
\end{equation}
where $i_1,i_2,j_1,j_2$ are integers satisfying $i_1 \equiv i_2 \equiv
j_1 \equiv j_2 \equiv a\ell \pmod{3}.$

A conjugacy class $[\al] \in \cC^\ell_{\ba\bb}$
with representative $\al \colon \pi_1 Z \to SU(3)$
extends to a representation
of $\pi_1 \Si = \pi_1 Z/\langle (xy)^rh^c\rangle$
if and only if $\al((xy)^r h^{c})=I.$
Setting $\bc=\la_\al((xy)^{-1}) \in Q^\ell_{\ba\bb},$
  we see that $\al$ extends if and only if
\begin{equation}\label{cform}
\bc = \left(\tfrac{k_1}{3r}, \tfrac{k_2}{3r},
\tfrac{-k_1-k_2}{3r}\right)
\end{equation}
for integers $k_1,k_2$ such that $k_1 \equiv k_2 \equiv c \ell \pmod{3}.$

In this way, we reduce the problem  of computing
$\tau_{SU(3)}(\Si)$ to one of counting lattice points of the form
(\ref{cform})
in the regions $Q^\ell_{\ba\bb}$, for all $\ba,\bb,\ell$
satisfying (\ref{abform}).
Of course, some lattice points contribute $+1$
and others contribute $+2$, depending on the topology of the
fiber of $\psi_{\ba\bb}$ (cf. Theorem \ref{fiber}).
This is a routine matter, as the topology of the fibers is constant
within the interior of $Q^\ell_{\ba\bb}.$

The same approach can be used to perform
computations for the entire family of Brieskorn spheres
$$\Si_n:= \Si(p,q,pqn+m), \ n\geq 0,$$ where
$p,q,m$ are fixed, pairwise relatively prime positive integers
with $ m < pq$.
We have described  $R(\Si_n,SU(3))$ as a disjoint union of
points and 2-spheres. Under the identification (\ref{identt}), each
point and 2-sphere corresponds to
a lattice point in  one of the regions $Q^\ell_{\ba\bb}$.
Observe that the regions $Q^\ell_{\ba\bb}$ are
   themselves independent of $n$;
the dependence
on $n$  is entirely through the denominators
of the lattice points via equation (\ref{cform}) and
$r=pqn+m$.

\begin{thm} \label{torusconj}
Suppose $p,q,m>0$ are pairwise relatively prime with $m < pq.$ Set
$\Si_n = \Si(p,q,pqn+m).$ Then
$\tau_{SU(3)}(\Si_n)$ is a quadratic polynomial in $n$ of the form
$$\tau_{SU(3)}(\Si_n) = An^2+Bn+C.$$
Obviously $C= \tau_{SU(3)}(\Si(p,q,m))$ and vanishes for $m=\pm 1$.
\end{thm}

Our proof uses Ehrhart's results on counting lattice points
in rational polytopes \cite{ehrhart},
so we begin by introducing notation and defer the proof to
the end of this subsection.

A {\it lattice polytope} $\cP$ in $\RR^N$ is a
convex polytope whose vertices lie on the standard
integer lattice $\La = \ZZ^N$, and a {\it rational polytope}
$\cQ$ in $\RR^N$ is one whose vertices have rational coordinates.
Equivalently, $\cQ$ is rational if
the dilated region $d\cQ = \{ dx \mid x \in \cQ \}$
is a lattice polytope
for some positive integer $d$.
For example, the 2-simplex
$\De$ of equation (\ref{tri})
is a rational polytope which, when dilated by $d=3$, is a lattice polytope.

We are interested in counting lattice points in
integral dilations $n \cP$ of such polytopes.
Denote by $f_\La(\cP,n)= \#\left(n\cP \cap \La \right)$,
the number of lattice points in $n  \cP.$
Ehrhart showed that if $\cP$ is a lattice polytope, then
$f_\La(\cP,n)$ is a polynomial in $n$ of degree $\dim \cP$.
Ehrhart also proved that if $\cQ$ is a rational polytope such that
$d \cP$ is a lattice polytope, then
$f_\La(\cQ,n)$ is a {\em quasi-polynomial}
of degree $\dim \cQ$ and periodicity $d$, where
(see \cite{ehrhart} or p.235 of \cite{St}).
Recall that a quasi-polynomial $f(n)$ of degree $j$ and periodicity $d$
is a function of the form
$$f(n) =  \sum_{i=0}^{j} a_i(n) n^i$$
whose coefficient functions $a_i(n)$ are periodic in $n$ of period $d$.

Fix $p,q,m$ and set $\Sigma_n:=\Sigma(p,q,pqn+m)$
as in the theorem. Choose
integers $a_n,c_n$ satisfying
\begin{equation}\label{tats}
  a_n(pqn+m)(p+q) + c_n pq =1
\end{equation}
  as in
Proposition \ref{forceit}.
Denote by
$Z_n$  the complement of a regular neighborhood of the
$(pqn+m)$-fiber in $\Si_n=\Si(p,q,pqn+m)$. The fundamental group
$\pi_1 Z_n$ has presentation
$\langle x,y,h \mid x^p = y^q = h^{a_n}, \ h \text{ central}\rangle.$
We will see that the Type I and II
components $\cC_{\ba\bb}^\ell$ of $R(Z_n,SU(3))$
are independent of $n$.
(Here, as established in Theorem \ref{fiber},
components of Types I and II have real dimension two and four,
respectively.)

We will identify components of $R^*(\Sigma_n,SU(3))$ with the union over
all $\ba,\bb$ of certain lattice
points in
$Q^\ell_{\ba \bb} \subset \RR^3,$
and a key point
is that these regions depend only on $\ba,\bb$ and not on $n$.

\begin{lem} \label{chwzac}
  The numbers $a_n, c_n$ can be chosen
so their values modulo three
are independent of $n$. Moreover:
\begin{enumerate}
\item[(i)] If both $p$ and $q$ are relatively prime to 3,
then we can choose $a_n, c_n$ so that
$a_n \equiv 0 \pmod{3}$ and
  $c_n \equiv pq \not\equiv 0  \pmod{3}.$
\item[(ii)] If either $p$ or $q$ is a multiple of 3,
then we can choose $a_n,c_n$ so that
$a_n \equiv  (p+q)m \not\equiv 0 \pmod{3}$ and $c_n \equiv -m \not\equiv
0 \pmod{3}$.
\end{enumerate}
\end{lem}

\begin{proof}
We start with $a_n,c_n$ satisfying (\ref{tats}) and
use the substitutions $a_n' = a_n + pqk$
and $c'_n = c_n - k(p+q)(pqn+m)$.
For example, in case (i), we can choose $k$
so that $a'_n$ is a multiple of $3$ since $pq$ is relatively
prime to 3. Reducing
equation (\ref{tats})
modulo 3 then implies that $c'_n \equiv pq \pmod{3}.$
In case (ii), the mod 3 reduction of equation
(\ref{tats}) gives that $a_n \equiv (p+q)m$ before (and after)
making any substitutions. Now since $(p+q)m$ is relatively prime to 3,
so is $(p+q)(pqn+m)$, and it follows that
we can substitute so that
$c'_n\equiv -m \pmod{3}$.
\end{proof}

\begin{remark} \label{rmkit}
In case (i), a consequence of Lemma \ref{chwzac} is that
$\ba$ has the form $(\frac{i_1}{p},\frac{i_2}{p},\frac{-i_1-i_2}{p})$
and $\bb$ has the form
$(\frac{j_1}{q},\frac{j_2}{q},\frac{-j_1-j_2}{q})$
when $p,q$ are both relatively prime to 3 (cf. equation (\ref{abform})).
In this case, the three components
$\cC^0_{\ba,\bb},\cC^1_{\ba,\bb},\cC^2_{\ba,\bb}$
have the same values for $\ba,\bb$.

In case (ii), we see that $\ell$ is completely determined by $\ba$ (or
$\bb$)
since $a_n \not \equiv 0 \pmod{3}$ when $p$ or $q$ is a multiple of 3.
In this case, different values of $\ell$
require different values of $\ba,\bb$.
\end{remark}

The next result gives an enumeration of the
number of Type I and Type II components in $R(Z_n,SU(3))$.
\begin{prop} \label{Numb}
Suppose $Z_n$ is the complement of the $(pqn+m)$-singular fiber
in $\Sigma(p,q,pqn+m)$. Then there are
\begin{eqnarray*}
N_{I} &=& \frac{(p-1)(q-1)(p+q-4)}{2} \qquad \text{and} \\
N_{II} &=&  \frac{(p-1)(p-2)(q-1)(q-2)}{12}
\end{eqnarray*}
components of Type I and Type II
  in $R(Z_n,SU(3))$, respectively.
\end{prop}

The next lemma is the key to proving this proposition.
\begin{lem} \label{nrpq}
Suppose $p \in \ZZ$ is a positive integer and $\ell \in \{0,1,2\}.$
Let $f_\ell(p)$  denote the number of conjugacy
classes of $p$-th
roots of $e^{2 \pi i \ell/3} I$ in $SU(3)$ with three distinct
eigenvalues, and let $g_{\ell}(p)$ denote the number of
conjugacy classes of $p$-th
roots of $e^{2 \pi i \ell/3} I$ in $SU(3)$ with two distinct
eigenvalues.
Then we have:
$$f_{\ell}(p)=\begin{cases}
\frac{1}{6}(p^2-3p+2) & \text{if $p$ is relatively prime to 3,} \\
\frac{1}{6}(p^2-3p+6) & \text{if $p$ is multiple of 3 and $\ell=0$,}\\
\frac{1}{6}(p^2-3p) & \text{if $p$ is multiple of 3 and $\ell=1,2$.}
\end{cases} $$

$$g_{\ell}(p)=\begin{cases}
p-1 & \text{if $p$ is relatively prime to 3,} \\
p-3 & \text{if $p$ is multiple of 3 and $\ell=0$,}\\
p & \text{if $p$ is multiple of 3 and $\ell=1,2$.}
\end{cases} $$
\end{lem}
Observe that
$\sum_{\ell=0}^2 f_{\ell}(p) = \frac{1}{2}(p-1)(p-2)$ and
$\sum_{\ell=0}^2 g_{\ell}(p) = 3p-3$ hold for all
$p$.

\begin{proof}
We begin by proving the stated formulas for $f_\ell(p)$ and
$g_\ell(p)$ under the assumption that $p$ is relatively prime to 3.

Consider the analogous problems for $U(3)$.
Set $\zeta=e^{2 \pi i/p}$ and notice that
a $p$-th root of unity in $U(3)$ has eigenvalues in the set
$\{ 1, \zeta, \zeta^2, \ldots, \zeta^{p-1} \}$.
Conjugacy classes in $U(3)$ are uniquely determined by their eigenvalues,
and it follows that there are  ${p \choose 3}$
conjugacy classes of $p$-th roots of unity in $U(3)$
with three distinct eigenvalues and that there are
$p(p-1)$ conjugacy classes of $p$-th roots of unity in $U(3)$
with two distinct eigenvalues

Multiplication by $\zeta$ defines a $\ZZ_p$ action
on these conjugacy classes. Using that
$\det \left( \zeta A\right) = \zeta^3 \det A,$ we see that
with respect to the map $\det\colon U(3) \to U(1)$,
the induced $\ZZ_p$ action downstairs on $U(1)$ has weight three.
If $(3,p)=1,$ the action
is effective on the image
$\det(\{ A \mid A^p = I\}) = \{ 1, \zeta, \zeta^2, \ldots, \zeta^{p-1} \}.$

Thus, if  $(3,p)=1$, the number of conjugacy classes
of $p$-th roots of unity in any fiber $\det^{-1}(\zeta^k)$
is  independent of $k$. Taking $k=0,$
it follows that $f_0(p)=\frac{1}{p} {p \choose 3}= (p-1)(p-2)/6$
and $g_0(p) = p-1$ if $(3,p)=1.$
Now multiplication by $e^{2 \pi i/3}$ shows that
$f_\ell(p) = f_{\ell+ p}(p)$ and $g_\ell(p) = g_{\ell+ p}(p)$.
Thus, if  $p$ is relatively prime to 3, it follows that
$f_\ell(p)$ and $g_\ell(p)$ are independent of $\ell \in \{0,1,2\}$
and are as stated in the lemma.

Now suppose $p$ is a multiple of 3 and notice that
 the $\ZZ_p$ action is no longer
effective on the image
$\det(\{ A \mid A^p = I\}) = \{ 1, \zeta, \zeta^2, \ldots, \zeta^{p-1} \}.$ Since the action
has weight three, there are precisely three orbits of the $\ZZ_p$
action, one orbit for each residue class of $k \pmod{3},$
where $\det A = \zeta^k.$

 \begin{claim} \label{pqclaim} If $p$ is a multiple of 3, then \\
 (i) $f_0(p)= \frac{p^2-3p+6}{6}$ and \\
(ii) $g_0(p) = p-3.$
\end{claim}

Establishing the claim proves the lemma, as we now explain.
Taking matrix inverses shows
that $f_1(p) = f_2(p)$ and $g_1(p) = g_2(p).$
As argued before,
the total number of $p$-th roots of unity in $U(3)$
with three distinct eigenvalues is ${p \choose 3}$,
and total number of $p$-th roots of unity in $U(3)$
with two distinct eigenvalues is $p(p-1).$
This gives the formulas
$${\textstyle \sum_{\ell=0}^2 f_\ell(p) = \tfrac{3}{p}{p \choose 3} = \frac{(p-1)(p-2)}{2}\quad \text{ and }
\quad \sum_{\ell=0}^2 g_\ell(p) = \frac{3}{p}(p^2-p)=3(p-1),}$$
which can then be used to solve for $f_1(p), g_1(p)$ in terms of $f_0(p),
g_0(p).$

Part (ii) of Claim \ref{pqclaim} can be proved directly. Every
conjugacy class is uniquely determined by its set of eigenvalues,
which for a $p$-th root of unity in $SU(3)$ with a double eigenvalue
is a set of the form $\{\zeta^k, \zeta^k, \zeta^{-2k}\}$
for $1 < k \leq p-1$ with
$k\neq m, 2m$. (Note: the conditions on $k$
ensure that $ \zeta^k \neq  \zeta^{-2k}.$) There are clearly $p-3$ such sets.

The direct argument for part (i) of Claim \ref{pqclaim}
is somewhat tedious,
so we argue indirectly as follows. Note that the total
number of conjugacy classes of
$p$-th roots of unity in $SU(3)$ includes
the three central matrices $I, e^{2 \pi i/3} I, e^{4 \pi i/3} I,$ as well as
the $p-3$ conjugacy classes
with two eigenvalues listed above.
The set
$$\{ (\zeta^a, \zeta^b,\zeta^{-a-b}) \mid 1 \leq a,b \leq p\}$$
of order $p^2$  lists all
possible eigenvalues of $p$-th roots of unity
as {\it ordered} sets. Subtracting 3 for the central roots
and $3(p-3)$ for the $p$-th roots of unity with two
distinct eigenvalues (each one being listed 3 times as ordered sets),
and dividing by the order of the symmetric group
$S_3$, we get that
$$f_0(p) = \tfrac{1}{6}(p^2 - 3(p-3) - 3) = \tfrac{p^2-3p+6}{6}$$
as claimed. This proves the claim and completes the proof of
the lemma.
\end{proof}

\noindent{\it Proof of Proposition \ref{Numb}. \ }
We consider the following two cases:

\smallskip\noindent
{\sc Case 1:} Both $p$ and $q$ are relatively prime to 3. \\{\sc Case
2:} One of $p$ or $q$ is a multiple of 3.
\smallskip\noindent

Assume 1 holds and choose $a_n\equiv 0 \pmod{3}$ as
in Lemma \ref{chwzac}.
Given an irreducible representation
$\al\colon\pi_1 Z_n \to SU(3),$ we have
$\al(h) =e^{2 \pi i \ell/3}$ for some $\ell \in \{0,1,2\}.$
Then for each $\ba, \bb \in \De$ with
$p\cdot \ba, q\cdot \bb \in \La=\ZZ^3$,
there are three isomorphic copies
of $\cC_{\ba\bb}^\ell$, one for each possible value
of $\ell$. Thus $N_I = 3(f_0(p)g_0(q)+g_0(p)f_0(q))$
and $N_{II }= 3f_0(p)f_0(q)$, and the formulas
for Lemma \ref{nrpq}
complete the argument in this case.

Now assume 2 holds, and note that
$a_n \not\equiv 0 \pmod{3}$ by Lemma \ref{chwzac}.
Without loss of generality,
we can assume that $p$ is a multiple of 3 and
that $q$ is relatively prime to 3.
The number of Type Ia components  is given by
summing over the possible values for $\ell \in \{0,1,2\}$,
and similarly for the number of Type IIa components.
Lemma \ref{chwzac} implies that
$f_{\ell}(q)=\frac{1}{6}(q-1)(q-2)$ and
$g_{\ell}(q) = q-1$
independent of $\ell$. It also gives that
$$\sum_{\ell=0}^2 f_{\ell}(p) = \tfrac{1}{2}(p-1)(p-2)\quad \text{ and
} \quad
\sum_{\ell=0}^2  g_{\ell}(p) = 3p-3.$$ Using these
formulas,
one computes that
\begin{eqnarray*}
N_I &=& \sum_{\ell=0}^2
f_{\ell}(p) g_{\ell}(q) + g_{\ell}(p) f_{\ell}(q)=
\tfrac{1}{2}(p-1)(q-1)(p+q-4),\\
N_{II }&=& \sum_{\ell=0}^2 f_{\ell}(p) f_{\ell}(q)=
\tfrac{1}{12}(p-1)(p-2)(q-1)(q-2),\\
\end{eqnarray*}
  completing the proof of the proposition.\hfill $\Box$

\bigskip
\noindent{\it Proof of Theorem \ref{torusconj}. \ }
It is enough to show that the contribution of
each component $\cC^\ell_{\ba\bb}$ in $R(Z_n,SU(3))$
to $\tau_{SU(3)} (\Si_n)$ is quadratic in $n$.
As with the proposition, there are two cases.

\medskip \noindent
{\sc Case 1:} Both $p$ and $q$ are  relatively prime to 3. \\
{\sc Case 2:} One of $p$ or $q$ is a multiple of 3.

\medskip
In order to apply Ehrhart's theorem, we consider
translations of the standard lattice and (in Case 2) of
the rational polytopes $Q^\ell_{\ba\bb}$.

Assume 1 holds and choose
$a_n \equiv 0 \pmod{3}$ and $c_n \equiv pq \pmod{3}$
as in Lemma \ref{chwzac}.
As noted in Remark \ref{rmkit}, the sets
$Q^\ell_{\ba\bb}$ are identical for the
different $\ell \in \{0,1,2\}$
corresponding to the different choices for
$\al(h) = e^{2 \pi i \ell/3}I$.
It follows from equations (\ref{TypeIvert}) and (\ref{TypeIIvert})
that $Q^\ell_{\ba\bb}$ is a rational polytope
whose dilation by $d=pq$ is a lattice polytope.

When $\ell=0,$ the component $\cC^0_{\ba\bb}$ contributes
$$f_\La(Q^0_{\ba\bb}, pqn+m)
=\# \left((pqn+m)Q^0_{\ba\bb} \cap \La\right)$$
to $\tau_{SU(3)} (\Si_n)$, where
$\La=\ZZ^3$ is the standard integer lattice in $\RR^3.$
By \cite{ehrhart},
$f_\La(Q^0_{\ba\bb},k)$ is a quasi-polynomial of periodicity $d=pq,$
and we see that $f_\La(Q^0_{\ba\bb},pqn+m)$ is polynomial in $n$
simply because the residue class of $pqn+m$ modulo $d=pq$
is constant.

This same idea should work for $\ell =1,2,$ but there are
difficulties adapting the argument to these cases individually.
Instead, we combine the three cases $\ell=0,1,2$ by superimposing
the three sets of lattice points.
This is possible here since $Q^0_{\ba\bb}=Q^1_{\ba\bb}=Q^2_{\ba\bb}$.
We denote this subset as $Q_{\ba\bb}$ for the remainder of this
argument.

Let $\La'$ be the 3-dimensional lattice in $\RR^3$ generated by the vectors $(1,0,0),(0,1,0),
(\frac{1}{3}, \frac{1}{3}, \frac{1}{3}).$ As a set, $\La'$ is the union $\La_0 \cup \La_1 \cup
\La_2$, where $$\La_\ell=\La+\left(\tfrac{\ell}{3},\tfrac{\ell}{3},\tfrac{\ell}{3}\right)$$ is the
translate of the standard integer lattice $\La$ by the vector
$\left(\frac{\ell}{3},\frac{\ell}{3},\frac{\ell}{3}\right).$ Alternatively, $\La'$ is the lattice
which intersects the unit cube $[0,1]^3$ at its vertices and at the interior points
$\left(\frac{1}{3}, \frac{1}{3}, \frac{1}{3}\right), \left(\frac{2}{3}, \frac{2}{3},
\frac{2}{3}\right).$ It is evident that $\La'$ contains the standard integer lattice as a
sublattice.

Given an arbitrary lattice $\La$ in $\RR^N$,
we call a convex polytope $\cP$ a {\sl $\La$-lattice polytope} if $\cP$
has vertices on $\La;$ and we call $\cQ$
a {\sl $\La$-rational polytope} if $d\cQ$
is a $\La$-lattice polytope for some dilation by a positive integer $d$.
Let $f_\La(\cP,n) = \#\left(n \cP \cap \La\right)$ be the number of
lattice
points in the dilated region. Ehrhart's theorems translate immediately
to this setting because the entire picture can be pulled back
to the standard situation by a linear map which takes $\La$
isomorphically
to the standard lattice.

Returning to our situation of the nonstandard lattice $\La'$ in
$\RR^3,$
for a fixed $\ell \in \{0,1,2\}$,
it follows from equation (\ref{cform}) with $r=pqn+m$ that
the contribution of the component
$\cC^\ell_{\ba\bb}$ to $\tau_{SU(3)} (\Si_n)$
is given by $\#\left((pqn+m) Q^\ell_{\ba\bb} \cap \La_\ell\right)$.
Summing over $\ell$, we see that the contributions
of the components
$\bigcup_{\ell=0}^2 \cC^\ell_{\ba\bb}$ to $\tau_{SU(3)} (\Si_n)$
are given  by
$f_{\La'}(Q_{\ba\bb}, pqn+m)$.
Note that $Q_{\ba\bb}$ is
a $\La'$-rational lattice with $d=pq$, so $f_{\La'}(Q_{\ba\bb},k)$
is a quasi-polynomial of periodicity $d=pq.$ Again, since the residue
class $pqn+m$ modulo $d=pq$
is constant, we conclude that $f_{\La'}(Q_{\ba\bb},pqn+m)$ is actually
polynomial in $n$, completing the proof of the theorem in
this case.

\medskip
Assume 2 holds and choose $a_n \equiv (p+q)m \pmod{3}$
and $c_n \equiv -m \pmod{3}$  as
in Lemma \ref{chwzac}.
If $\ell =0,$ then $Q^0_{\ba\bb}$
is a rational polytope with $d=pq$ and the contribution of
$\cC^0_{\ba\bb}$ to $\tau_{SU(3)}(\Si_n)$
is given by $f_\La(Q^0_{\ba\bb}, pqn+m).$ Since $f_\La(Q^0_{\ba\bb}, k)$
is a quasi-polynomial of periodicity $d=pq,$ and
since the residue class of $pqn+m$ modulo $pq$ is constant,
it follows that the contribution of $\cC^0_{\ba\bb}$ to
$\tau_{SU(3)}(\Si_n)$ is a quadratic polynomial in $n$.

If $\ell = 1$ or $2$, then $Q^\ell_{\ba\bb}$ is a rational polytope with $d=3pq$, but that is not
sufficient for our needs. Notice from equations (\ref{TypeIvert}) and (\ref{TypeIIvert}) that the
dilation $pq Q^\ell_{\ba\bb}$ has vertices on the translate
$$\La_\ep=\La+\left(\tfrac{\ep}{3},\tfrac{\ep}{3},\tfrac{\ep}{3}\right)$$
of the standard integer
lattice $\La$, where $\ep \in \{0,1,2\}$ is given by $\ep \equiv -m \ell \pmod{3}.$ Further,the
contribution of $\cC^\ell_{\ba\bb}$ to $\tau_{SU(3)}(\Si_n)$ is given by
$\#\left((pqn+m)Q^\ell_{\ba\bb} \cap \La_{\ep} \right)$ (because $\ep \equiv -m \ell \equiv c_n
\ell \pmod{3}$). Although $\La_\ep$ is not really a lattice, we can translate the entire situation by
subtracting $\left(\tfrac{\ep}{3},\tfrac{\ep}{3},\tfrac{\ep}{3}\right)$ from $\La_\ep$ and
subtracting $\left(\tfrac{\ep}{3pq},\tfrac{\ep}{3pq},\tfrac{\ep}{3pq}\right)$ from
$Q^\ell_{\ba\bb}$. The resulting region, denoted here $\widetilde{Q}^\ell_{\ba\bb}$, is a rational
polytope with $d=pq.$ Moreover, $$f_\La(\widetilde{Q}^\ell_{\ba\bb},pqn+m)
=\#\left((pqn+m)\widetilde{Q}^\ell_{\ba\bb} \cap \La\right) =\#\left((pqn+m)Q^\ell_{\ba\bb} \cap
\La_{\ep} \right),$$ the contribution of $\cC^\ell_{\ba\bb}$ to $\tau_{SU(3)}(\Si_n)$. Now since
$f_\La(\widetilde{Q}^\ell_{\ba\bb},k)$ is a quasi-polynomial of periodicity $d=pq,$ we obtain the
desired conclusion and this completes the proof.  \hfill $\Box$

\subsection{Concluding remarks} \label{conclud}
Table \ref{table1} gives some computations of the integer valued Casson
invariant
$\tau_{SU(3)}$ for Brieskorn spheres $\Si(p,q,r)$. This extends the
computations
given in \cite{BHK}, where it was assumed that $p=2.$
\begin{center}
\begin{table}[ht]
\renewcommand\arraystretch{1.15}
\noindent\[\begin{array}{|c|c|}   \hline
    \Si   & \qquad    \tau_{SU(3)} (\Si) \qquad \\  \hline    \hline
\Si(2,3,6n\pm 1) &  3 n^2 \pm n \\ \hline
\Si(2,5,10n\pm 1) &  33 n^2 \pm 9n  \\
\Si(2,5,10n\pm 3) &  33 n^2 \pm 19n + 2 \\ \hline
\Si(2,7,14n\pm 1) &  138 n^2 \pm 26n  \\
\Si(2,7,14n\pm 3) &  138 n^2 \pm  62n + 4  \\
\Si(2,7,14n\pm 5) &  138 n^2 \pm 102n + 16  \\ \hline
\Si(2,9,18n\pm 1) &  390 n^2 \pm 58n \\
\Si(2,9,18n\pm 5) &  390 n^2 \pm 210n + 24 \\ \Si(2,9,18n\pm 7) &  390
n^2 \pm 298n +52 \\ \hline
\Si(3,4,12n\pm 1) &  105 n^2 \pm 21n \\
\Si(3,4,12n\pm 5) &  105 n^2 \pm 87n+16 \\ \hline
\Si(3,5,15n\pm 1) &  276 n^2 \pm 40 n \\
\Si(3,5,15n\pm 2) &  276 n^2 \pm 74 n + 2\\
\Si(3,5,15n\pm 4) &  276 n^2 \pm 148 n + 16\\
\;  \Si(3,5,15n\pm 7)\;  &  \;  276 n^2 \pm 254 n + 56  \;    \\
\hline
\end{array}  \]
\medskip
\caption{{\sc Calculations of the integer valued SU(3) Casson invariant
    for some Brieskorn spheres $\Si(p,q,r)$.}}
\label{table1}
\end{table}
\end{center}
\renewcommand\arraystretch{1}
Let $K_{p,q}$ be the $(p,q)$ torus knot and set
    $X_n=1/n$ Dehn surgery on $K_{p,q}$.
Then $X_n = \pm \Si(p,q,r)$ for $r = |pqn -1|$.   Table \ref{table2}
gives the
value of $\tau_{SU(3)}(X_n)$ for various $p,q$.These computations were
performed using MAPLE.

For surgeries on torus knots, Theorem \ref{torusconj} asserts that
$$\tau_{SU(3)} (X_n) = A(K_{p,q}) n^2 - B(K_{p,q}) n, $$
where $A(K_{p,q})$ and $B(K_{p,q})$ depend only on $K_{p,q}$.
There is a  pattern for the leading coefficient  $A(K_{p,q})$
present in Table \ref{table2}. If $\De_K(z) =  \sum_{i \geq 0}
c_{2i}(K) z^{2i}$
denotes the Conway polynomial of $K$,
we conjecture generally that $\tau_{SU(3)}(X_n)$ has quadratic growth
in $n$ with leading coefficient
\begin{equation}\label{A2formula}
A(K) = 6c_4(K) + 3c_2(K)^2.
\end{equation}
This is what one would expect from  Frohman's work \cite{frohman} on
$SU(n)$
Casson knot invariants in the case of $n=3,$
at least for fibered knots (cf.~\cite{frohman-nicas, boden-nicas}).
It gives the formula
$$A(K_{p,q}) = {\displaystyle \frac{(p^2-1)(q^2-1)(2p^2
q^2-3p^2-3q^2-3)}{240}},$$
which  agrees with the data in Table \ref{table2}.
The coefficient $B(K)$ of the linear term
is not as well-behaved.
For example, interpolating the data from  Table \ref{table2}, we get
the formulas
$$B(K_{2,q})=\begin{cases}
  \frac{1}{12} (q^3-4q+3) & \text{if $q \equiv 1 \pmod 4$,} \\
\frac{1}{12} (q^3-4q-3)  & \text{if $q \equiv 3 \pmod 4$,}
\end{cases}$$

$$B(K_{3,q})=\begin{cases}
  \frac{1}{54} (20q^3+3q^2-48q+25) & \text{if $q \equiv 1 \pmod 6$,} \\
  \frac{1}{54} (20q^3-3q^2-48q+2) & \text{if $q \equiv 2 \pmod 6$,}  \\
  \frac{1}{54} (20q^3+3q^2-48q-2) & \text{if $q \equiv 4 \pmod 6$,}  \\
  \frac{1}{54} (20q^3-3q^2-48q-25)& \text{if $q \equiv 5 \pmod 6$,}
\end{cases}$$
and
$$B(K_{4,q})=\begin{cases}
  \frac{1}{16} (16q^3+q^2-42q + 25) & \text{if $q \equiv 1 \pmod 8$,} \\
  \frac{1}{16}(16q^3-q^2-42q + 39)& \text{if $q \equiv 3\pmod 8$,}  \\
  \frac{1}{16}(16q^3+q^2-42q - 39) & \text{if $q \equiv 5 \pmod 8$,}  \\
  \frac{1}{16}(16q^3-q^2-42q - 25)& \text{if $q \equiv 7 \pmod 8$.}
\end{cases}$$

\begin{center}
\begin{table}[t]
\renewcommand\arraystretch{1.15}
\noindent\[\begin{array}{|c|c||c|c||c|c|}     \hline p=2 &
\tau_{SU(3)}(X_n) & p=3 & \tau_{SU(3)}(X_n) &p=4 & \tau_{SU(3)}(X_n) \\
\hline
K_{2,3} & 3n^2-n & K_{3,4} &   105n^2-21n&K_{4,5} &   1011n^2-111n \\
K_{2,5} & 33n^2-9n& K_{3,5} &   276n^2-40n& K_{4,7} &   4110n^2-320n\\
K_{2,7} & 138n^2-26n& K_{3,7} &   1128n^2-124n&K_{4,9} &
11\,490n^2-712n \\
K_{2,9} & 390n^2-58n& K_{3,8} &   1953n^2-179n & K_{4,11}    &
25\,935n^2-1297n\\
K_{2,11} & 885n^2-107n& K_{3,10}    &   4851n^2-367n &K_{4,13}    &
50\,925n^2-2171n \\
K_{2,13}    &   1743n^2-179n&K_{3,11}    &   7140n^2-476n& K_{4,15}
&   90\,636n^2-3320n \\
K_{2,15}    &   3108n^2-276n& K_{3,13}    &   14\,028n^2-812n& K_{4,17}
    &   149\,940n^2-4888n\\
K_{2,17}    &   5148n^2-404n& K_{3,14}    &   18\,915n^2-993n& K_{4,19}
    &   234\,405n^2-6789n\\  K_{2,19}    &   8055n^2-565n& K_{3,16}    &
   32\,385n^2-1517n& K_{4,21}    &  350\,295n^2-9231n \\
K_{2,21}    &   12\,045n^2-765n& K_{3,17}    &   41\,328n^2-1788n&
K_{4,23}    & 504\,570n^2-12\,072n  \\
K_{2,23}    &   17\,358n^2-1006n& K_{3,19}  & 64\,620 n^2 - 2544 n  &
K_{4,25}    & 704\,886 n^2 - 15\,600 n \\
K_{2,25}    &   24\,258n^2-1294n& K_{3,20}  & 79\,401 n^2 - 2923 n  &
K_{4,27}    & 959\,595 n^2 - 19\,569 n \\
\,K_{2,27}\,    &   \,33\,033n^2-1631n\,& \,K_{3,22} \, & \,116\,403
n^2 - 3951 n  \, & \, K_{4,29}   \, &\, 1\,277\,745 n^2 - 24\,363 n\,
\\ \hline
\end{array}  \]
\medskip
\caption{{\sc The integer valued SU(3) Casson invariant for homology
3-spheres $X_n$
obtained by  $1/n$ surgery on  $K_{p,q}$}}
\label{table2}
\end{table}
\end{center}\renewcommand\arraystretch{1}

  The increasing complexity of these formulas makes it difficult to guess
a general formula for $B(K)$ in terms of classical invariants of the
knot.
Nevertheless, it provides a negative answer to
the question of whether $\tau_{SU(3)}$ is a finite type invariant.
For suppose $\tau_{SU(3)}$ were a finite type invariant.
Then, as explained to us by Stavros Garoufalidis,
$B(K_{p,q})$ would necessarily be a polynomial in $p$ and $q$.
Since $B(K_{p,q})$ is obviously not a polynomial in $p$ and $q$, it
follows that
$\tau_{SU(3)}$ is not a finite type invariant of any order.

Notice that $\tau_{SU(3)}(X)$ is even in all known computations.
Further, a simple argument using the involution
on $\cM_{SU(3)}$ induced by complex conjugation proves
evenness of  $\tau_{SU(3)}(X)$
under the hypothesis that $H^1_\al(X;su(3))=0$ for every nontrivial
representation
$\al\colon\pi_1 X \to SU(3)$.
We conjecture  that $\tau_{SU(3)}(X)$ is even for all homology
3-spheres.

\vfill\eject

\end{document}